\documentclass{article}

\usepackage{arxiv}

\usepackage[utf8]{inputenc} 
\usepackage[T1]{fontenc}    
\usepackage{hyperref}       
\usepackage{url}            
\usepackage{booktabs}       
\usepackage{amsfonts}       
\usepackage{nicefrac}       
\usepackage{microtype}      
\usepackage{lipsum}
\usepackage{graphicx}
\graphicspath{ {./images/} }

\usepackage[T1]{fontenc}
\usepackage[utf8]{inputenc}
\usepackage[british]{babel}

\usepackage[monochrome]{xcolor}
\usepackage{tabularx}
\usepackage{multirow}
\usepackage{subcaption}
\usepackage{enumitem}
\usepackage{cite}
\usepackage{float}

\usepackage{authblk} 
\usepackage{amsmath} 
\usepackage{pgffor}
\usepackage{alphalph}

\title{An FEA surrogate model with Boundary Oriented Graph Embedding approach}

\author[a]{Xingyu Fu}
\author[a]{Fengfeng Zhou}
\author[b]{Dheeraj Peddireddy}
\author[c]{Zhengyang Kang}

\author[a,d]{Martin Byung-Guk Jun* \thanks{Corresponding author: Martin Byung-Guk Jun, mbgjun@purdue.edu}}
\author[b]{Vaneet Aggarwal}

\affil[a]{School of Mechanical Engineering, Purdue University, 585 Purdue Mall, West Lafayette, IN 47907 , USA}
\affil[b]{School of Industrial Engineering, Purdue University, 585 Purdue Mall, West Lafayette, IN 47907 , USA}
\affil[c]{School of Mechanical and Power Engineering, Nanjing Tech University, Nanjing, 211800, Jiangsu, China}
\affil[d]{Indiana Manufacturing Competitiveness Center (IN-MaC), West Lafayette, IN, 47907, USA}

\begin{document}

\maketitle

\begin{abstract}
In this work, we present a Boundary Oriented Graph Embedding (BOGE) approach for the Graph Neural Network (GNN) to serve as a general surrogate model for regressing physical fields and solving boundary value problems. Providing shortcuts for both boundary elements and local neighbor elements, the BOGE approach can embed structured mesh elements into the graph and performs an efficient regression on large-scale triangular-mesh-based FEA results, which cannot be realized by other machine-learning-based surrogate methods. Focusing on the cantilever beam problem, our BOGE approach cannot only fit the distribution of stress fields but also regresses the topological optimization results, which show its potential of realizing abstract decision-making design process. The BOGE approach with 3-layer DeepGCN model \textcolor{blue}{achieves the regression with MSE of 0.011706 (2.41\% MAPE) for stress field prediction and 0.002735 MSE (with 1.58\% elements having error larger than 0.01) for topological optimization.} The overall concept of the BOGE approach paves the way for a general and efficient deep-learning-based FEA simulator that will benefit both industry and design-related areas.
\end{abstract}


\section{Introduction}\label{intro}

Finite Element Analysis (FEA), which plays a critical role in modern engineering design and quality evaluations, provides numerical solutions for boundary value problems and simulates the distribution of physical fields with high precision \cite{axelsson2001finite}. However, FEA generally requires high-performance computational resources and has an expensive time cost which becomes a bottleneck in the present manufacturing industry, especially for advanced time-sensitive techniques: inverse modeling, agile manufacturing, generative design system, etc. In order to solve this problem, many researchers have employed machine learning technologies to regress the distribution of physical fields to obtain instant simulation results while compromising the accuracy. Some researchers employed conventional machine learning methods (Support Vector Regression-SVG \cite{capuano2019smart}, Artificial Neural Network - ANN \cite{tamaddon2020data}) to approximate the physical field distribution and solve the boundary value problem. However, most methods are for specific boundary conditions and require extra training ahead of the application. This limits their capability of serving as a general FEA surrogate model, and their ability to solve abstract decision-making problems.

\begin{figure}[htp]
\vspace{-0.12in}
\centering
	\includegraphics[keepaspectratio,width=0.7\textwidth]{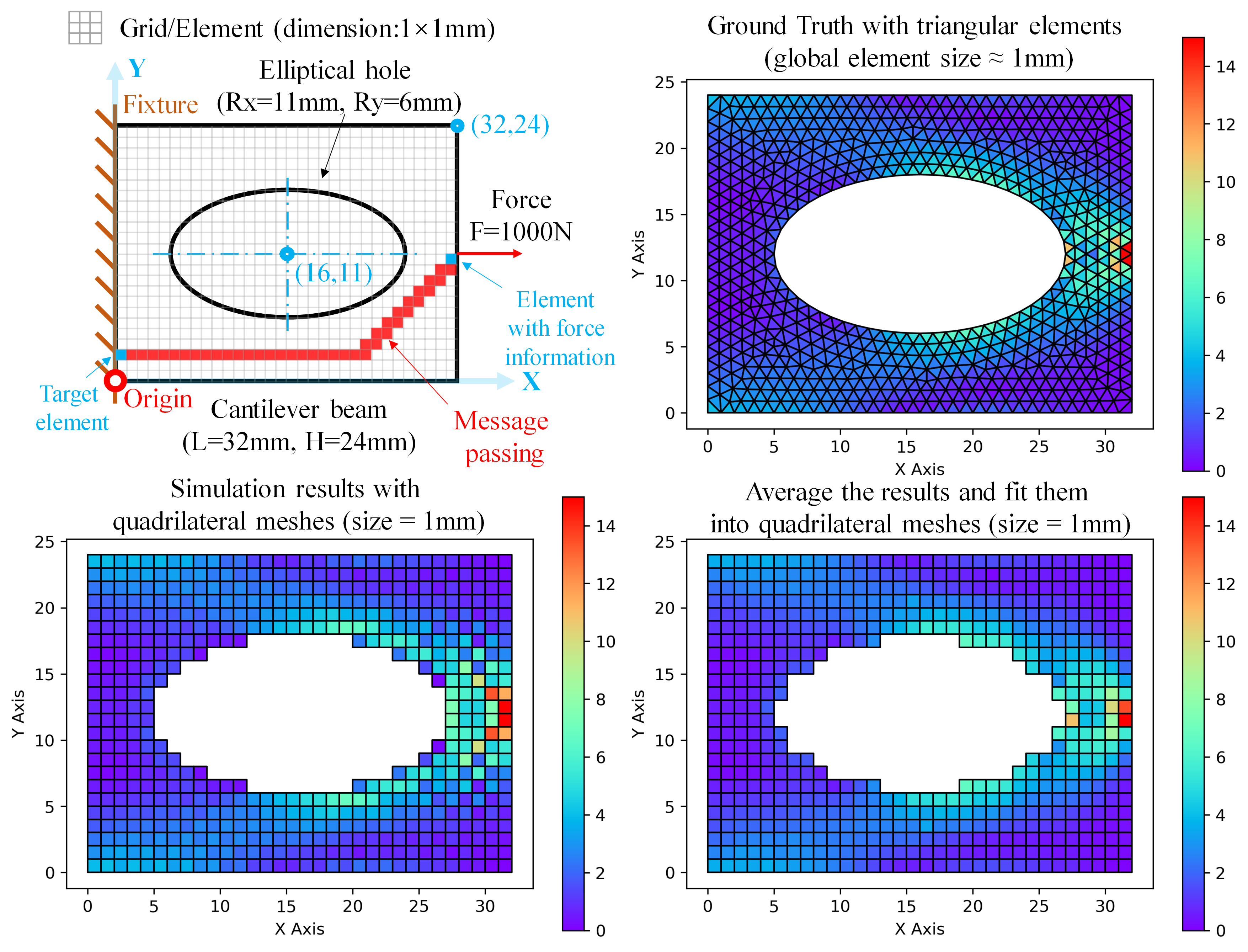}
	\caption{Comparison between different FEM simulation and prediction results}
	\label{fig:0_1_FemCompare}
	\vspace{-0.12in}
\end{figure}

Recently, with the rapid development of Deep Learning (DL) techniques, some researchers employed Convolutional Neural Network (CNN) as the backbone for solving physical field prediction problems, and found out that CNN based surrogate models cannot only solve basic boundary value problems (e.g., solid mechanics \cite{nie2020stress, kantzos2019design, khadilkar2019deep}, fluid dynamics \cite{kutz2017deep, guo2016convolutional, zhang2018application}, etc.) but also perform efficiently on some abstract design problems based on predicted physical field, for instance, topology optimization \cite{lee2020cnn, nie2021topologygan,zhang2019deep,banga20183d}. However, CNN-based surrogate models still have inevitable disadvantages. Firstly, CNN can only accept input tensor with fixed input size, which becomes an obstacle for input tensor encoding, especially for some complex FEA conditions. Also, since CNN can only perform calculations on Euclidean coordinates, most CNN-based surrogate models employ square-shaped or cube-shaped meshes (2D quadrilateral or 3D hexahedron meshes with uniform side lengths) that require homogeneous grid generation on the input FEA model. The grid resolution sometimes cannot reach a sufficiently high level since the CNN usually cannot handle a large tensor input with hundreds of meshes along X/Y/Z axes due to the limitation of modern computation resources. Lastly, most FEA run on structured meshes (e.g., triangular meshes, or polygon meshes with different side lengths) that cannot directly fit into the CNN input tensor. Some researchers have presented alternate solutions to handle this problem, which is illustrated in Fig.~\ref{fig:0_1_FemCompare}. We provide a von Mises stress distribution prediction for a 2D cantilever beam based on simulation results from ABAQUS \cite{abaqus2011abaqus}. The beam is made of A36 steel with Young’s modulus 200 GPa and Poisson’s ratio of 0.32 and has an elliptical hole inside. Its left edge is fixed and an external force 1000N is exerted on the center of its right edge. The ground truth simulation result with triangular meshes (with an approximate global mesh size of 1.0 mm) is shown in Fig.~\ref{fig:0_1_FemCompare}. One type of method for using CNN to regress the physical field is to approximate the original input model's geometry with a high-resolution voxelized representation, shown in "Simulation results with quadrilateral meshes" in Fig.~\ref{fig:0_1_FemCompare} \cite{nie2020stress}. However, the resolution is still limited by computational resources, and the simulation result can significantly change since the change in model shape can largely influence the global boundary condition, shown in Fig.~\ref{fig:0_1_FemCompare}. The other type of method usually uses some post-process tasks, for instance, averaging the physical fields to fit into the square/cube-shaped meshes. This will still lead to extra approximation error for the output results, shown in "Average the results" in Fig.~\ref{fig:0_1_FemCompare}. All of those defects call for a more efficient DL surrogate model for FEA. The Graph Neural Network (GNN) comes into sight.

The GNN is a generalized version of CNN for non-Euclidean graphs, which has been widely applied to abstract decision-making and data-regression problems on different graph data structures, for instance: point cloud, triangular meshes, and molecular structures \cite{alet2019graph, wang2020kalibre, belbute2020combining, guo2020semi, ogoke2020graph,  pfaff2020learning, sanchez2020learning, gilmer2017neural}. Since the structured meshes can be naturally embedded by an undirected graph $G$, many researchers investigated GNN-based surrogate models to solving FEA regression problems, show in Table~\ref{tab:0_1_Compare}. It is assumed that an undirected graph $G=(V, E)$ represents the embedded structured meshes in an FEA model. Then, generally, $V={1, \cdots, N}$ can be the set of N graph vertices (or graph nodes) denoting N geometrical node points in FEA while $ E \subseteq V \times V$ is the set of graph links (or graph edges) representing the connection of neighborhood geometrical nodes. Graph vertex $v$ contains features $f_v$ encoding the properties of a geometrical node point (e.g., coordinates, material properties, boundary conditions). Graph links generally employ adjacent matrix $E \in [0,1]^{N \times N}$ to describe the connection between geometrical node points. Nowadays, most applications employed GNN layers with Message Passing Neural Network (MPNN) framework \cite{gilmer2017neural} whose updated message $m_v^{t+1}$ of graph vertex $v$ at layer $t+1$ can be represented by:

\begin{equation}
m_v^{t+1} = \sum_{\omega \in N(v)}^{} M_t (h_v^{t}, h_\omega^{t}, e_{v\omega})
\label{mpnn1}
\end{equation}
\begin{equation}
h_v^{t+1}=U_t(h_v^{t}, m_v^{t+1})
\label{mpnn2}
\end{equation}
\begin{equation}
\hat{p}=R( \left \{ h_v^{T}|v \in G \right \})
\label{mpnn3}
\end{equation} 
\newline Where $M_t$, $U_t$ and $R$ are message function, vertex update function and readout function, respectively; $h$ represents the hidden state in GNN layers (e.g., $h_v^{t}$ is the hidden state of $v$ at layer $t$); $N(v)$ denotes the neighbors of $v$ in graph $G$; $e_{v\omega}$ represents the graph link features; $\hat{p}$ is the final prediction results. From the Eq.~(\ref{mpnn2}), it can be seen that GNN layers can only pass the information of the embedded feature to its neighbor vertices. If the distance on a graph for passing an important message is relatively long, with current graph embedding approaches, more layers of GNN are needed. For instance, this can happen when GNN is employed to pass the information from the element with boundary condition information (external force or fixture, e.g., element marked by "Element with force information" in Fig.~\ref{fig:0_1_FemCompare}) to the element whose stress field is required to be solved ("Target element" in Fig.~\ref{fig:0_1_FemCompare}). One of the shortest message propagation paths from the boundary element to the target element is marked with red shown in Fig.~\ref{fig:0_1_FemCompare}, which requires several layers of GNN to pass the information. However, deeply stacking the GNN layers can suffer from over-smoothing problems which cause indistinguishable graph-vertex embeddings and make GNN eventually not predict anything \cite{li2020deepergcn, yang2020revisiting}. In general, this problem restricts most GNN models to around 3 layers with only short-range graph-vertex interactions and limited message passing distance.


\begin{table}[htp]
\centering
\resizebox{0.7\textwidth}{!} {
\begin{tabular}{cccccc}
\hline
Model type                                                & Graph size & \begin{tabular}[c]{@{}c@{}}Graph vertices\\interaction scale\end{tabular} & Generalization & References           \\\hline
Conventional methods                                      & -          & -                                      & Limited        & \cite{capuano2019smart, tamaddon2020data}  \\
\begin{tabular}[c]{@{}c@{}}CNN with quadrilateral\\ or hexahedron meshes\end{tabular}& -          & -                                      & Good           & \cite{nie2020stress,nie2021topologygan,koeppe2020intelligent,khadilkar2019deep, kutz2017deep, guo2016convolutional,zhang2018application}     &  \\ 
\begin{tabular}[c]{@{}c@{}}GNN for small graph\\ with triangular meshes\end{tabular} &Small    & Short range                          & Good           & \cite{alet2019graph, wang2020kalibre}   \\ 
\begin{tabular}[c]{@{}c@{}}GNN for large graph\\ with triangular meshes\end{tabular} &Large    & Long range                            & Limited        & \cite{belbute2020combining,guo2020semi}       \\ 
Other GNN models                                          &Large    & Short range                          & Good           & \cite{ogoke2020graph, pfaff2020learning,sanchez2020learning}     \\  \hline          
\end{tabular}}
\caption{Researches for FEA surrogate models}
\label{tab:0_1_Compare}
\end{table}


Given the limitations of GNN, from Table~\ref{tab:0_1_Compare}, it can be seen that the GNN is effective on two types of physical field regression problems: The first type is the problems which input a small-scale graph data into GNN with generally no more than 100 graph vertices, and usually requires data pre-processing or model simplification \cite{alet2019graph, wang2020kalibre}; the other type can solve the large-scale-graph problem, however, has only short-range graph-vertex interactions. This means the graph distance for the propagation of boundary information is relatively short since the elements embedding boundary information is relatively close to the target element whose physical field is required to be solved. For instance, Pfaff et al. \cite{pfaff2020learning} developed an efficient GNN model for dynamic simulations based on triangular meshes. However, since the dataset contained sequential physical fields (displacement field, velocity field, stress field, etc.) and the input physical field for the current frame (or simulation step) was derived from the previous frame, the boundary condition information can slowly propagate from the boundary elements to the elements inside the model frame by frame. Then, the physical field of the target element was only determined by its surrounding elements that embed boundary information, which only requires the short-range graph-vertex interaction. The same happens when the regression model works on sequential physical field predictions(e.g., \cite{sanchez2020learning}). Ogoke et al. \cite{ogoke2020graph} predicted drag force on a large graph using velocity fields , but only elements around the streamlined object contributed to the total force. Although some researchers successfully developed efficient GNN based prediction model on large graphs, their applications were still limited. Belbute-Peres et al. \cite{ belbute2020combining} developed a GNN based fluid dynamics prediction model on large graphs by interpolating PDE solver’s results, which had a relatively short-range graph distance for the propagation of boundary-condition information since the result from the PDE solver carried the model's boundary condition. Guo et al. \cite{ guo2020semi} developed a semi-supervised GNN model to predict optimized architected material structure on a large graph, but the selected training dataset was chosen from inside of the material which embedded boundary conditions. Also, the training needed to be conducted before each application in \cite{ guo2020semi}.

From above, it can be seen that there’s no ready-to-use generalized machine learning approach for solving structured-mesh-based boundary value problems, especially when long-range graph-vertex interactions are involved. FEA generally requires structured meshes to accurately approximate the boundary of the input model geometry. It usually provides finer meshes to the model's delicate structures to improve the computation accuracy and employ coarse meshes to save the computational resource. The CNN-based FEA surrogate models have difficulty approximating the structured-mesh-based results, which becomes a bottleneck for the efficiency of related algorithms. The GNN-based approach shows its unique advantage of regressing structured-mesh-based results, however, a precise FEA model usually requires fine meshing to provide a high-resolution result, which can largely increase the graph size and therefore extends the propagation distance of the boundary information. This requires a deep GNN model and usually can cause the oversmoothing problem that is difficult to solve. Also, when GNN is applied to solve the static boundary value problems (different from the dynamic problems using the sequential physical field as its input with multiple frames\cite{pfaff2020learning, sanchez2020learning}), the target element needs to gather all the boundary condition information in one frame (or simulation step) and then provides the prediction result. For instance, in the cantilever beam problem, the target element stress field is determined by all the elements with the boundary information (fixture, model's geometrical edges, external forces, shown in Fig.~\ref{fig:0_1_FemCompare}), and this cannot be solved by any of the current approaches. Therefore, in order to solve these problems and provide a potential general DL-based surrogate regression model, we focus on the most fundamental solid mechanics problem - cantilever beam problem and develop a Boundary Oriented Graph Embedding (BOGE) approach. With the novel BOGE method, the GNN model works effectively on high-resolution triangular-mesh-based FEA simulation and predicts the accurate distribution of von Miss stress. The topological optimization predicted by GNN shows the capability of  our model to conduct abstract decision-making results based on the model’s boundary condition. In this paper, Section \ref{data} introduces our novel data representation which provides both multi-hop local adjacent information and shortcuts for global boundary information. Section \ref{dataset} describes the details of our dataset preparation and the parametric designs for the FEA simulation results. Section \ref{train} shows the model’s training results that \textcolor{blue}{achieves the regression with MSE of 0.011706 (2.41\% MAPE) for stress field prediction and 0.002735 MSE (with 1.58\% elements having error larger than 0.01) for topological optimization.} The outperformance of our model's prediction accuracy as well as its irreplaceable properties on regression of FEA results shows our models’ un-substitutable significance on fast FEA simulations. We summarize the advantages of our GNN based surrogate model as follows:
\begin{itemize}[leftmargin=*]
\item We propose a new BOGE (Boundary Oriented Graph Embedding) approach for regression of physical fields based on structured meshes (especially those based on triangular meshes) resulting from boundary value problems. The BOGE approach not only provides shortcuts for neighborhood meshes to smooth prediction results, but also offers shortcuts of global boundary conditions to each internal element, which circumvents the message passing limitation of the Message Passing Neural Network (MPNN) framework. Also, the BOGE largely compresses the graph size derived from the simulation model, enabling one graph vertex embed both mesh properties and geometrical edge properties, which leads to more efficient message propagation in GNN.
\item We employ the parametric design method to prepare and standardize the synthesized dataset for cantilever beam simulation. The parametric dataset can be expanded with more design models to enhance the generalisation capacity of our model. With the same design approach, the dataset for other types of physical-field simulators (for instance, temperature, velocity, etc.) can be developed.
\item The BOGE is compatible with most GNN based MPNN frameworks. With the simple combination of BOGE and other state of art GNN models, our GNN based surrogate model shows its outstanding performance of approximating both physical fields and highly abstract optimization results. The model has the potential to be applied to most physical-field fitting problems with any type of polygon meshes or mixed types of meshes. The GNN based surrogate approach paves the way for further deep learning-based physical field simulators and will benefit all computational simulation fields.

\end{itemize}

\section{Data representation}\label{data}

The physical field distribution is generally derived by resolving partial differential equations under the given boundary condition. For the stress analysis of the cantilever beam problem, denote the cantilever beam's geometry as $\Omega$. Then, the global boundary condition $\partial \Omega$ of the problem can be expressed by:

\begin{equation}
\partial \Omega = S_u + S_p
\label{bc}
\end{equation}
\newline Where $S_u$ and $S_p$ respectively represent the displacement (essential) and force/stress (natural) boundary condition. In an FEA model, the algorithm usually employs discretized approach to decompose the input geometry $\Omega$ into finite elements so that the complex boundary value problem can be simplified into solving the partial differential equation for each element. For the linear elasticity problem on isotropic materials, the canonical form of the relationship between mesh nodal displacement $\mathbf{U}$ and the global boundary condition $\mathbf{F}$ can be expressed by:

\begin{equation}
\mathbf{K}\mathbf{U}=\mathbf{F}
\label{hooke}
\end{equation}
\newline Where $\mathbf{K}$ is the global stiffness matrix. Here, the global stiffness matrix $\mathbf{K}$ is assembled by the elemental stiffness matrices $\mathbf{k_e}$, which can be represented by an integration function of the element area $A$:

\begin{equation}
\mathbf{k_e} = \int_{A}^{}\mathbf{B_e^T}\mathbf{C_e}\mathbf{B_e} \textup{d}A
\label{ke}
\end{equation}
\newline Where $\mathbf{B_e}$ is the displacement differentiation matrix obtained by differentiation of displacements expressed through shape functions and nodal displacements; $\mathbf{C_e}$ is the elasticity matrix constructed by Young’s modulus ($E$) and Poisson’s ratio ($\nu$). Then, the nodal displacement field $\mathbf{U}$ can be solved. Based on $\mathbf{U}$, the stress field $\mathbf{\sigma} = [\sigma_x, \sigma_y, \tau_{xy}, \tau_{yx}]^T$ can be expressed by:

\begin{equation}
\mathbf{\sigma} = \mathbf{C_g}\mathbf{B}\mathbf{u}
\label{stress}
\end{equation}
\newline Where $\mathbf{B}$ is the assembly of $\mathbf{B_e}$ matrices; the block-diagonal matrix $\mathbf{C_g}$ is constructed from elemental elasticity tensor $\mathbf{C}$ that can be derived by $\mathbf{C_e}$ \cite{nie2020stress}. From the equation, it can be seen that the element nodal stress state and the nodal displacement are not only determined by the global boundary condition $\partial \Omega$ but also influenced by local boundaries determined by its neighbor elements due to the assembly of $\mathbf{K}$. Specifically, if a GNN model is applied to solve the cantilever beam problem, the target element requires both long-range global boundary condition information and the short-range neighbor elements' properties to generate its stress field. Also, in a 2D cantilever-beam FEA simulation, three key items need to be fully defined to determine the final simulation solution: global boundary condition ($\partial \Omega$), predefined elements for the model geometry ($\Omega$), and material properties. This section mainly introduces the BOGE approach which can comprehensively feed that information from the FEA model to the GNN.

\subsection{Graph vertex embedding method}\label{chs}

Most GNN-based FEA surrogate models embed mesh nodes as graph vertices for physical field simulation \cite{belbute2020combining, guo2020semi, ogoke2020graph, pfaff2020learning, sanchez2020learning}. However, boundary conditions can be cast on mesh edges or mesh bodies, which is not convenient to be represented using element nodes. Also, it’s natural to encode the material property inside the mesh body rather than mesh nodes, especially for the conditions when material property changes spatially, for instance, composite materials and non-homogeneous material. The mesh nodes are between the mesh bodies, which is hard to embed the change of mesh body properties. In order to solve this problem, according to the convention from the previous research \cite{pfaff2020learning}, the mesh edge, as well as the mesh body and the mesh node, are considered as graph vertices. However, compared with the condition when only element nodes are considered, embedding each mesh node, mesh edge, and mesh body as separate graph vertices will expand the graph size 3 times for triangular-mesh-based simulations, which can run out of the computational resources or lead to an extremely long message passing distance for the propagation of boundary conditions. Therefore, we combine the properties of the mesh node, mesh edge, and mesh body into one graph vertex and embed the overall mesh features into the graph vertex to provide a compressed graph representation.

\begin{figure}[htp]
\vspace{-0.12in}
\centering
	\includegraphics[keepaspectratio,width=0.7\textwidth]{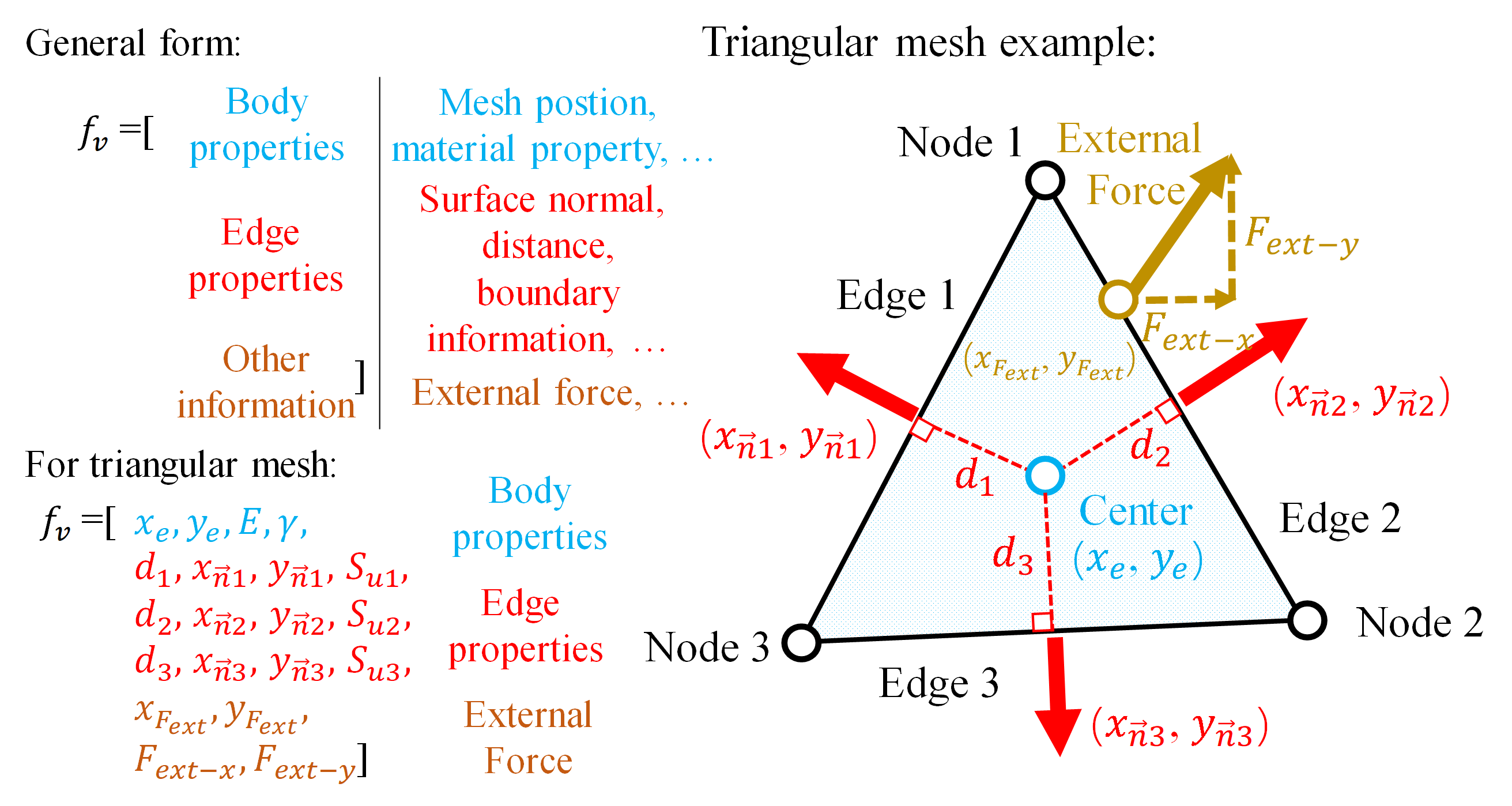}
	\caption{Encoding approach for triangular mesh}
	\label{fig:1_1_Encoding}
	\vspace{-0.12in}
\end{figure}

The format of our graph embedding features $f_v$ is illustrated in Fig.~\ref{fig:1_1_Encoding}. The first several channels encode the center of the mesh which provides the approximate position information of the mesh body. Material properties follow those channels which provide mesh body properties. After that, we encode the geometrical properties of mesh edges using the edge's normal vector and its distance to the mesh center to fully define the edge's geometrical position. Each edge's geometrical property follows the edge's physical property that encodes the global boundary condition of the simulation model. Following the edge properties, additional point-based information can be added to the data representation, for instance, the external force cast on the mesh.

In Fig.~\ref{fig:1_1_Encoding}, we provide an example of our data representation for solving the cantilever beam problem (details of the problem definition are explained in Sec.~\ref{dataset}). For a 2D cantilever beam problem based on triangular-mesh simulation, the first two channels encode the mesh center $(x_e, y_e)$,  which is calculated by the average of the coordinates for all the mesh nodes. The following body property for a stress analysis simulation can be Young’s modulus ($E$) and Poisson’s ratio ($\nu$) of the mesh material. After that, three mesh edges with their distances from the center ($d_k$) and their normal vectors ($(x_{\vec{n}k}, y_{\vec{n}k})$) are encoded (where k=1, 2, 3 representing each edge of the mesh). Following each edge geometry, the boundary condition is encoded by its boundary state ($S_{uk}$). We set the boundary state as 1 when the edge is on the outer contour of the input model geometry, while $S_{uk}=2$ for edges on fixtures. After that, the data representation encodes the external 2D force ($S_p$) by attaching its relative position ($(x_{F_{ext}}, y_{F_{ext}})$, calculated according to the mesh center) and its magnitude ($(F_{ext-x}, F_{ext-y})$) to the graph vertex features. Meshes without the external force need to be padded with zero to meet the requirement of graph embedding. Then, all the necessary information from a triangular mesh has been embedded into the graph vertex.

\subsection{Adjacency information}\label{adj}

The GNN requires the input of graph links $E$ to figure out the message passing direction, which is encoded by the adjacency matrix. Conventionally, the mesh-based graph embedding is bounded by the geometry of the input mesh model  ($\Omega$), which only considers the graph links between the mesh and its neighboring meshes. However, in order to shorten the message passing distance to use shallow layers of GNN to regress the physical fields, we provide shortcuts for both the global boundary condition and the local neighboring meshes.

\begin{figure}[htp]
\vspace{-0.12in}
\centering
	\includegraphics[keepaspectratio,width=0.7\textwidth]{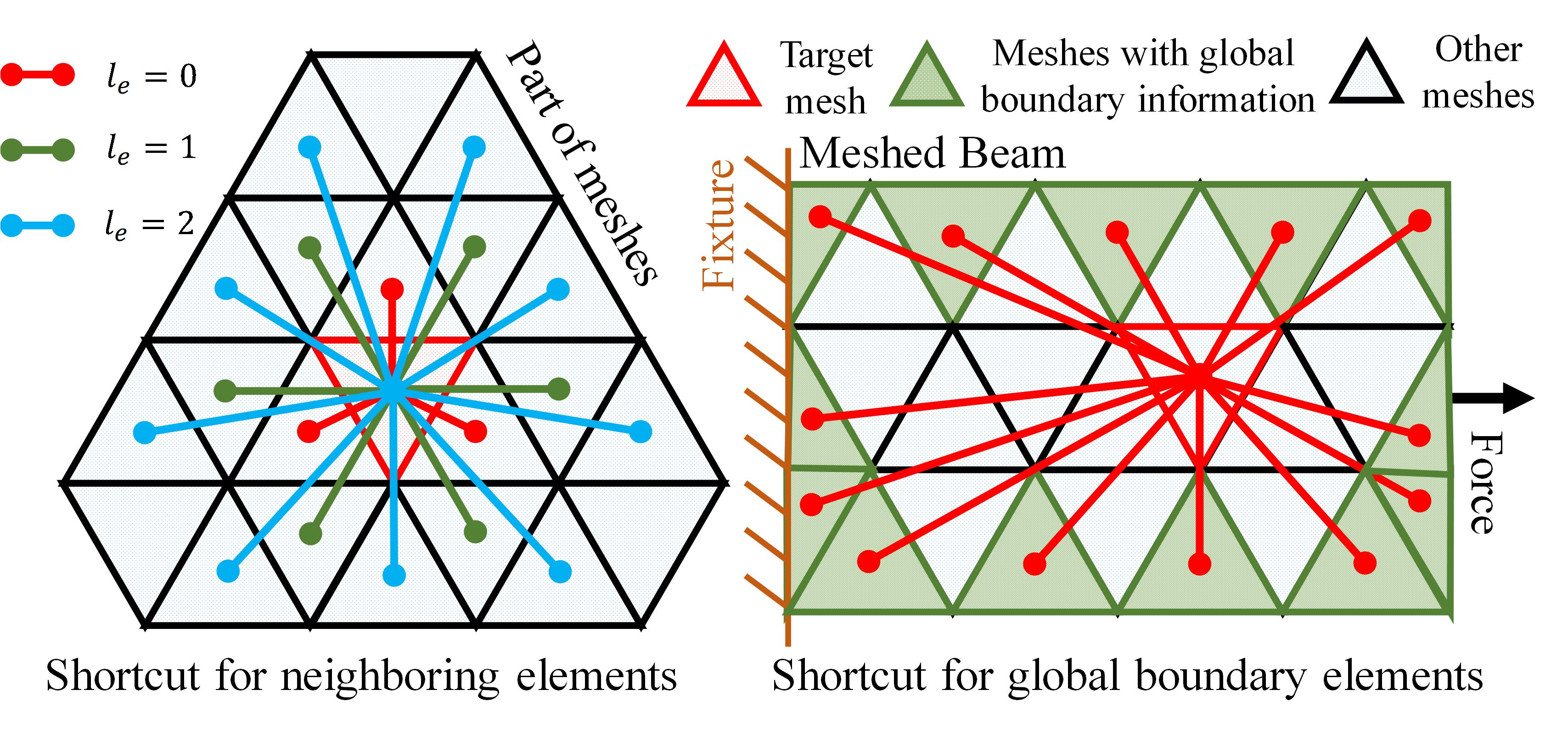}
	\caption{Encoding for graph connectivity and shortcuts}
	\label{fig:1_2_Shortcut}
	\vspace{-0.12in}
\end{figure}

Most researchers focus on one-hop graph vertices, which are the meshes straightly connecting to the target mesh. Denote the distance between the target element and its neighbor element on the graph as $l_e$, then we represent the graph links between one-hop elements and the target element as $l_e=0$, shown in Fig.~\ref{fig:1_2_Shortcut} with red graph links. However, in order to provide more local information to the target mesh and smooth the physical-field results, we also consider two-hop ($l_e=1$) and three-hop ($l_e=2$) graph links (marked with green and blue graph links respectively in Fig.~\ref{fig:1_2_Shortcut}) and encode those to the adjacency matrix. Adding information about meshes' local connectivity assists to provide more details of local boundary conditions, whose performance is presented in Sec.~\ref{train}.

Then, we provide mesh connectivity between the target mesh and all the meshes that have the information of global boundary conditions ($\partial \Omega$) in order to provide the shortcut for solving the large-scale graph problem. In the cantilever beam problem, all the meshes with information of the model's geometrical edges, fixtures, and external forces, are connected to every internal mesh inside the beam, which is illustrated in Fig.~\ref{fig:1_2_Shortcut}. With these graph linking shortcuts, the GNN serves for collecting weighted boundary conditions from the entire simulation model, other than only passing the mesh information to its neighborhood meshes with short-range distance. Therefore, the shallow-layer GNN structure can still regress high-resolution structured-mesh-based simulation results. The performance of this improvement is shown in detail in Sec.~\ref{train}.

In this section, we present the BOGE approach to embed triangular-mesh-based FEA problems into the graph structure. This concept can be generalized to other polygon-mesh simulations with more mesh-edge information. Padding the channels with zeros to keep the uniform graph-vertex feature size for all the meshes, data representation can encode simulation with mixed type of structured mesh (e.g., triangular meshes and quadrilateral meshes) and will incur only a small amount of extra computational resources since GNN generally employs sparse tensor production during training and computing. Our BOGE approach largely compresses the graph size and provides shortcuts for the boundary conditions, which can improve the long-range graph-vertex interactions in boundary value problems. It is compatible with most MPNN-based GNN structures and can be generalized to solve other boundary value problems.
\section{Dataset preparation}\label{dataset}

In order to validate the efficiency of this approach, we focus on two aspects of the cantilever beam problem: stress field prediction and topology optimization. We summarized all the cantilever beam models provided in \cite{nie2020stress} and created a parameterized synthetic FEA simulation dataset. According to \cite{zhang2018featurenet, peddireddy2021identifying}, we provide uniform random values to provide design parameters for cantilever beams to improve the repeatability and portability of the dataset. Details on dataset preparation are explained in this section.

\subsection{Stress field simulation}\label{stress}

\begin{figure}[htp]
\vspace{-0.12in}
\centering
	\includegraphics[keepaspectratio,width=0.7\textwidth]{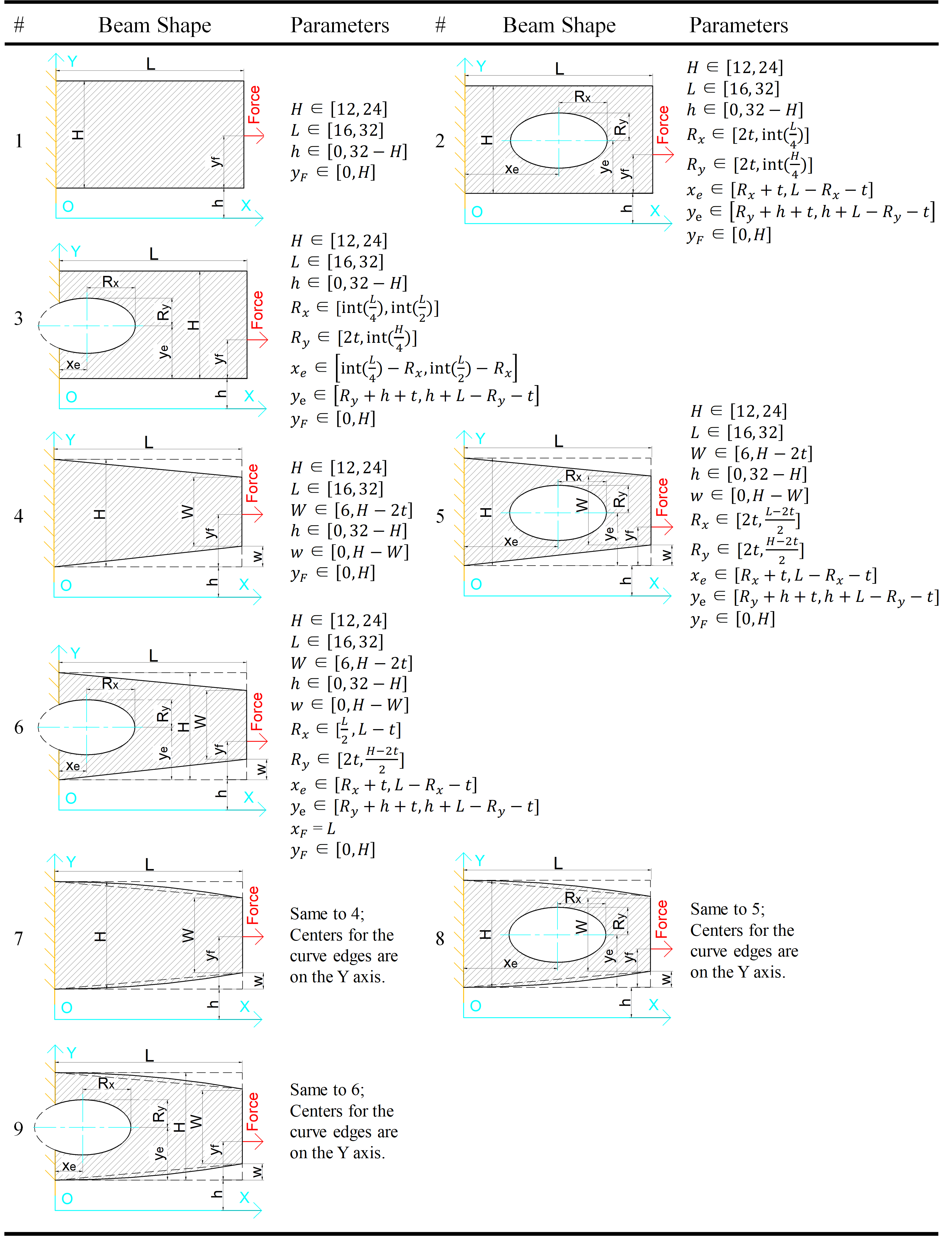}
	\caption{Dataset design parameters}
	\label{fig:3_1_BeamDataset}
	\vspace{-0.12in}
\end{figure}

According to \cite{nie2020stress}, a balanced FEA simulation dataset with 9 types of cantilever beam models is summarized and shown in Fig.~\ref{fig:3_1_BeamDataset} which includes all the beam geometries from \cite{nie2020stress}. All the design parameters for the cantilever beams geometry are generated uniformly-randomly according to the range presented in Fig.~\ref{fig:3_1_BeamDataset} with a minimum increment of 0.1mm. In the case of beam structures with holes, certain constraints have been considered to prevent thin walls with a thickness less than 1mm. The fixtures ($S_u=0$) in this dataset are always on the left edge of the cantilever beam while the external force ($S_p$) is exerted on a random position selected from the right edge of the beam. The magnitude and the angle of the force are uniform-randomly generated with a minimum increment of 100N and $\frac{\pi}{6}$ respectively, ranging from 100N to 1000N and from 0 to $2\pi$. Although the generated cantilevered structures cannot represent all the beam structures, the selected ones represent the most common structures in the mechanical areas \cite{nie2020stress}. Also, following the similar parametric design pattern, other features can be considered and added into the dataset in the future according to previous synthetic CAD design researches \cite{zhang2018featurenet, peddireddy2021identifying}. 

The ground truth of the stress field distribution is calculated by a commercial FEA software - ABAQUS \cite{abaqus2011abaqus}. The chosen material for the simulations is the A36 steel with a Young’s modulus of 200 GPa and a Poisson’s ratio of 0.32. The process of triangular tessellation for the cantilever beam structures is carried out by ABAQUS built-in functions with the approximate global mesh size of 1.0 mm. The settings for deviation factor and minimum size factor are both chosen to be the default value 0.1mm to guarantee a fine triangular-mesh approximation. ABAQUS conducts the FEA algorithm explained in Sec.~\ref{data} to calculate the displacement field and the stress field for the cantilever beams. We choose von Mises stress ($\sigma_{\textup{vm}}$) as the target stress field for the final regression, which can be represented by:

\begin{equation}
\sigma_{\textup{vm}} = 
\sqrt{\sigma^2_{x} + \sigma^2_{y} -
\sigma_{x} \sigma_{y} +
3 \, \tau^2_{xy} }
\label{vm}
\end{equation}

We encode the triangular meshes with their boundary condition into the graph according to the BOGE approach explained in Sec.~\ref{data}. The external force has been normalized in vertex features and the material properties have not been included since only one material is considered. Table~\ref{tab:2_1_Dataset} provides details regarding the generated dataset.

\subsection{Topology optimization}\label{topo}

Based on the simulated stress field, we calculate topological optimization results to validate the abstract decision-making capability of our approach. ABAQUS employs density-based Simple Isotropic Material with Penalization (SIMP) method to optimize the cantilever beam structure. The optimization process works toward minimizing the compliance function $C(z)$ which is the sum of the strain energy of all the elements:

\begin{equation}
\left.\begin{matrix}
x_\textup{min} & :C(z)=\mathbf{U^{T}}\mathbf{K}\mathbf{U} & =\sum_{e=1}^{N} (z_e)^{p} \mathbf{u_e^{T}}\mathbf{k_e}\mathbf{u_e} \\ 
\textup{s.t.} & \frac{V(z)}{V_0}=V_f & \\ 
 & \mathbf{K}\mathbf{U}=\mathbf{F} & \\ 
 & 0 < z_{min} & 
\end{matrix}\right\}
\label{simp}
\end{equation}
\newline Where $V(z)$ and $V_0$ represent the optimized material volume and the model’s design domain volume; $V_f$ is the volume friction, which has been chosen as 0.3 for the final optimization target. Variable $z$ represents the element material density that varies from $z_{min}$ (void if 0, but generally nonzero to avoid singularity) to 1 (the full solid); $p$ is the penalization power whose value is usually 3. \cite{lee2020cnn, nie2021topologygan} Using ABAQUS built-in optimization algorithm, we set the maximum design cycle as 25 and solve the SIMP problem. The optimized cantilever beams have a final $V_f\in (0.29,0.3]$ which satisfies the required design target. The optimized $z_e$ for each element has been exported as the regression ground truth for GNN to be trained.

We finish all the dataset preparation tasks on a personal computer with an Intel i7-10700k CPU. The data generation algorithm runs in less than 10 threads with the total CPU occupation less than 90\%, which has not influenced the evaluation of the computation efficiency of ABAQUS. The average computation time for the stress prediction is around 11.5s while that for the topological optimization is around 489s. We provide a balanced dataset with 45k simulation results (5k simulations for each class from Table~\ref{fig:3_1_BeamDataset} with 3.5k/0.75k/0.75k simulations for training/validating/testing subsets) for both stress field distribution predictions and topological optimizations. Another balanced dataset with 180k simulation results (20k simulations for each class from Table~\ref{fig:3_1_BeamDataset} with 14k/3k/3k simulations for training/validating/testing subsets) for stress field predictions only is generated to investigate the influence of the dataset size.

\begin{table}[]
\centering
\resizebox{\textwidth}{!}{
\begin{tabular}{cccccccc}
\hline
\multirow{2}{*}{Dataset}    
& \multirow{2}{*}{\begin{tabular}[c]{@{}c@{}}Number of data in training/ \\ validating/testing subsets\end{tabular} }
& \multicolumn{3}{c}{Physical quantity}                                                                 
& \multirow{2}{*}{\begin{tabular}[c]{@{}c@{}}Avg. number\\of meshes\end{tabular}}  
& \multirow{2}{*}{\begin{tabular}[c]{@{}c@{}}Computation\\time (s)\end{tabular}} \\

 & & Physical quantity                                                                      & Avg.±Std.    & (Min.,Max.) & \multicolumn{2}{c}{}                                       &                                                                                 \\ \hline
\begin{tabular}[c]{@{}c@{}}Stress prediction\\(45k)\end{tabular}       & 31.5k/6.75k/6.75k                                                                                                 & \begin{tabular}[c]{@{}c@{}}Von Mises stress \\ ($\sigma_\text{vm}$) - MPa\end{tabular} & 2.064.±2.239 & (0,74.09)   & 857                                   & 11.53                                                                           \\
\begin{tabular}[c]{@{}c@{}}Stress prediction\\(180k)\end{tabular}    & 126k/27k/27k                                                                                                      & \begin{tabular}[c]{@{}c@{}}Von Mises stress \\ ($\sigma_\text{vm}$) - MPa\end{tabular} & 2.057.±2.232 & (0,74.09)   & 857                                    & 11.58                                                                           \\
\begin{tabular}[c]{@{}c@{}}Topology optimization\\(45k)\end{tabular} & 31.5k/6.75k/6.75k                                                                                                 & \begin{tabular}[c]{@{}c@{}}Volume density \\ ($z_e$)   \end{tabular}                  & 0.298.±0.389 & (0,1)       & 857                                    & 489.09                                                                          \\ \hline
\end{tabular}}
\caption{Details for the dataset}
\label{tab:2_1_Dataset}
\end{table}

The parameterized synthetic dataset can be generalized to other fields related to the FEA surrogate DL model. The synthetic dataset sometimes has an expensive time cost that can become a bottleneck for our approach to be a general-purpose surrogate model. However, the dataset prepared by physical experiments can also serve the training task. The stress fields can be measured by installing strain gauges on real cantilever beams on the area whose position corresponds to that of the predefined meshes in the simulation. Multiple data can be generated within a short time when the external forces change. In this condition, the physical experiment can serve as a fast PDE solver which assists to generate sufficient data for the surrogate model. All of these ideas contribute to the realization of a generalized model with high-precision predictions under various working conditions.

\section{Neural Network Training and Results}\label{train}

In order to present the efficiency of the BOGE approach, we compare its performance with the conventional graph embedding methods (without any graph shortcuts) on different state of art GNN structures (GCN \cite{kipf2016semi}, GAT \cite{velivckovic2017graph}, UNet \cite{gao2019graph}, DeepGCN \cite{li2019deepgcns, li2020deepergcn}). We choose the MSE (Mean Square Error) function as the loss function and the evaluation metric to quantify our prediction performance \cite{nie2020stress, koeppe2020intelligent}. The predicted outputs – the magnitude of von Mises stress $\sigma_\textup{vm}$ and the volume density $z_e$, are all 1D tensors which can be written as $\hat{p}=(\hat{p}_1, \hat{p}_1,..., \hat{p}_n)$, while the ground truth results can be represented by $p=(p_1, p_2, ...,p_n)$. The MSE can be represented as:

\begin{equation}
\mathrm{MSE}=\frac{1}{n}\sum_{i=1}^{n}(p_i-\hat{p}_i)^2
\label{mse}
\end{equation}

According to \cite{koeppe2020intelligent}, we also employ the MAPE (Mean Absolute Percentage Error) to present the error rate in percentage, which can be expressed by:

\begin{equation}
\mathrm{MAPE}=\frac{100}{n}\sum_{i=1}^{n}|\frac{p_i-\hat{p}_i}{p_i+\varepsilon}|
\label{mape}
\end{equation}
\newline Where we set $\varepsilon=0.01$ to avoid computation error when the ground truth result is zero or extremely small.

All the codes have been written in Pytorch (shared in Github\footnote{\url{https://github.com/stxyfu/Stress\_GNN}}) with Pytorch Geometric library \cite{fey2019fast} and run on an NVIDIA RTX 3090 GPU and an Intel i7-10700k CPU. For all training experiments, the ADAM optimization algorithm has been employed to adjust the hyperparameters in the model. The training batch size is set to 8 to make maximal use of GPU memory and provide a precise training gradient while the learning rate is 0.01 obtained through the grid search to provide the best performance of the final GNN models. 500 training epochs have been employed for each experiment to ensure the GNNs are all fully trained with suitable prediction accuracy. The entire training process for the selected GNN models (Table~\ref{tab:3_3_stress} and Table~\ref{tab:3_4_topo}) takes around 125 hours on 45k datasets and 493 hours on the 150k dataset. The training time can be largely reduced if multiple high-performance GPUs are employed. 

\subsection{Stress field prediction with conventional graph embedding}\label{deepGNN}

We first investigate the regression capability of GNNs with conventional graph embedding methods on the 45k stress field prediction dataset, shown in Table~\ref{tab:2_1_Dataset}. Conventional graph embedding methods only consider one-hop graph links (with $l_e = 0$ in our condition) whose adjacency matrix does not contain any shortcuts for boundary elements. Using this graph embedding method, shallow layer GNN ("GCN(3 layers)" in Table~\ref{tab:3_1_results}) cannot pass the boundary information to all the internal elements which cannot perform any regression at the end of the training. Fig.~\ref{fig:ConvRes}(\subref{fig:Conv3GCN}) illustrates that, with the 3-layer GCN, the boundary information only propagates through limited meshes near the model's geometrical boundary and cannot reach the internal elements inside the material. The deep-layer GNN can save the weighted features passed from other elements in its hidden state $h$ according to Eq.~(\ref{mpnn1}) and provides a relatively long-range message passing. According to our statistics, the largest graph distance between two structured meshes is around $l_e=250$ which requires at least 8 layers of GNN for general MPNN layers. However, deep GNN usually incurs oversmoothing that makes the embedded features similar and prevents the further regression process. For example, in Fig.~\ref{fig:ConvRes}(\subref{fig:Conv8GCN}), the 8-layer GCN ("GCN(8 layers)" in Table~\ref{tab:3_1_results}) cannot predict anything but only outputs uniform values after 500 training epochs, which illustrates this extreme oversmoothing phenomenon. In order to solve this problem, researchers have employed different methods to mitigate the influence of oversmoothing and realized a few state-of-art GNN structures with deep GNN layers. From those methods, the most popular anti-oversmoothing techniques are DeepGCN \cite{li2019deepgcns, li2020deepergcn}, attention layers(e.g., GAT  \cite{velivckovic2017graph}) and DropEdge \cite{rong2019dropedge}. Some other GNN structures, though have been developed for other purposes, can also reach deep layers with sufficient accuracy, for instance, the g-U-net \cite{gao2019graph}. However, most of these methods mainly work on classification problems. For the regression problem in our condition, extra validations are needed. Therefore, we have tested those models using the conventional graph embedding approach and shown those results in Table~\ref{tab:3_1_results} and Fig.~\ref{fig:ConvRes}.

\begin{table}[!h]
\centering
\includegraphics[keepaspectratio, width=0.7\columnwidth]{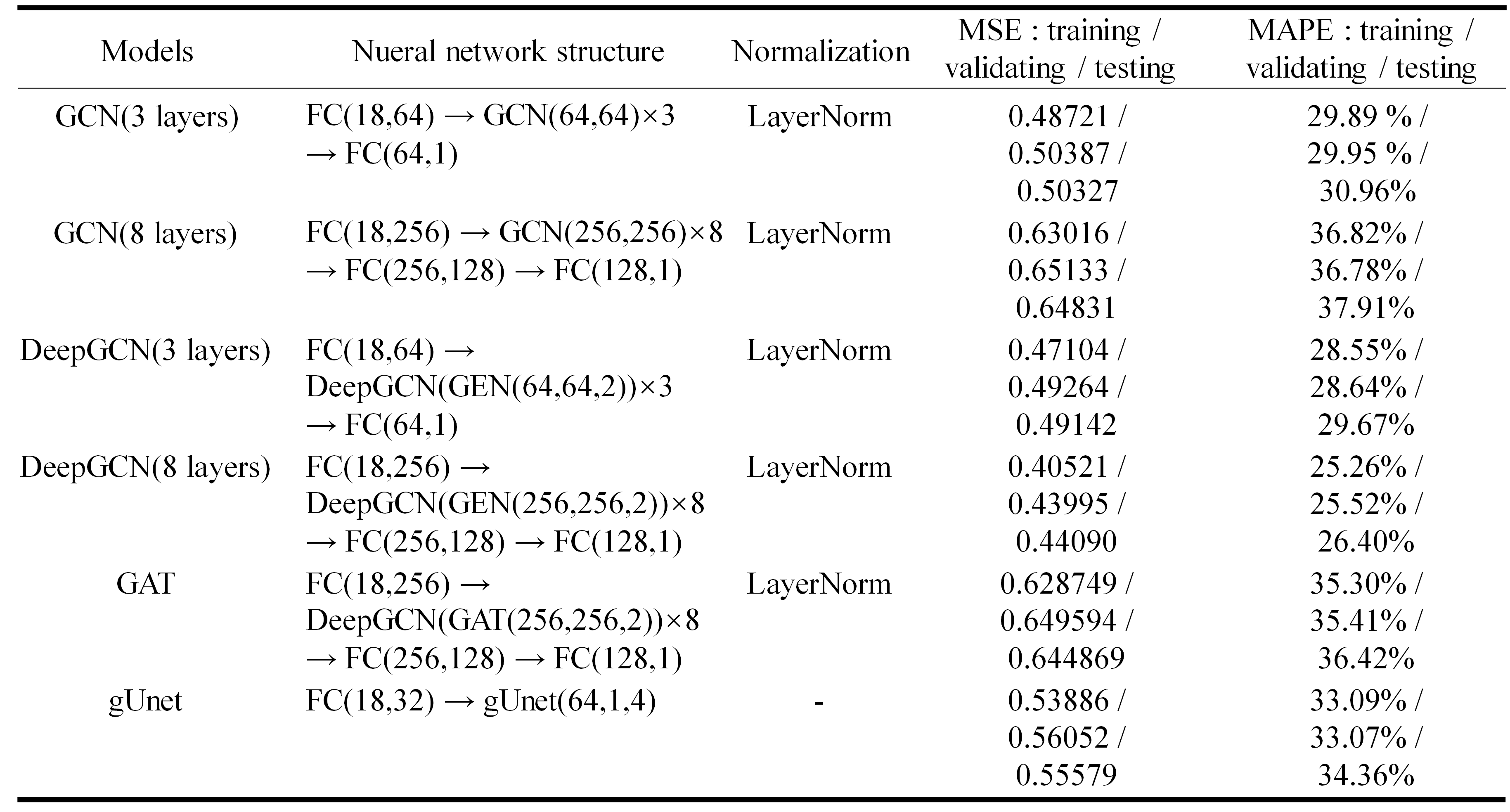}


\resizebox{0.7\textwidth}{!}{

\begin{tabular}[l]{ll}
\hline
FC($f_{in}$, $f_{out}$): Fully connected layer;                                        & $f_{in}$: Number of input channels;\\
GCN($f_{in}$, $f_{out}$): Graph convolutional layer;\cite{kipf2016semi}                & $f_{out}$: Number of output channels;\\
GAT($f_{in}$, $f_{out}$, $H_{GAT}$): Graph attentional layer;\cite{velivckovic2017graph}  & $H_{GAT}$:  Number of multi-head-attentions;\\
gUnet($f_{in}$, $f_{out}$, $D_{gUnet}$): Graph U-net;\cite{gao2019graph}                   & $D_{gUnet}$:  The depth of the U-Net architecture;\\
GEN($f_{in}$, $f_{out}$, $L_{MPL}$): GENeralized Graph Convolution layer;\cite{li2020deepergcn}  & $L_{MPL}$:  The number of MLP layers;\\
\multicolumn{2}{l}{DeepGCN(*): DeepGCN layer (\text{Graph convolutional layer}$\to$\text{Normalization}$\to$\text{Activation}$\to$
\text{ResGCN}$\to$\text{Dropout});\cite{li2019deepgcns} }\\
\hline
\end{tabular}}
\caption{Training results with different GNN structures and graph embedding methods (training/validating/testing accuracy are presented without regularization)}
\label{tab:3_1_results}
\end{table}

DeepGCN employs the GEN (GENeralized Graph Convolution layers \cite{li2020deepergcn}), layerNorm \cite{li2019deepgcns} and ResNet structure \cite{he2016deep}, making GNN sufficiently robust against oversmoothing and can reach 56 layers with high classification accuracy. However, 8-layers DeepGCN ("DeepGCN(8 layers)" in Table~\ref{tab:3_1_results}) with the ResGCN backbone \cite{li2019deepgcns} still cannot reach a high accuracy and performs some oversmoothing shown in Fig.~\ref{fig:ConvRes}(\subref{fig:Conv8DeepGCN}). GAT employs the attention mechanism that provides weighted attention coefficients to various graph features to ensure that GNN can pass the important information through deep layers. We employ the same DeepGCN backbone, but substitute GAT (with 4 attention heads) for GEN and try to improve the prediction accuracy, shown as "GAT" in Table~\ref{tab:3_1_results}. According to our tests, the combination of GAT and the DeepGCN backbone performs better than simply stacking GAT layers on our dataset. However, the training result still appears to be incapable of obtaining acceptable results and shows an extremely oversmoothed result, shown in Fig.~\ref{fig:ConvRes}(\subref{fig:Conv8GAT}). The g-U-Net can encode and decode both high-level features and local features which can reach 5 layers in \cite{gao2019graph}. However, each graph pooling layer drops both graph vertices (mesh elements) and graph links. When only one-hop graph links ($l_e=0$) are considered, with a very sparse adjacency matrix, stacking pooling layers in the g-U-net makes the adjacent matrix too sparse, so that only a few graph links remain in high-level feature scope, which cannot pass ample information to other elements through GNN and ruin the g-U-net's prediction capability. Also, a deeper g-U-net can still incur the extreme oversmoothing problem shown in Fig.~\ref{fig:ConvRes}(\subref{fig:ConvgUnet}) (with the training result shown as "gUnet" in Table~\ref{tab:3_1_results}). Other anti-oversmoothing techniques like DropEdge can ruin the integrity of the graph, which is not suitable for our problem as we have tested. The overall training results show that deep-layer GNNs with the conventional graph embedding approach cannot simply regress the boundary value problem with long-range graph-vertex interactions.



\def \figWidth{0.24} 
\def \linWidthRatio{0.7} 

\begin{figure}[!htbp]
\centering
\includegraphics[width=\linWidthRatio\linewidth]{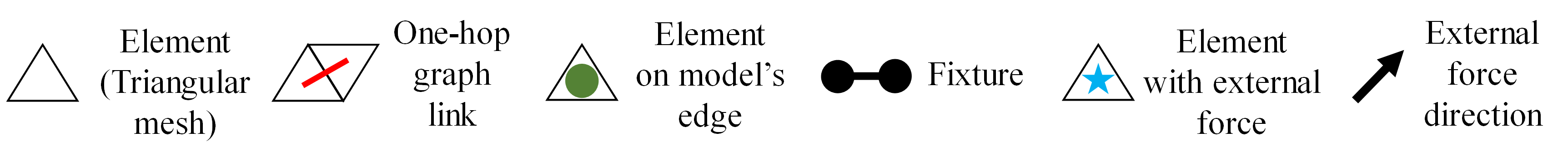}
\newline
\begin{subfigure}{\figWidth\textwidth}
  \centering
  \includegraphics[width=\linWidthRatio\linewidth]{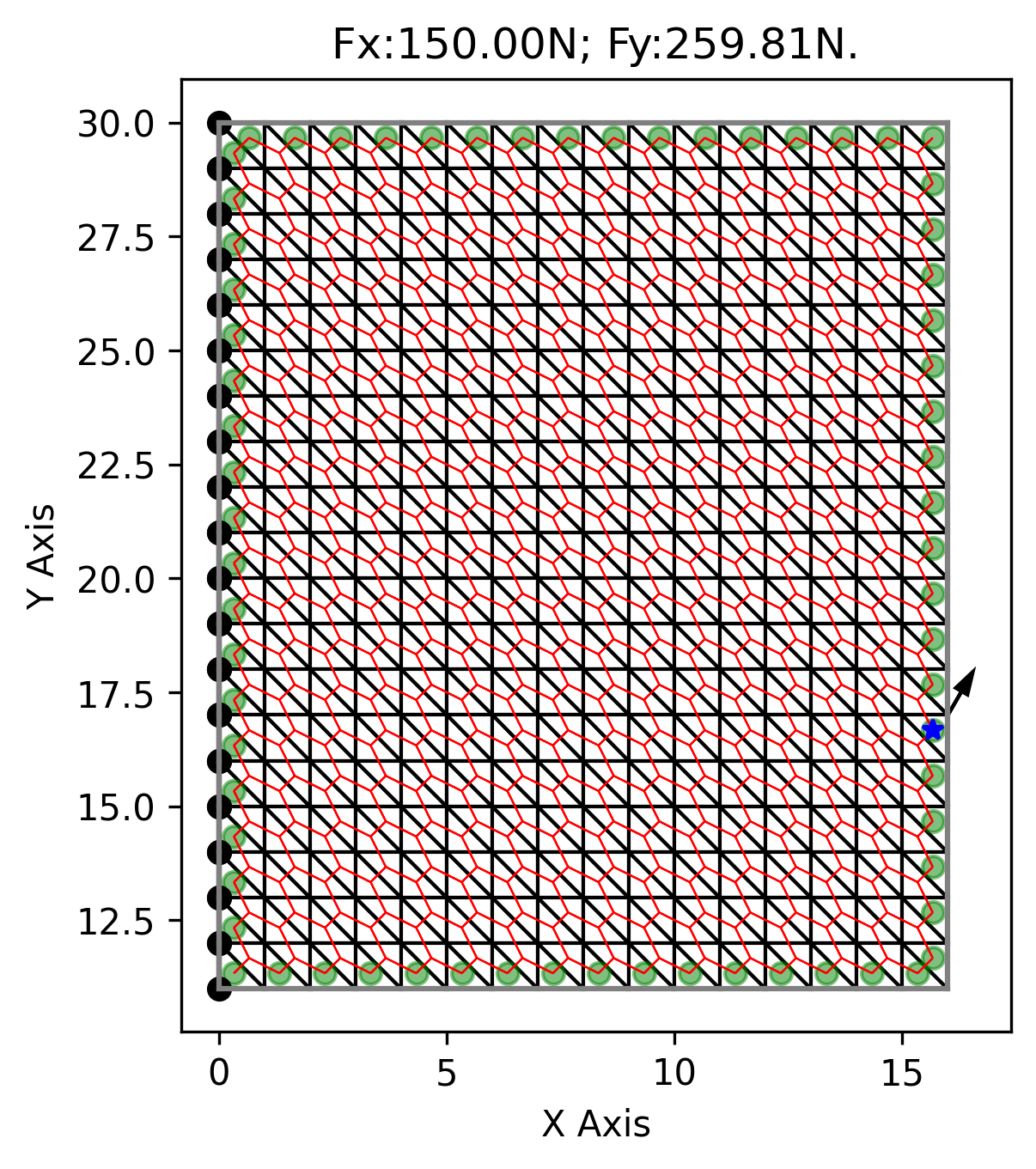}
  \caption{Simulation settings}
  \label{fig:ConvSim}
\end{subfigure}%
\begin{subfigure}{\figWidth\textwidth}
  \centering
  \includegraphics[width=\linWidthRatio\linewidth]{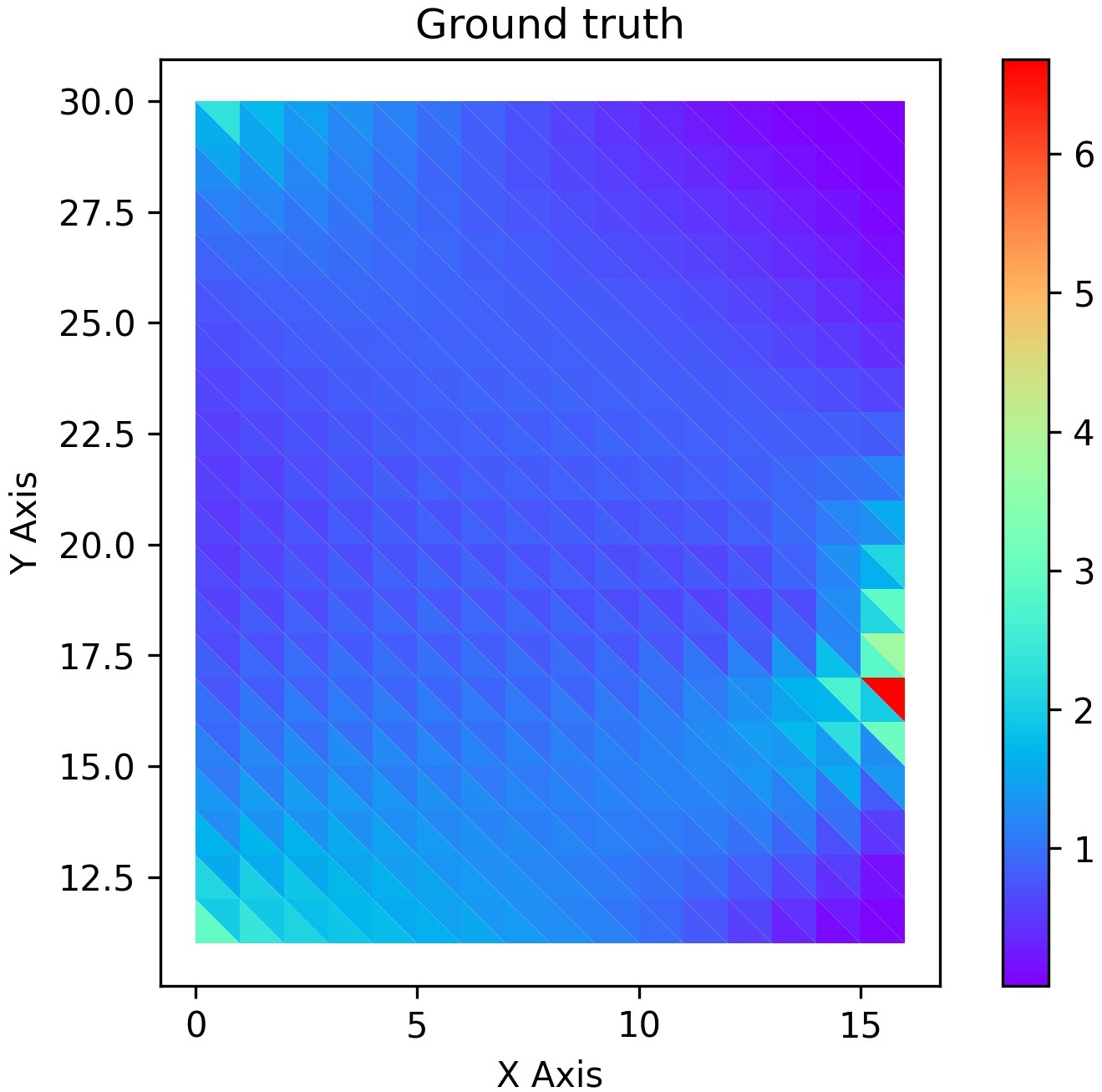}
  \caption{Ground truth}
  \label{fig:ConvGround}
\end{subfigure}
\begin{subfigure}{\figWidth\textwidth}
  \centering
  \includegraphics[width=\linWidthRatio\linewidth]{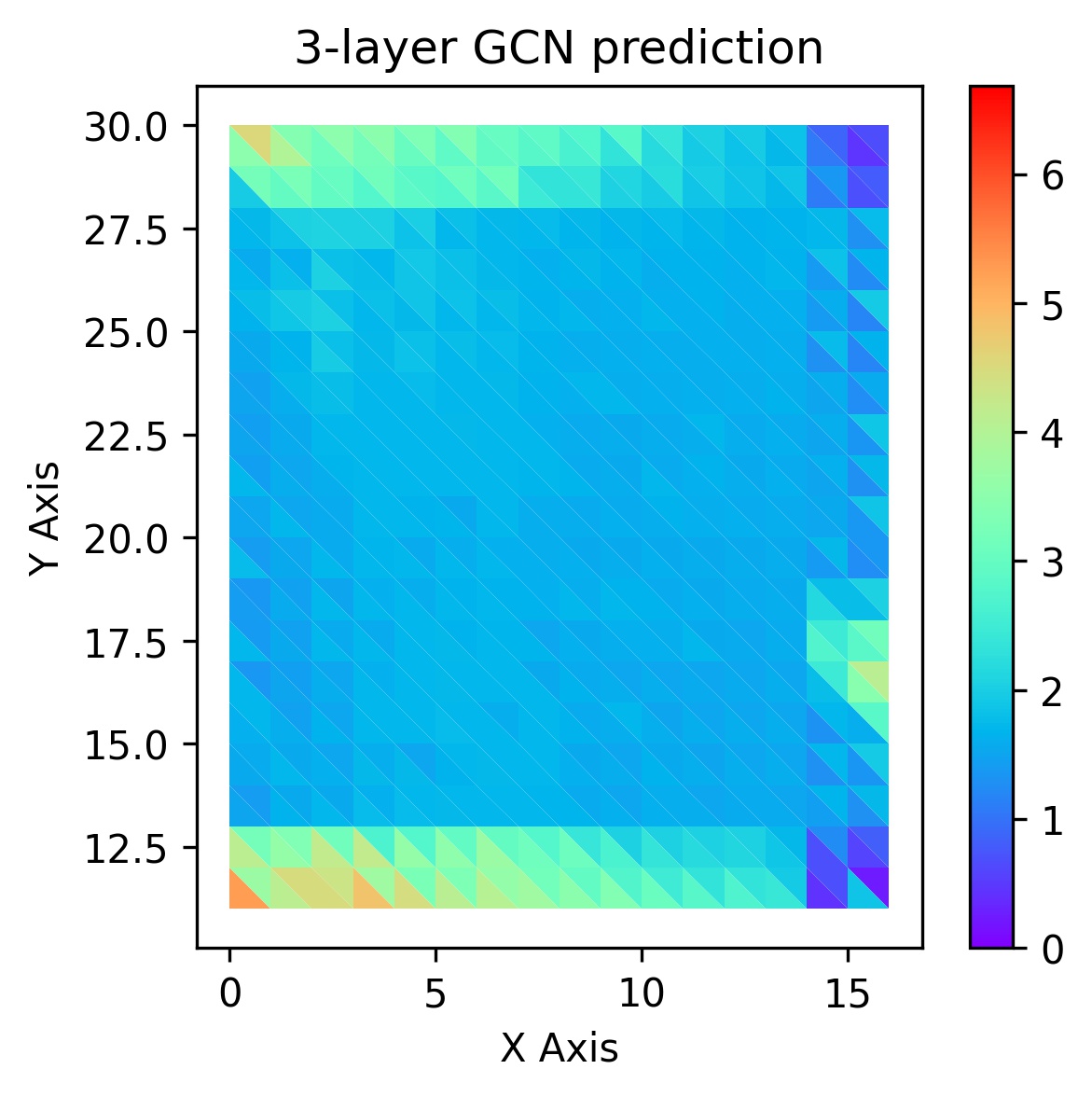}
  \caption{3-layer GCN prediction}
  \label{fig:Conv3GCN}
\end{subfigure}
\begin{subfigure}{\figWidth\textwidth}
  \centering
  \includegraphics[width=\linWidthRatio\linewidth]{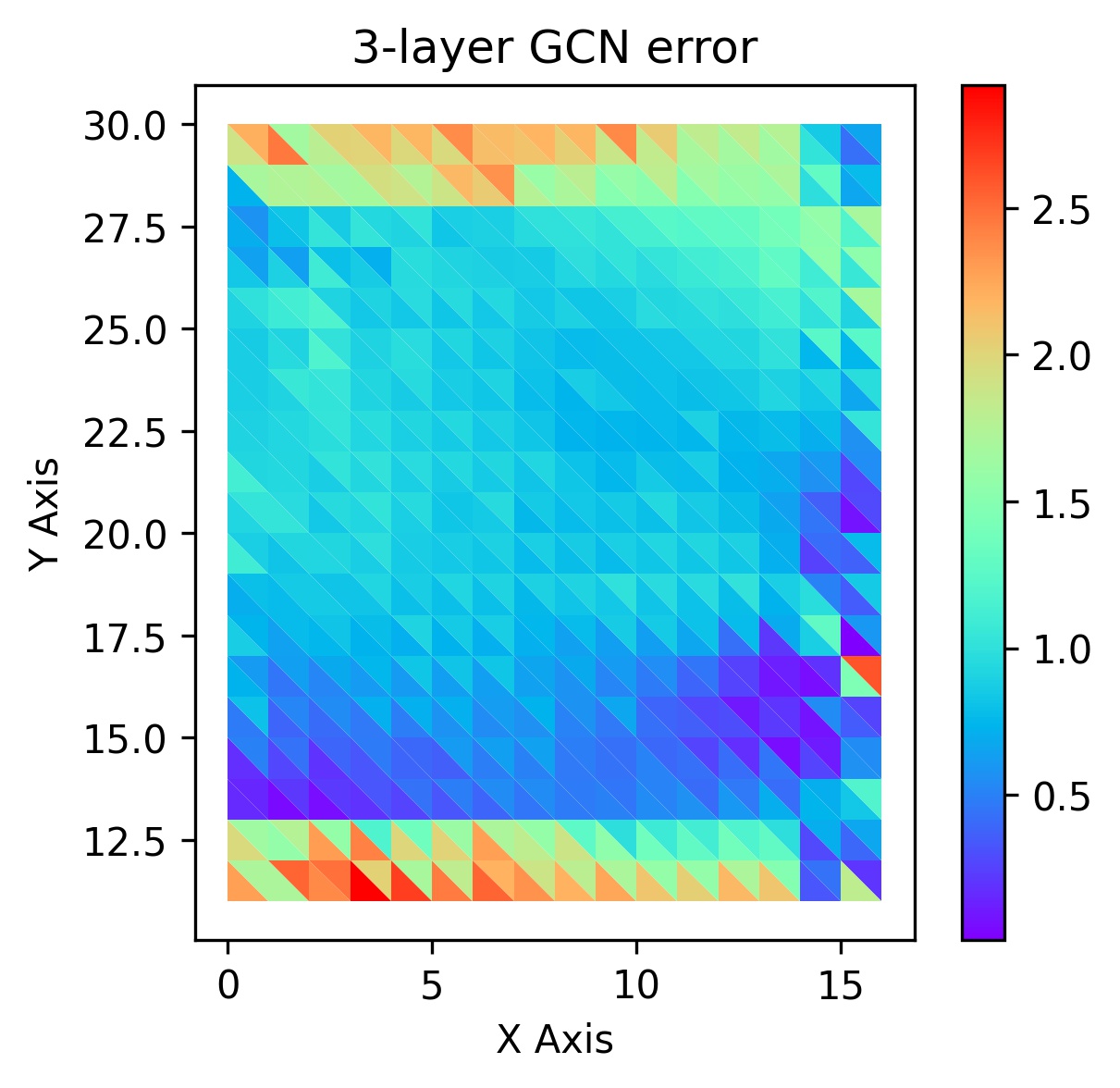}
  \caption{3-layer GCN error}
  \label{fig:Conv3GCNerr}
\end{subfigure}
\end{figure}
\begin{figure}[!htbp]\ContinuedFloat
\begin{subfigure}{\figWidth\textwidth}
  \centering
  \includegraphics[width=\linWidthRatio\linewidth]{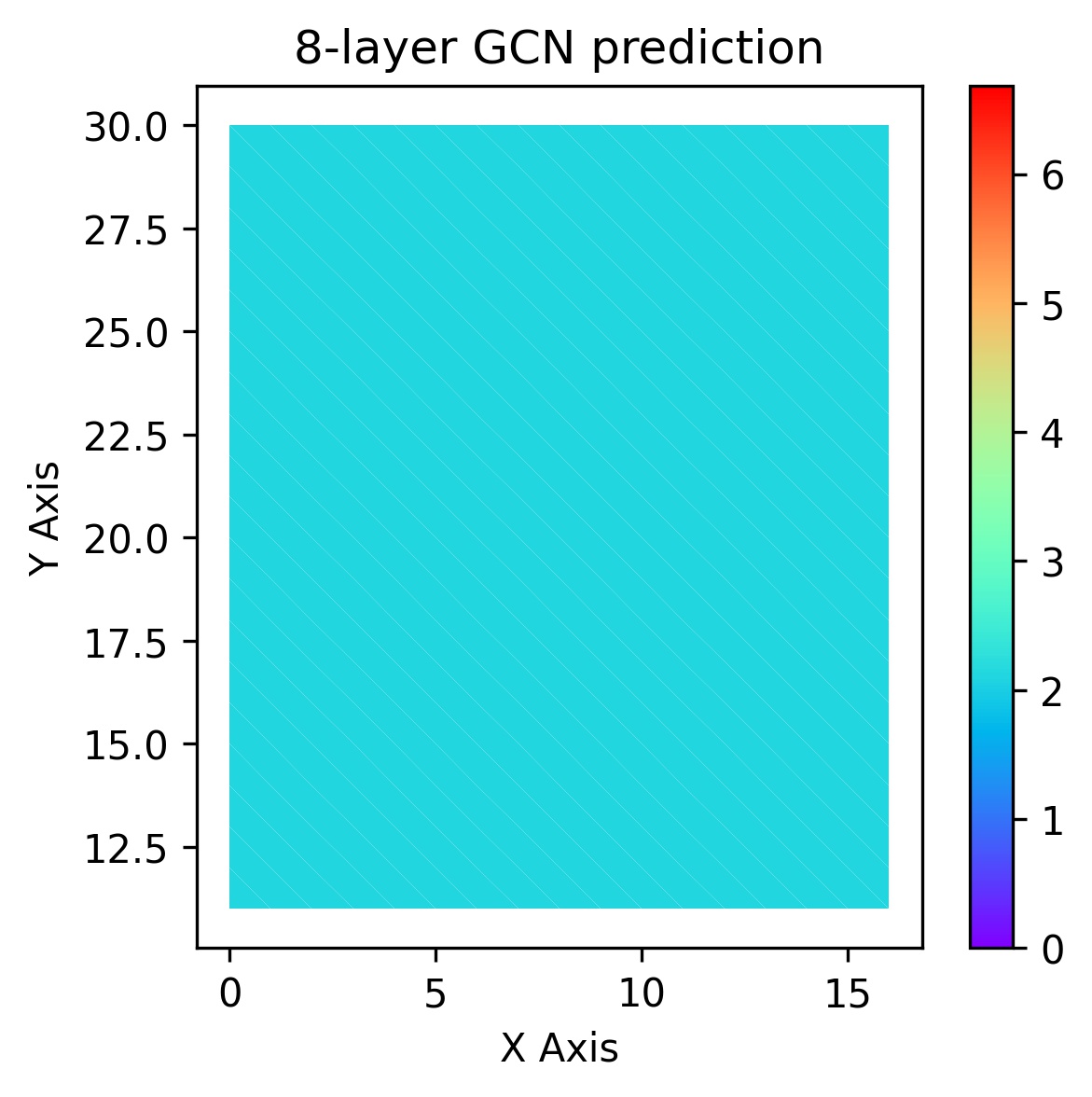}
  \caption{8-layer GCN prediction}
  \label{fig:Conv8GCN}
\end{subfigure}
\begin{subfigure}{\figWidth\textwidth}
  \centering
  \includegraphics[width=\linWidthRatio\linewidth]{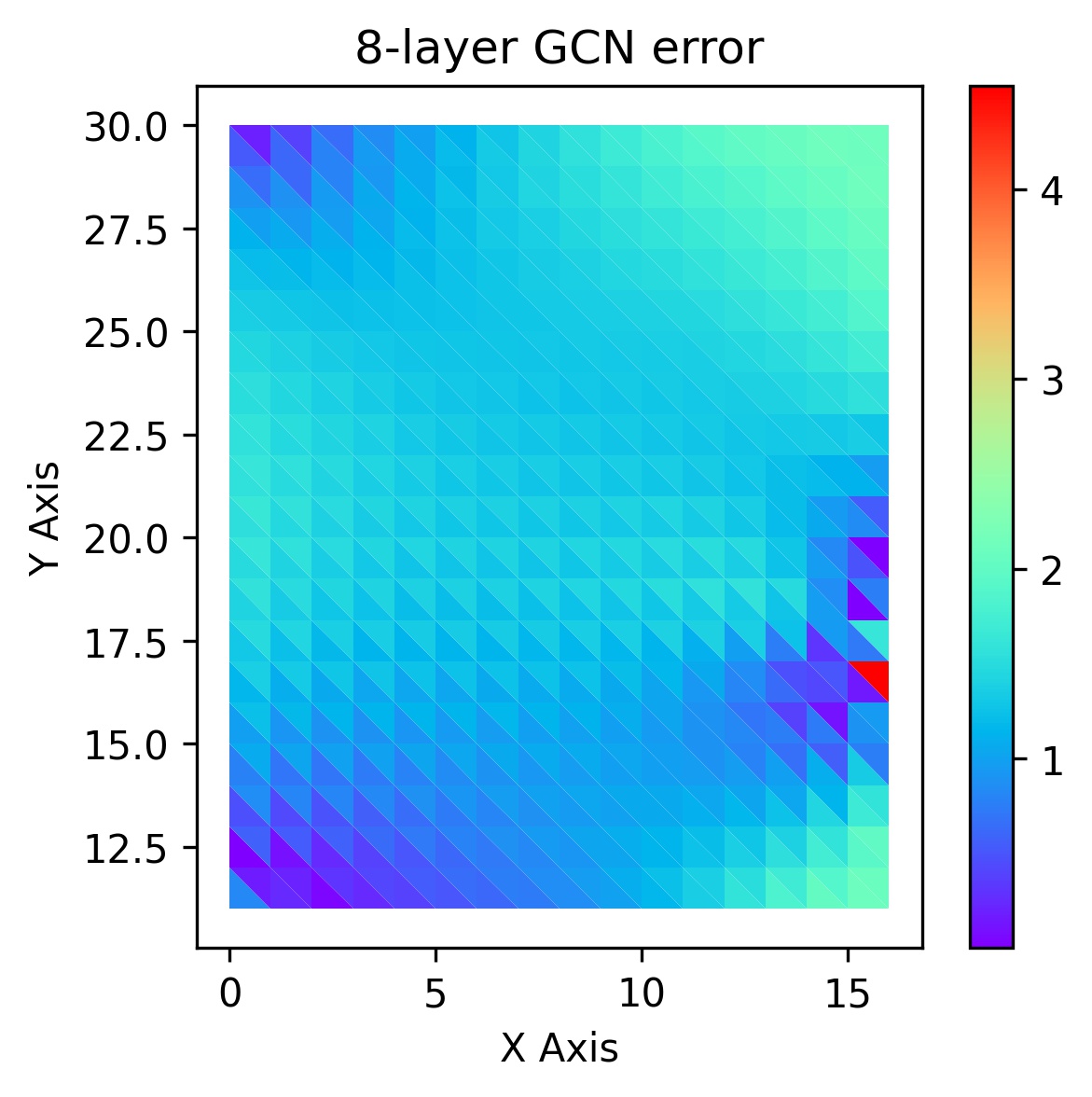}
  \caption{8-layer GCN error}
  \label{fig:Conv8GCNerr}
\end{subfigure}
\begin{subfigure}{\figWidth\textwidth}
  \centering
  \includegraphics[width=\linWidthRatio\linewidth]{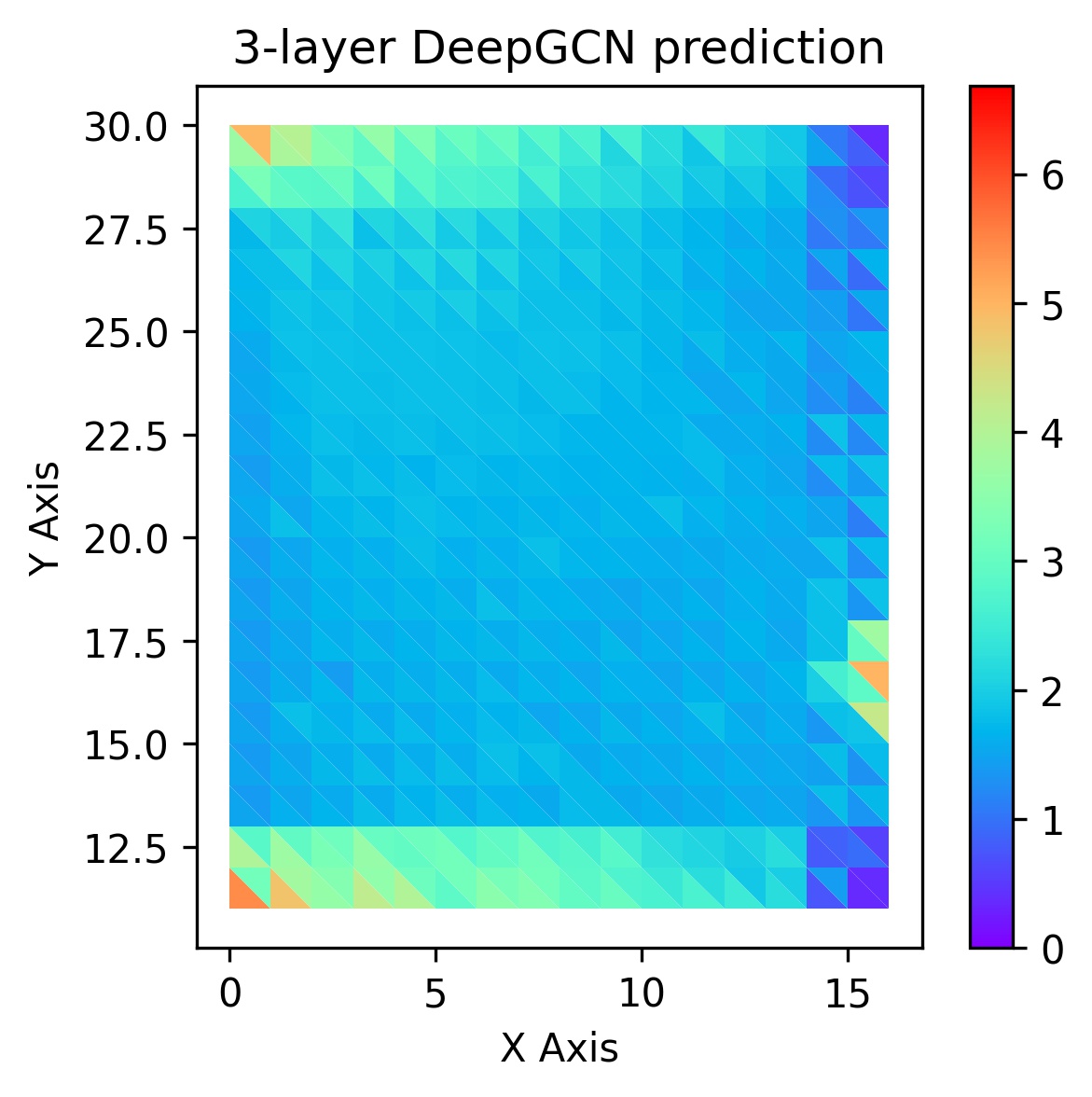}
  \caption{3-layer DeepGCN prediction}
  \label{fig:Conv3DeepGCN}
\end{subfigure}
\begin{subfigure}{\figWidth\textwidth}
  \centering
  \includegraphics[width=\linWidthRatio\linewidth]{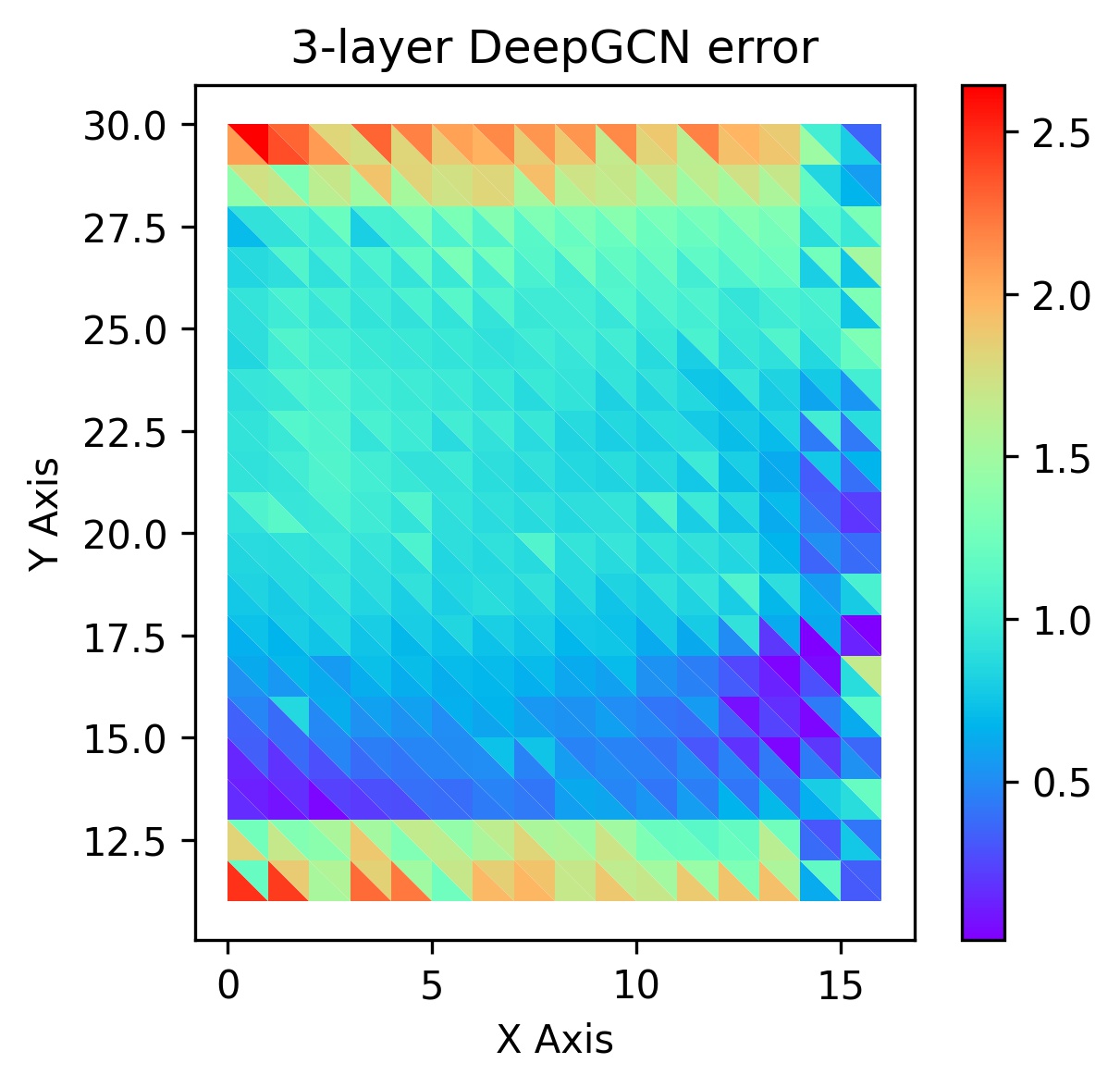}
  \caption{3-layer DeepGCN error}
  \label{fig:Conv3DeepGCNerr}
\end{subfigure}
\end{figure}
\begin{figure}[!htbp]\ContinuedFloat
\begin{subfigure}{\figWidth\textwidth}
  \centering
  \includegraphics[width=\linWidthRatio\linewidth]{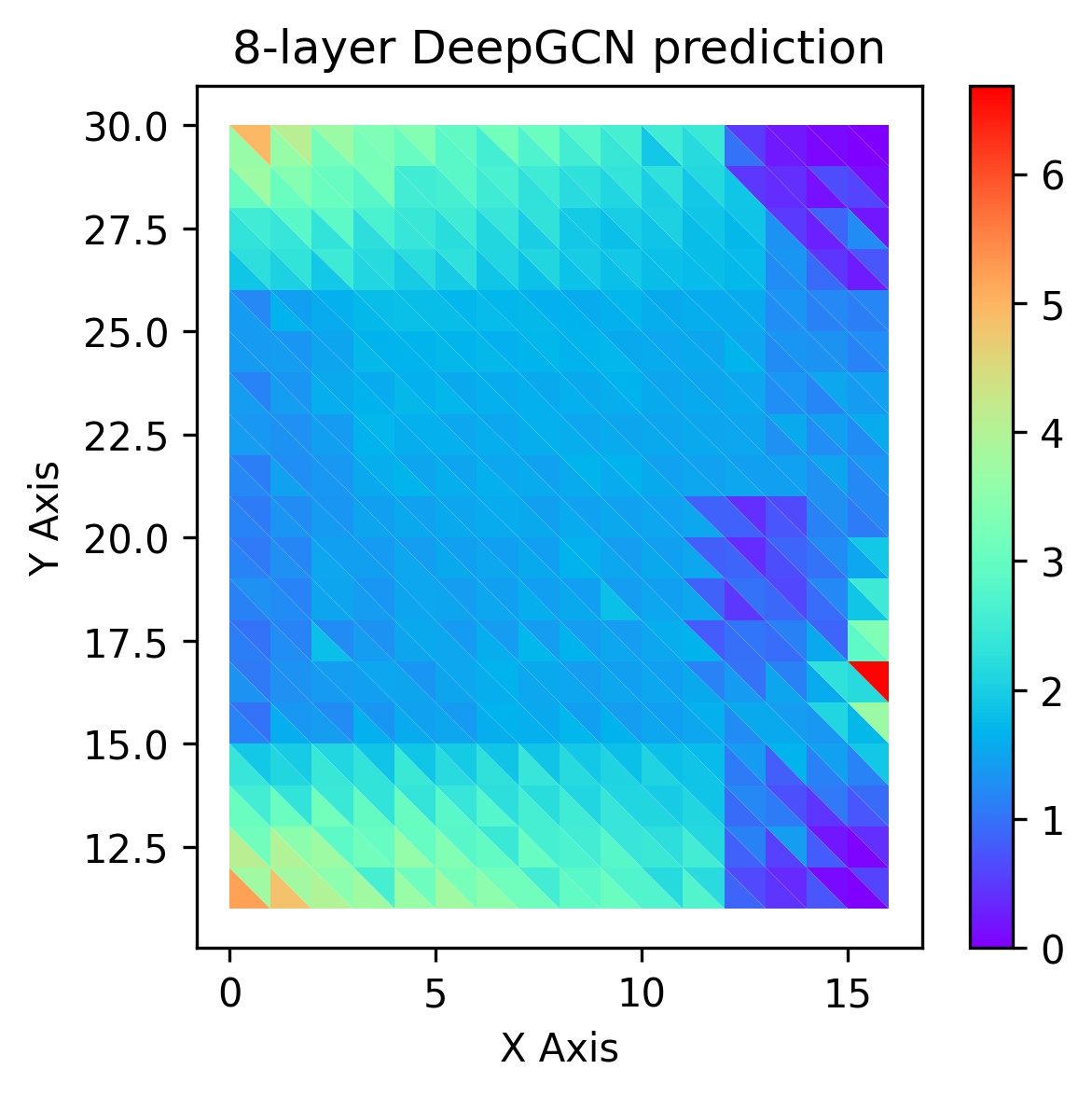}
  \caption{8-layer DeepGCN prediction}
  \label{fig:Conv8DeepGCN}
\end{subfigure}
\begin{subfigure}{\figWidth\textwidth}
  \centering
  \includegraphics[width=\linWidthRatio\linewidth]{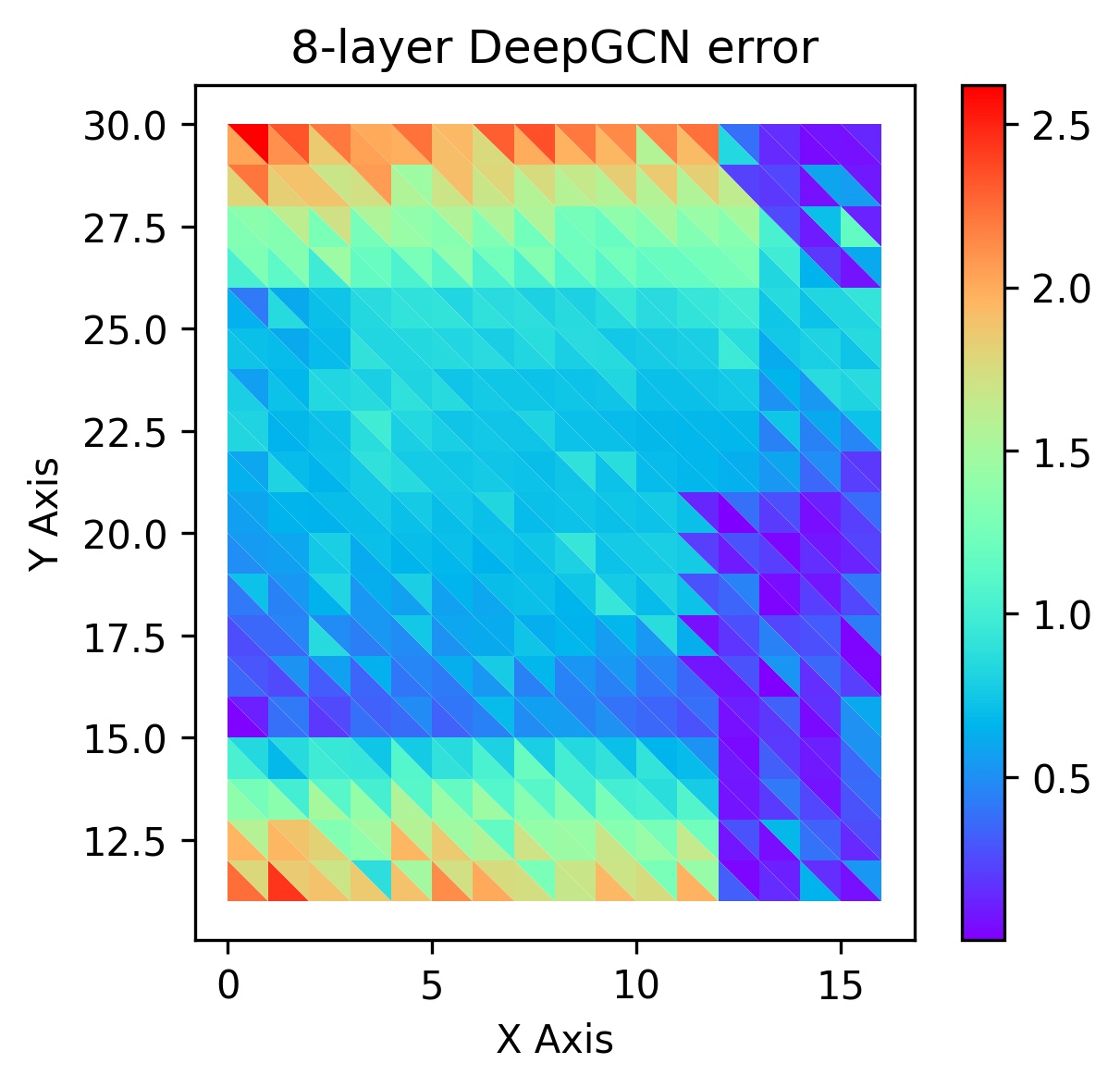}
  \caption{8-layer DeepGCN error}
  \label{fig:Conv8DeepGCNerr}
\end{subfigure}
\begin{subfigure}{\figWidth\textwidth}
  \centering
  \includegraphics[width=\linWidthRatio\linewidth]{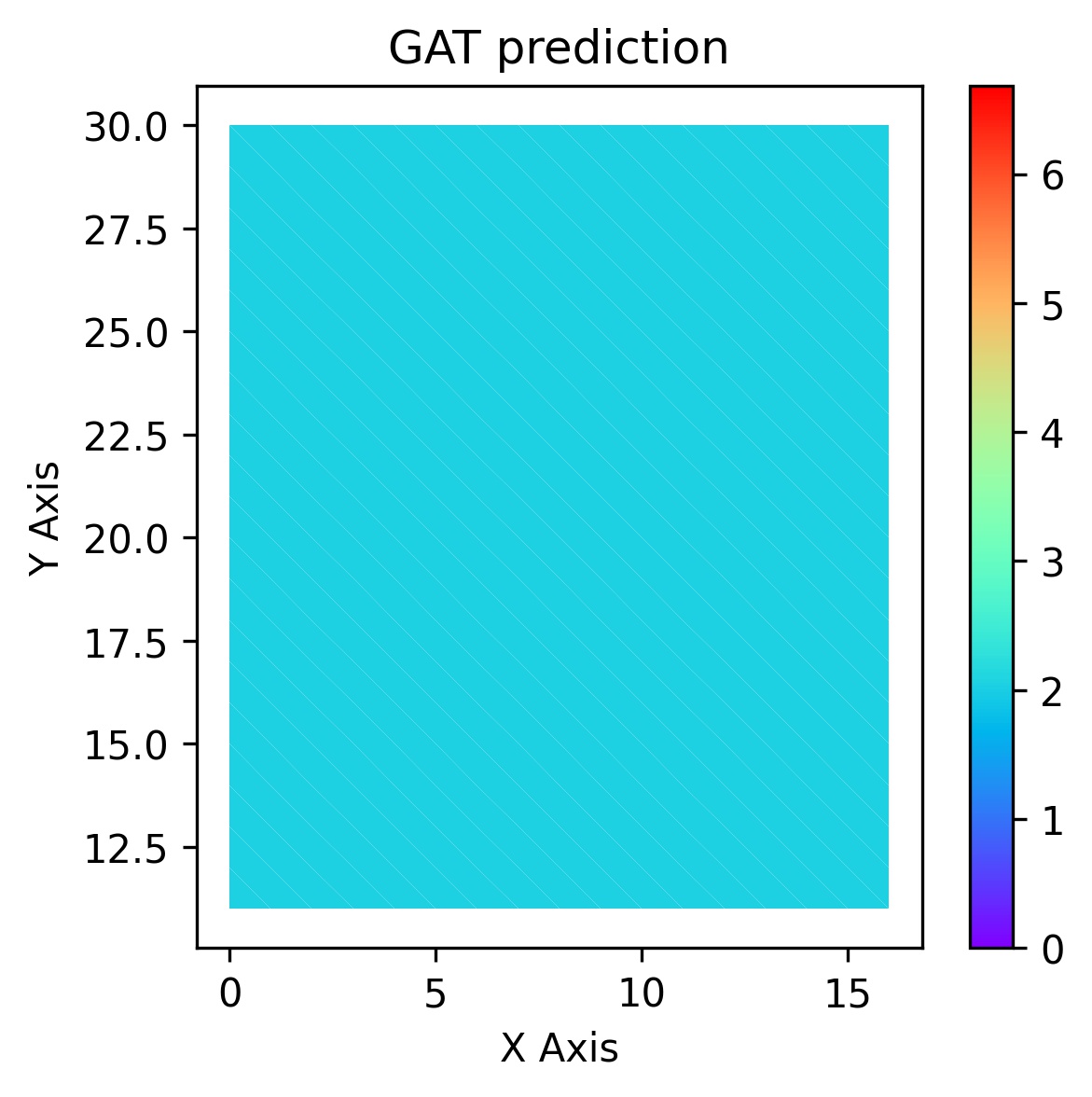}
  \caption{8-layer GAT prediction}
  \label{fig:Conv8GAT}
\end{subfigure}
\begin{subfigure}{\figWidth\textwidth}
  \centering
  \includegraphics[width=\linWidthRatio\linewidth]{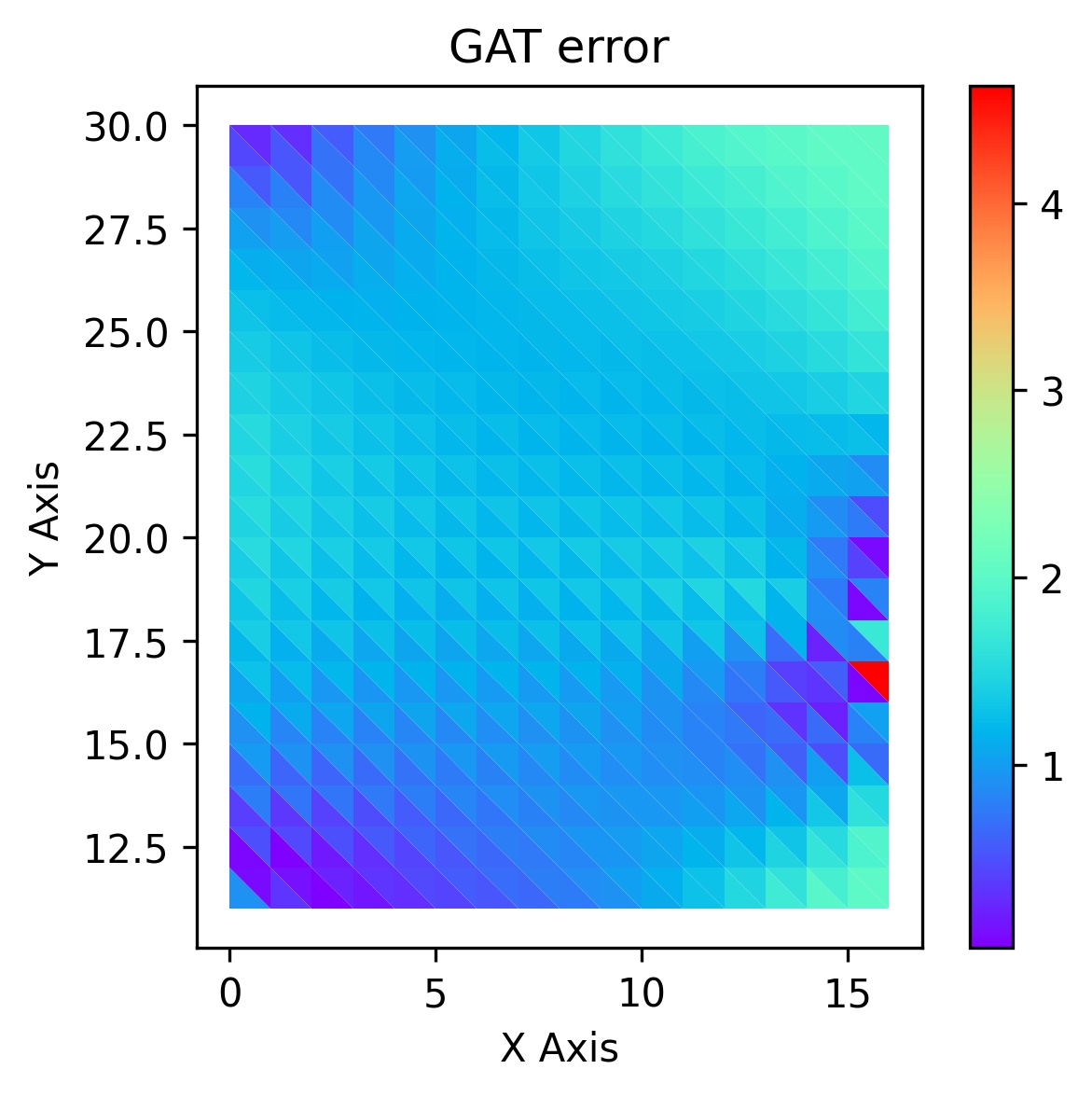}
  \caption{8-layer GAT error}
  \label{fig:Conv8GATerr}
\end{subfigure}
\end{figure}
\begin{figure}[!htbp]\ContinuedFloat
\begin{subfigure}{\figWidth\textwidth}
  \centering
  \includegraphics[width=\linWidthRatio\linewidth]{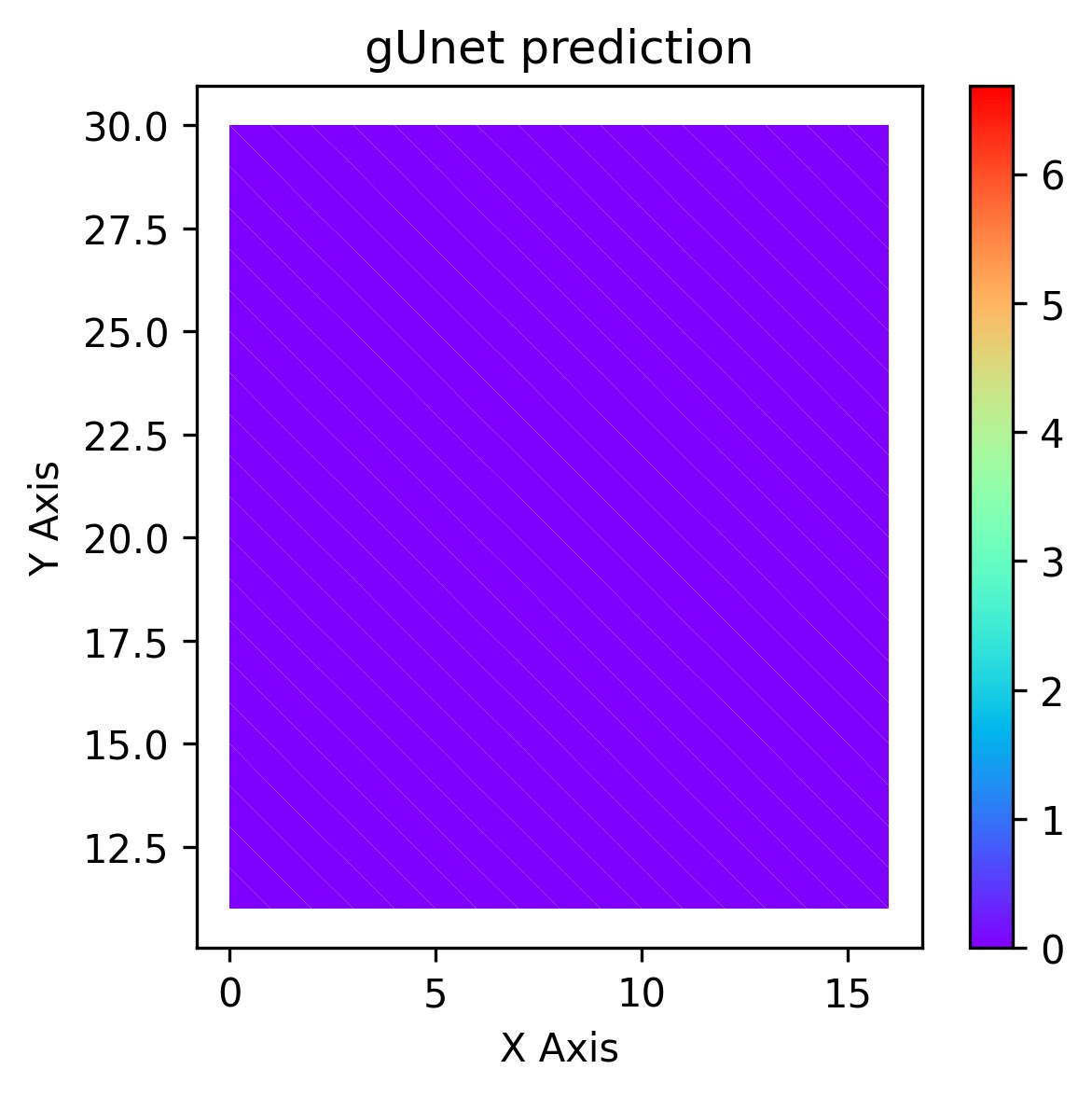}
  \caption{The g-U-net prediction}
  \label{fig:ConvgUnet}
\end{subfigure}
\begin{subfigure}{\figWidth\textwidth}
  \centering
  \includegraphics[width=\linWidthRatio\linewidth]{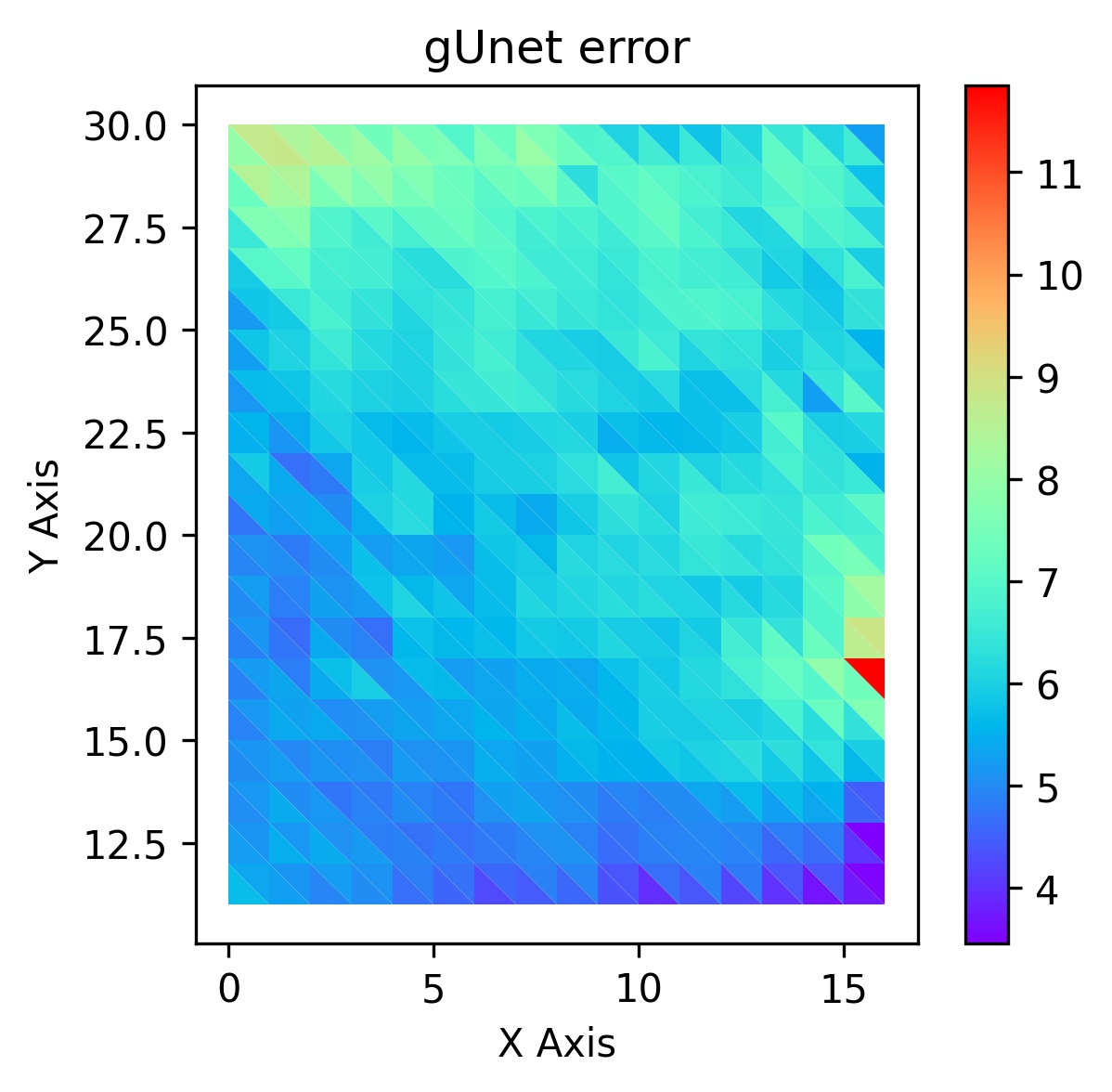}
  \caption{The g-U-net error}
  \label{fig:ConvgUneterr}
\end{subfigure}
\caption{Prediction results for conventional graph embedding methods}
\label{fig:ConvRes}
\end{figure}

\subsection{Stress field prediction with BOGE}\label{stress}

The BOGE approach shows its unique benefits for regressing the results of large-scale boundary value problems with simple shallow-layer GNNs. By providing shortcuts for both global boundary conditions and local information of neighboring elements, GNNs with BOGE bypass the message passing limitation of the MPNN framework, and all the internal elements can directly get the boundary condition information to solve the boundary value problem. We have employed the 3-layers DeepGCN and BOGE with $l_e=1,2,3$ graph links to validate its performance on the 45k stress field prediction dataset. The results are shown in Table~\ref{tab:3_3_stress}. 

Compared with the predicted results from "DeepGCN(3 layers)" in Table~\ref{tab:3_1_results}, BOGE with $l_e=0$ graph links (model \#1 in Table~\ref{tab:3_3_stress}) largely improves the prediction accuracy and reaches 0.016383 testing accuracy with MAPE of 2.85\%, which is accurate enough to serve as an efficient FEA surrogate prediction model. The testing results has been compared for BOGE with different adjacency matrices that encodes local graph links with $l_e=1,2$. From Table~\ref{tab:3_3_stress}, it can be seen that additional local mesh information cannot improve the prediction accuracy. The testing accuracy for $l_e=1$ (model \#2 in Table~\ref{tab:3_3_stress}) drops to 0.021737 with MAPE of 3.08\% while that for $l_e=2$ (model \#3 in Table~\ref{tab:3_3_stress}) drops to 0.027302 with MAPE of 3.44\%. More local graph edges provide worse prediction accuracy, which can be caused by message congestion \cite{loukas2019graph} since more elements enlarge the size of the message passing at each training epoch. Also, the adjacency matrix with more graph links can lead to a denser tensor product which takes more computational resources. Therefore, we only employ BOGE with $l_e=0$ for the rest of the training task to provide a satisfying prediction result.

\begin{table}[!h]
\centering
\includegraphics[keepaspectratio, width=\columnwidth]{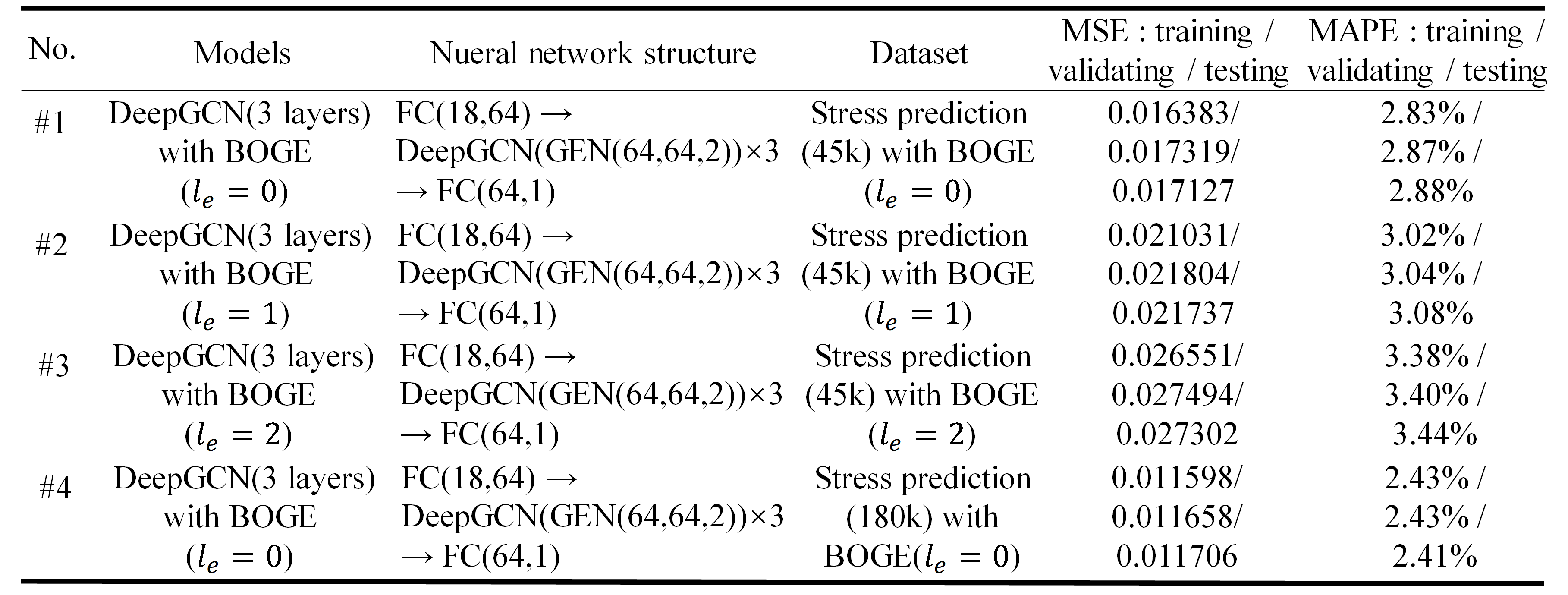}
\caption{Stress prediction results (notations for neural network layers are explained in Table~\ref{tab:3_1_results}; training/validating/testing accuracy are presented without regularization)}
\label{tab:3_3_stress}
\end{table}

\begin{figure}[!htbp]
\centering
\begin{subfigure}{0.4\textwidth}
  \centering
  \includegraphics[width=\linewidth]{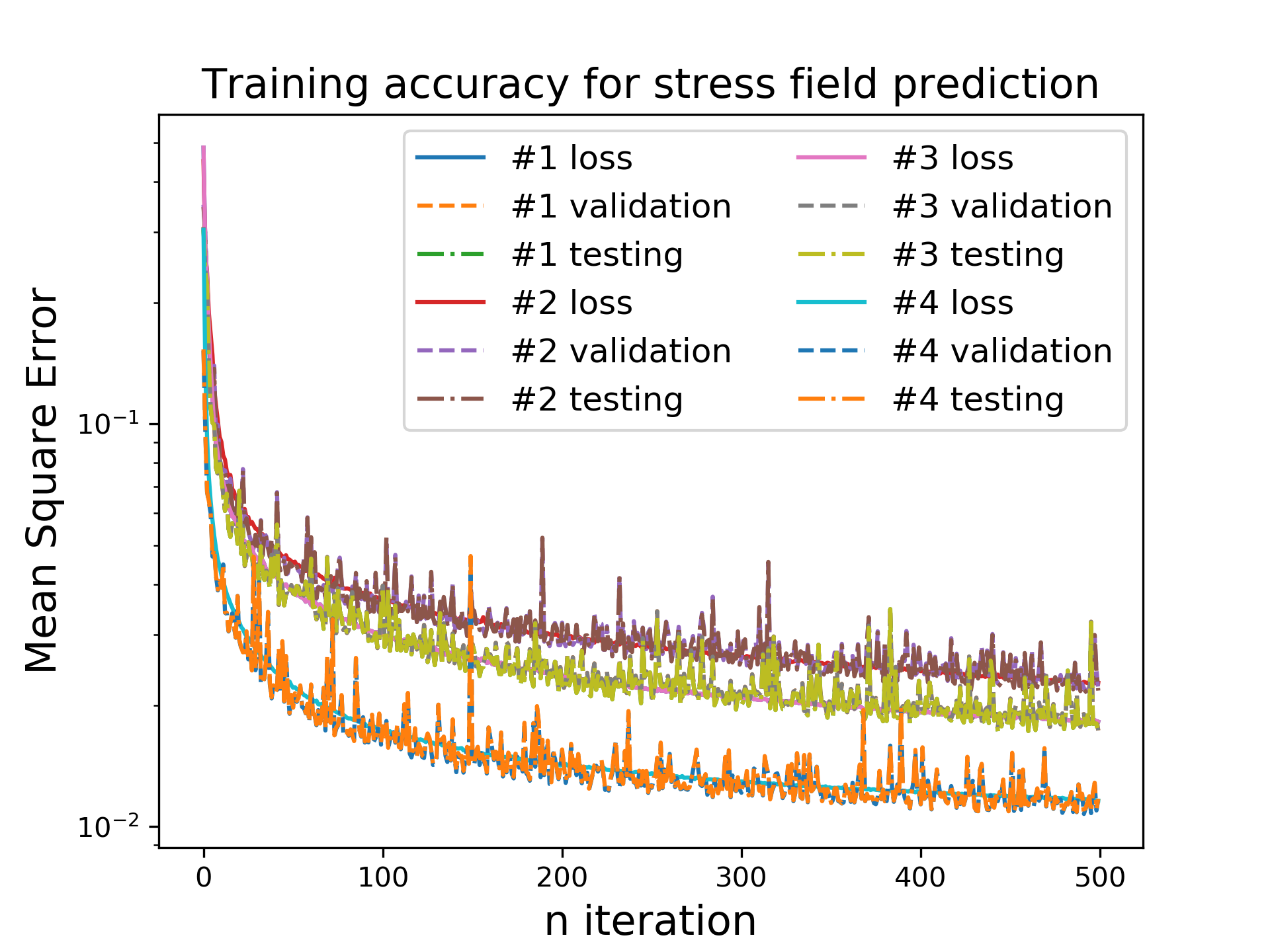}
  \caption{Stress prediction}
  \label{fig:StressLoss}
\end{subfigure}%
\begin{subfigure}{0.4\textwidth}
  \centering
  \includegraphics[width=\linewidth]{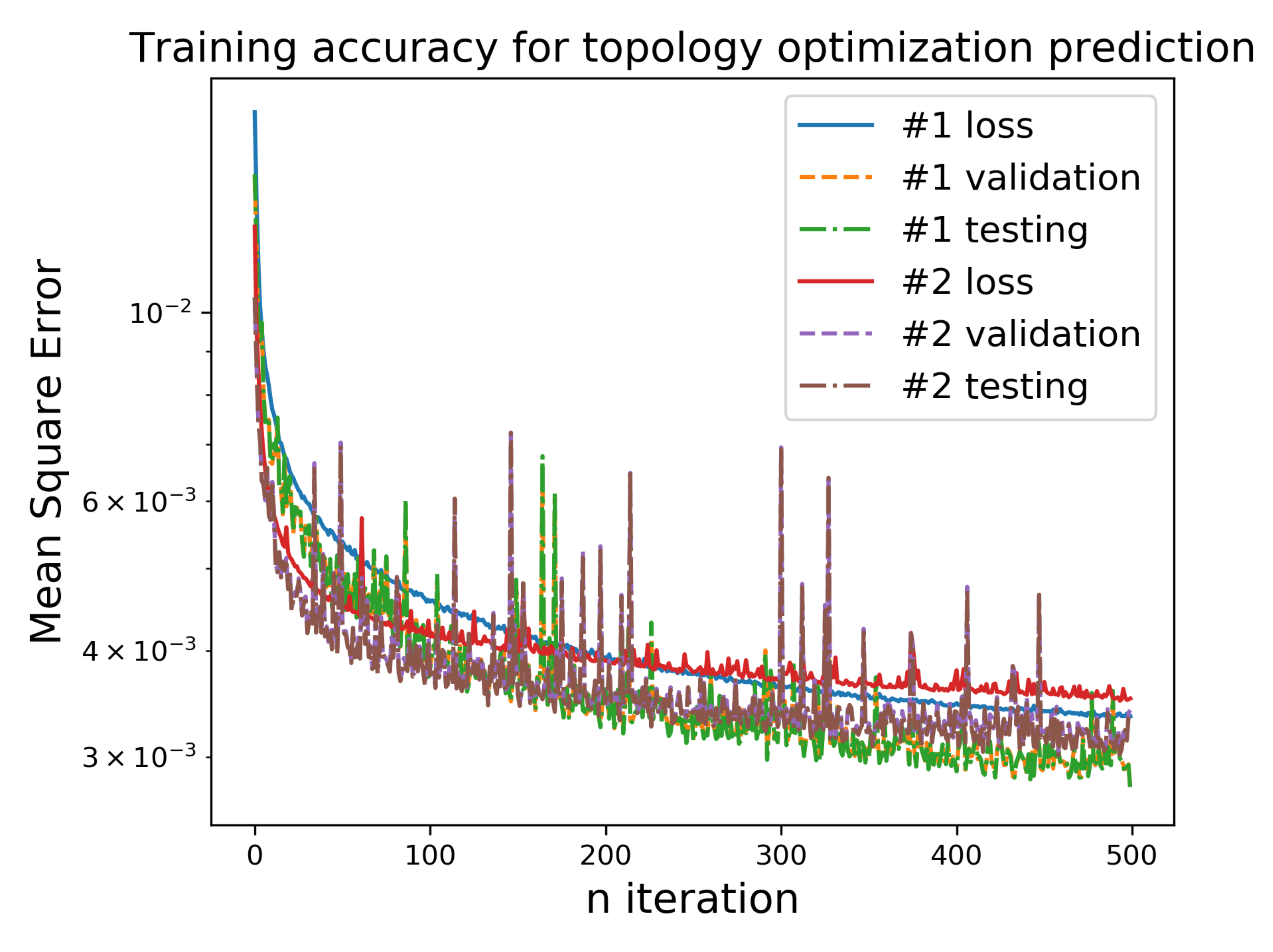}
  \caption{Topology optimization prediction}
  \label{fig:TopoLoss}
\end{subfigure}
\caption{Training loss, validation accuracy, and testing accuracy}
\label{fig:loss}
\end{figure}

We also consider the effect of the dataset size. Another training experiment with 180k simulation dataset ("Stress prediction (180k)" in Table~\ref{tab:2_1_Dataset}) is generated by the same 3-layer DeepGCN structure (model \#4 in Table~\ref{tab:3_3_stress}). The testing accuracy increases to 0.011706 which shows a better training performance, but only slightly improves the MAPE that reaches 2.41\%. This indicates that while the larger training dataset provides better results, a proper smaller dataset with around 45k data can still provide appropriate training results. Some of the predicted results from model \#4 in Table~\ref{tab:3_3_stress} are presented in Fig.~\ref{fig:StressResults} which illustrates the effectiveness of the BOGE approach. The average prediction time for the GNN model (model \#4 in Table~\ref{tab:3_3_stress}) with GPU is 0.012ms which is far less than the ABAQUS computation time (around 11.5s in Table~\ref{tab:2_1_Dataset}). Though the computation time for the GNN model does not include the time for data loading, input checking, and other unknown processes running in ABAQUS, the large difference of the order of magnitude validates the efficiency of our GNN surrogate model. The training loss and the validation/testing accuracy for the stress predictions is shown in Fig.~\ref{fig:loss}(\subref{fig:StressLoss})

\def \figWidth{0.24} 
\def \linWidthRatio{0.7} 

\begin{figure}[!htbp]
\centering
\includegraphics[width=\linWidthRatio\linewidth]{fig/3_0_Legend.png}
\newline
\begin{subfigure}{\figWidth\textwidth}
  \centering
  \includegraphics[width=\linWidthRatio\linewidth]{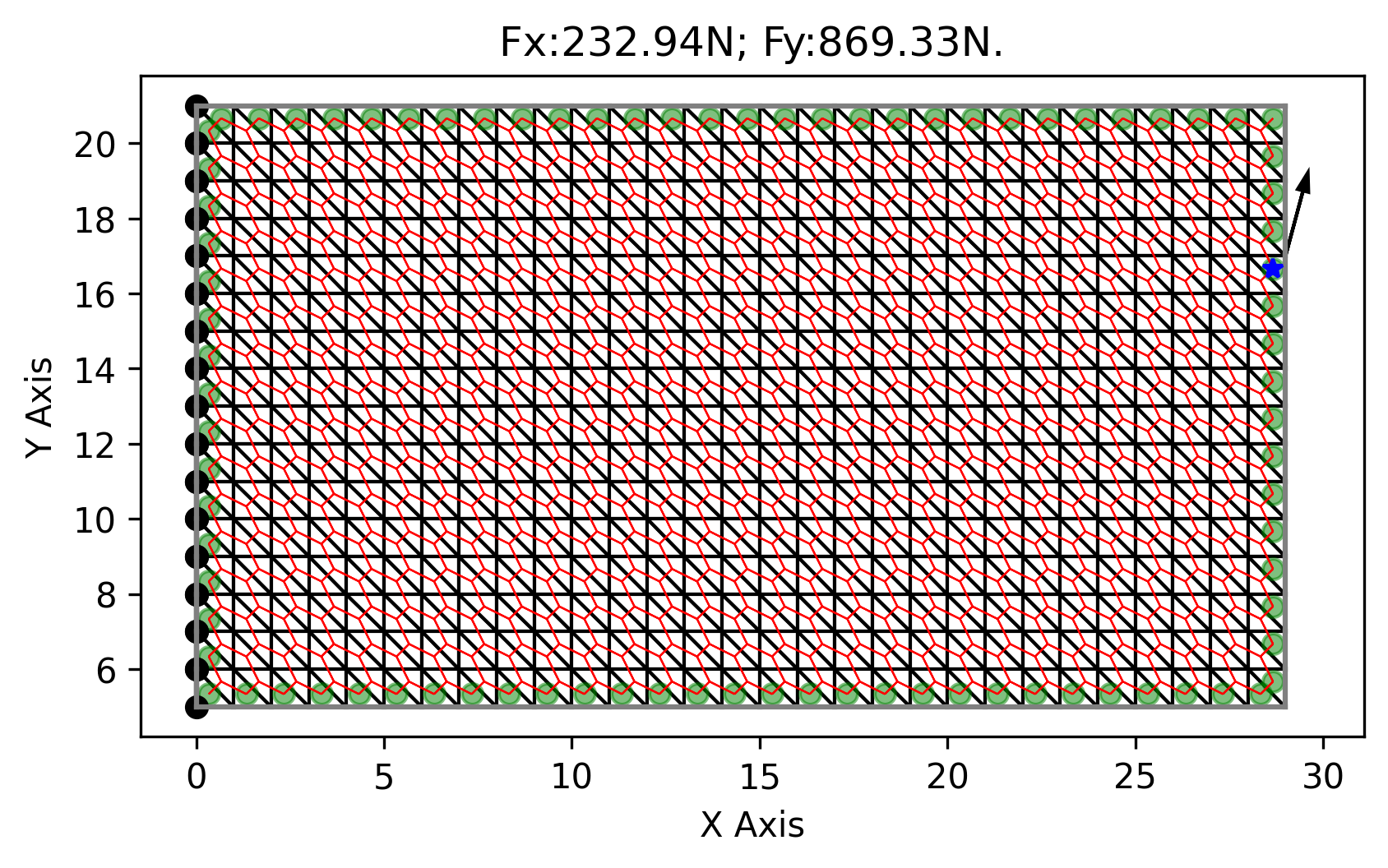}
  \caption{Shape \#1 simulation settings}
\end{subfigure}%
\begin{subfigure}{\figWidth\textwidth}
  \centering
  \includegraphics[width=\linWidthRatio\linewidth]{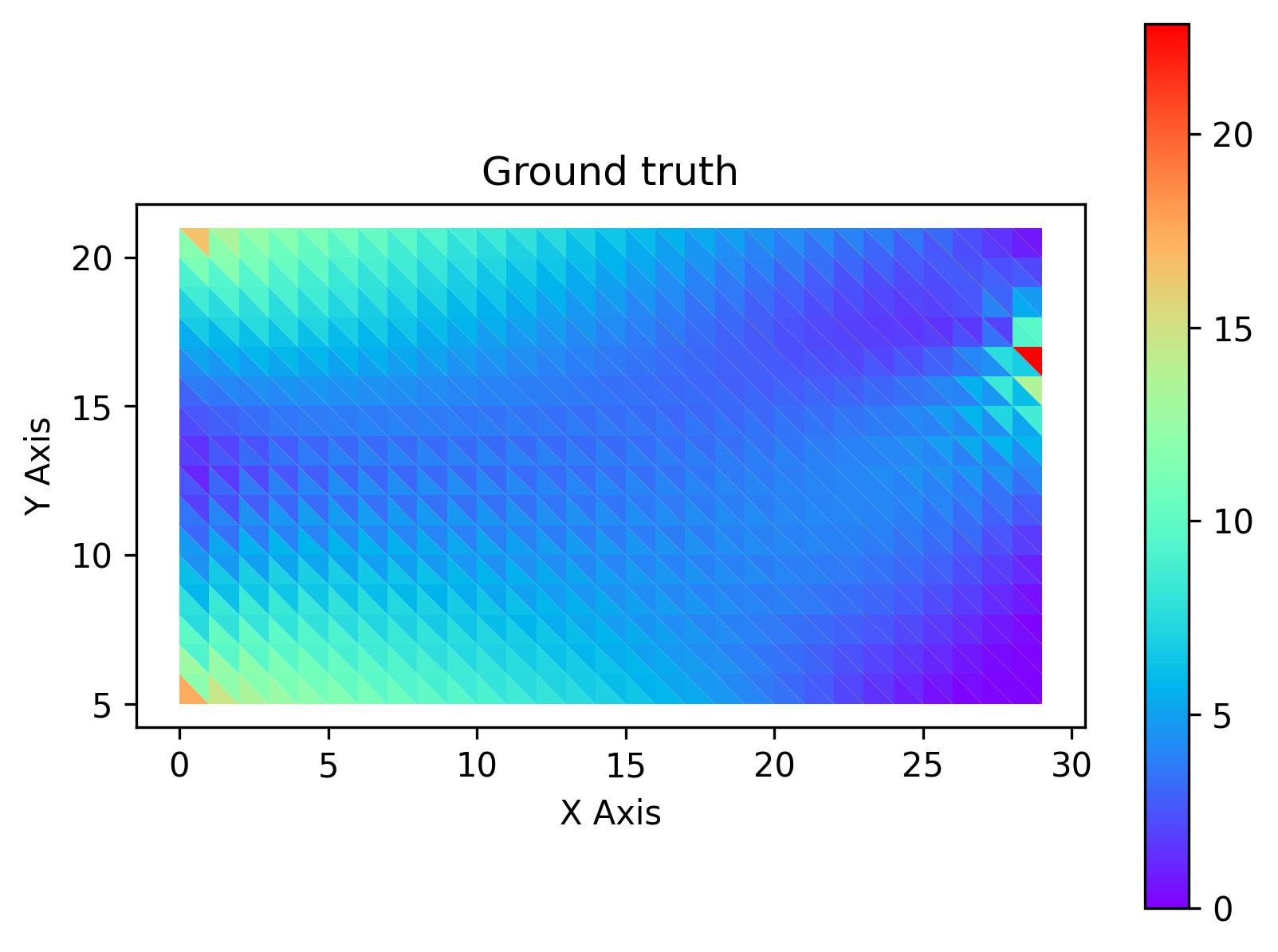}
  \caption{Shape \#1 ground truth}
\end{subfigure}
\begin{subfigure}{\figWidth\textwidth}
  \centering
  \includegraphics[width=\linWidthRatio\linewidth]{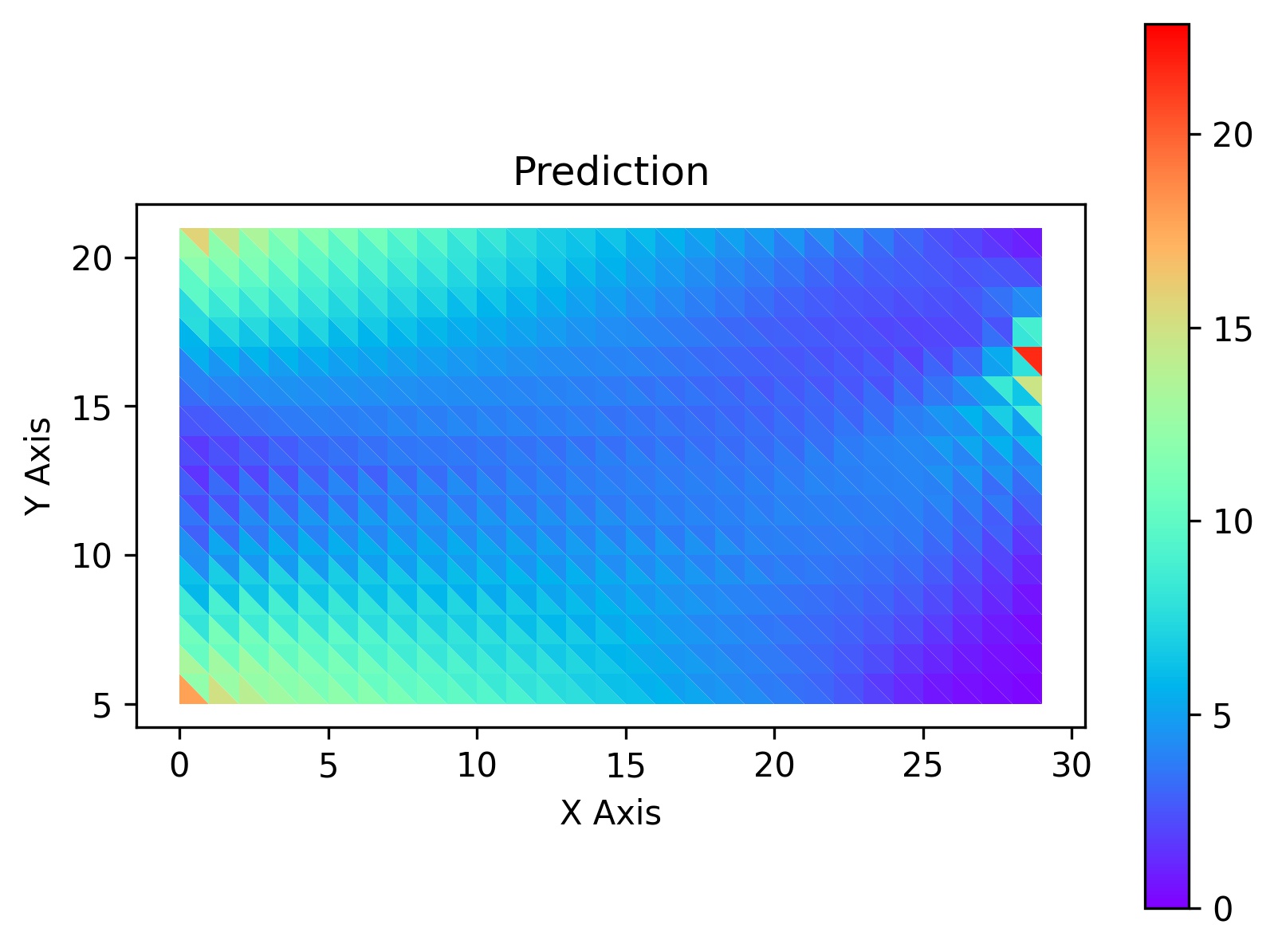}
  \caption{Shape \#1 prediction}
\end{subfigure}
\begin{subfigure}{\figWidth\textwidth}
  \centering
  \includegraphics[width=\linWidthRatio\linewidth]{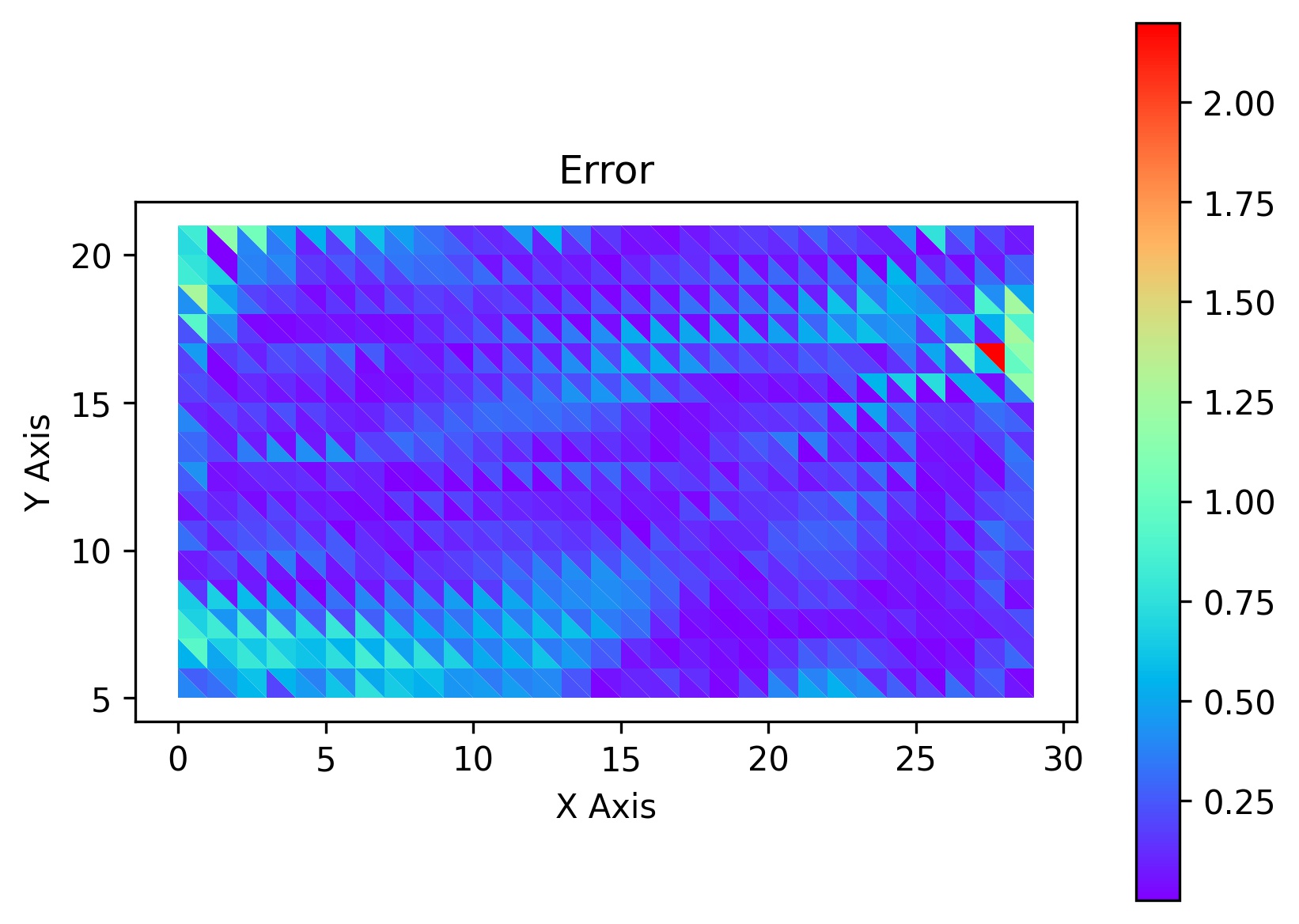}
  \caption{Shape \#1 error}
\end{subfigure}
\end{figure}

\begin{figure}[!htbp]\ContinuedFloat
\begin{subfigure}{\figWidth\textwidth}
  \centering
  \includegraphics[width=\linWidthRatio\linewidth]{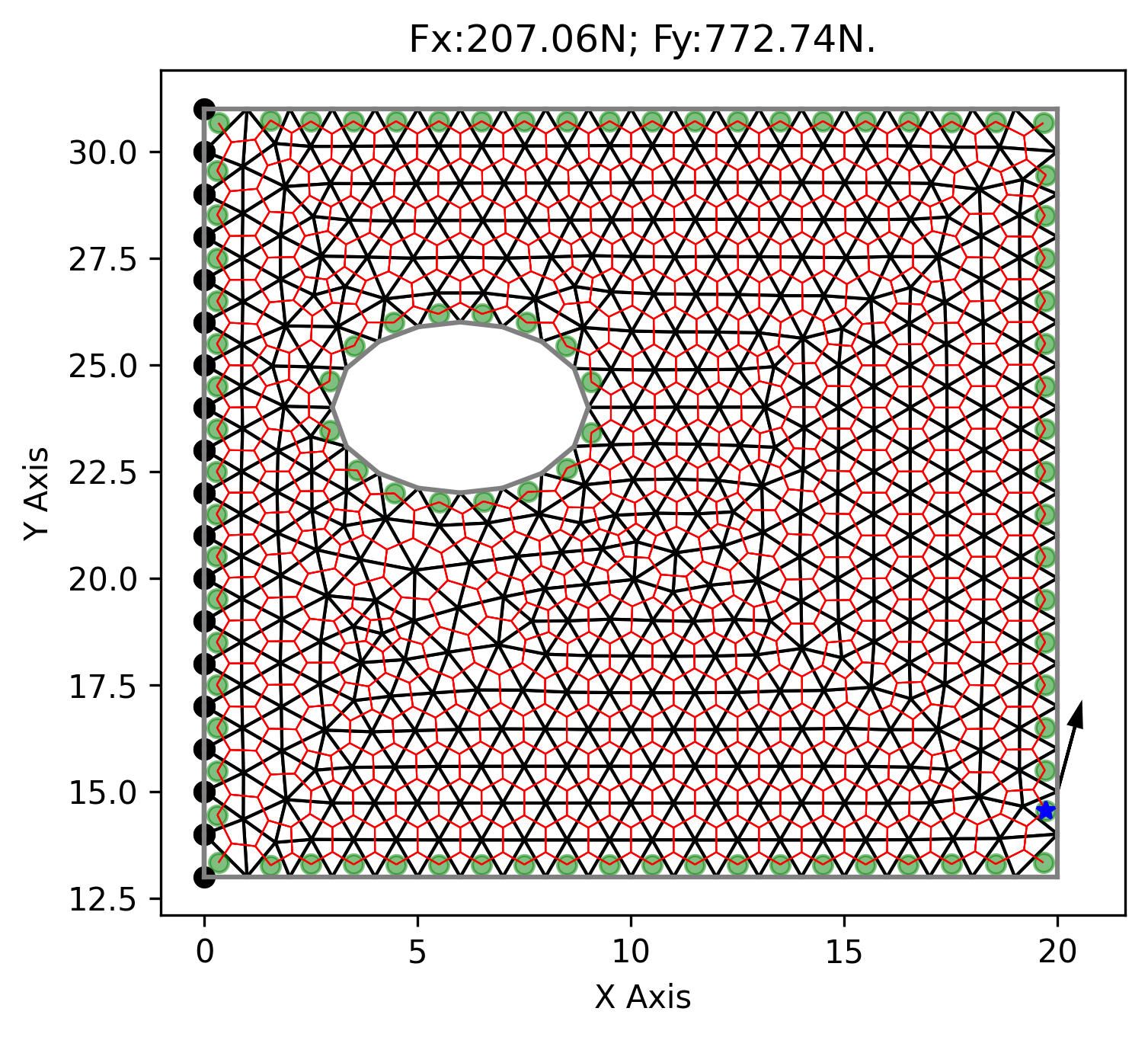}
  \caption{Shape \#2 simulation settings}
\end{subfigure}%
\begin{subfigure}{\figWidth\textwidth}
  \centering
  \includegraphics[width=\linWidthRatio\linewidth]{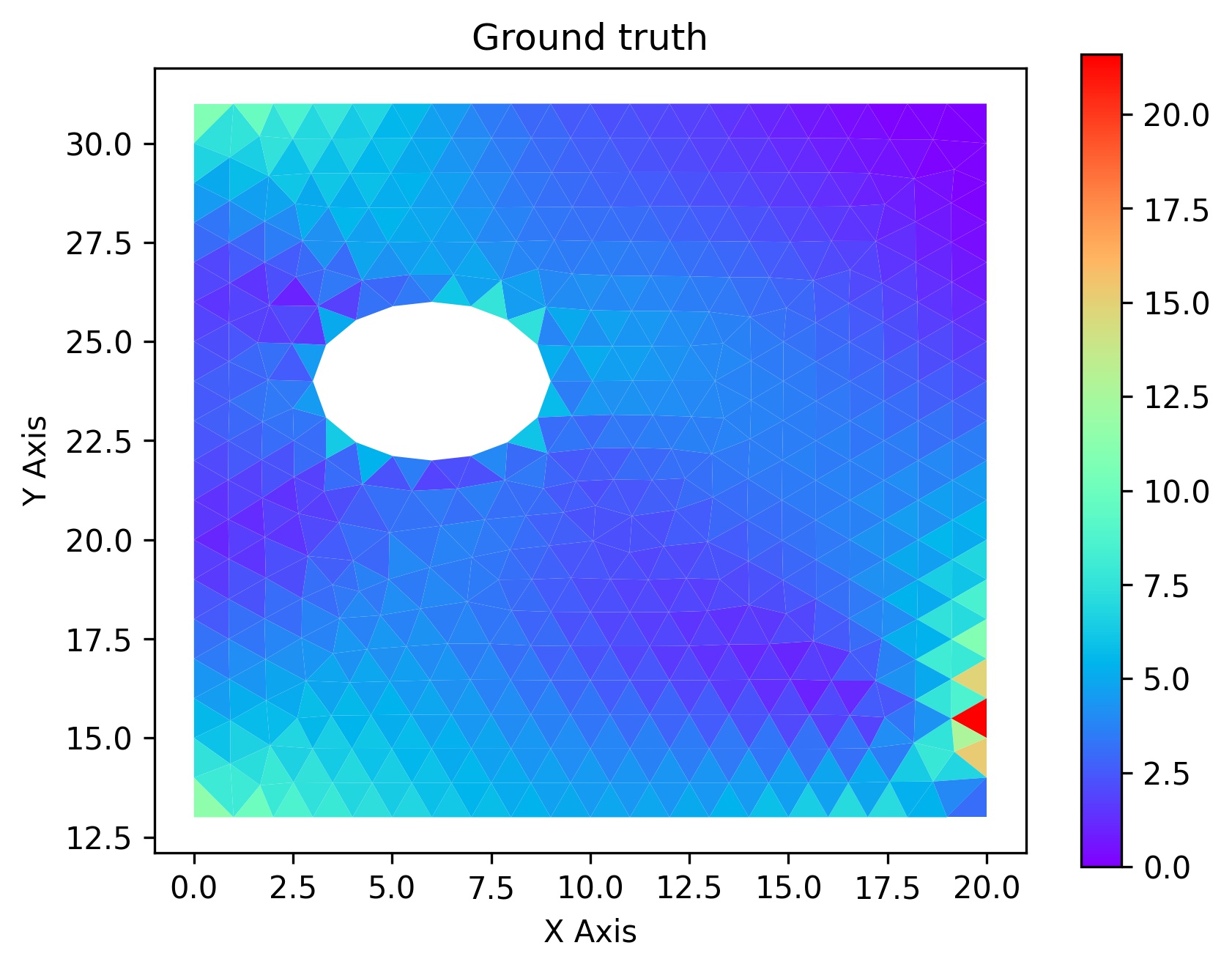}
  \caption{Shape \#2 ground truth}
\end{subfigure}
\begin{subfigure}{\figWidth\textwidth}
  \centering
  \includegraphics[width=\linWidthRatio\linewidth]{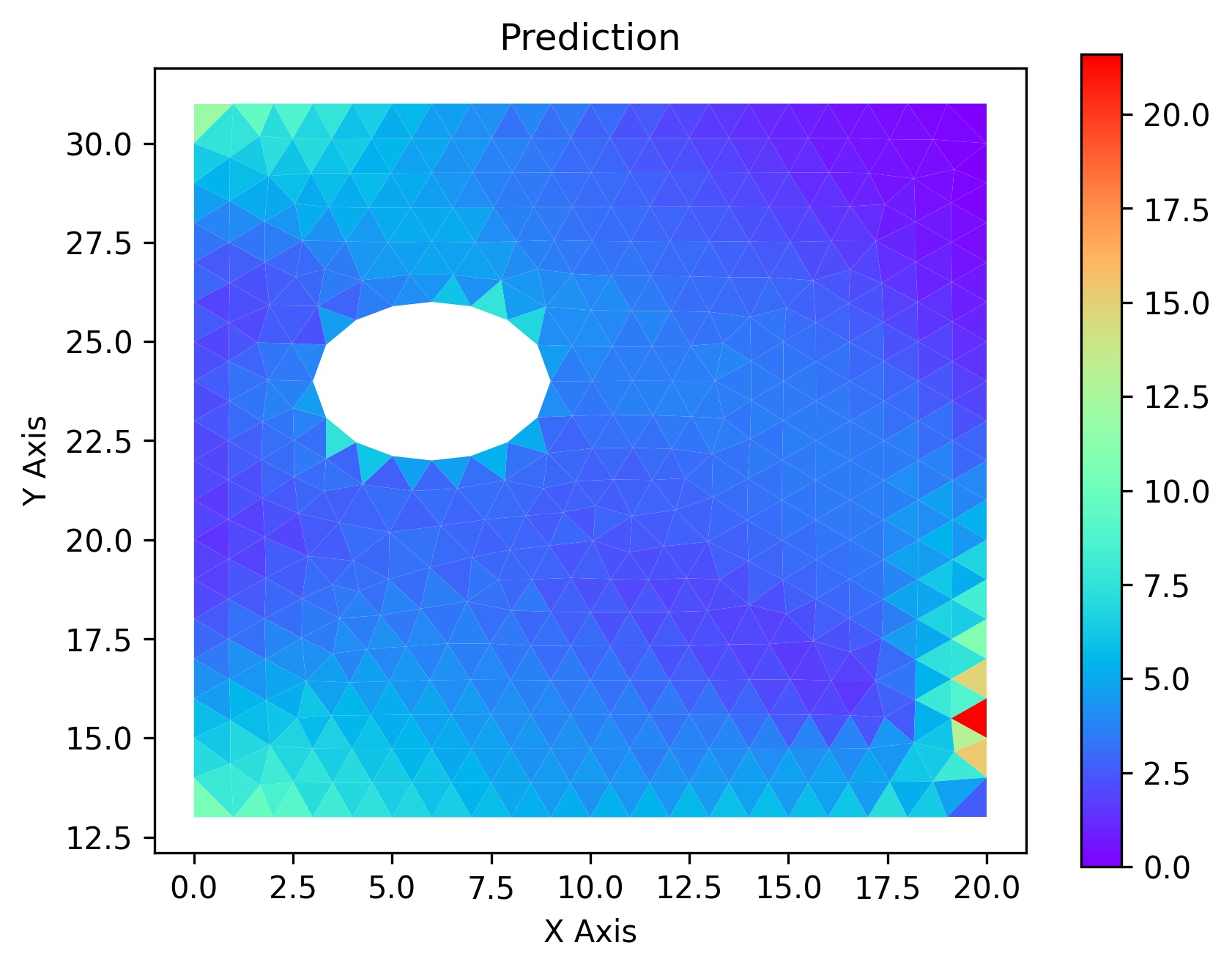}
  \caption{Shape \#2 prediction}
\end{subfigure}
\begin{subfigure}{\figWidth\textwidth}
  \centering
  \includegraphics[width=\linWidthRatio\linewidth]{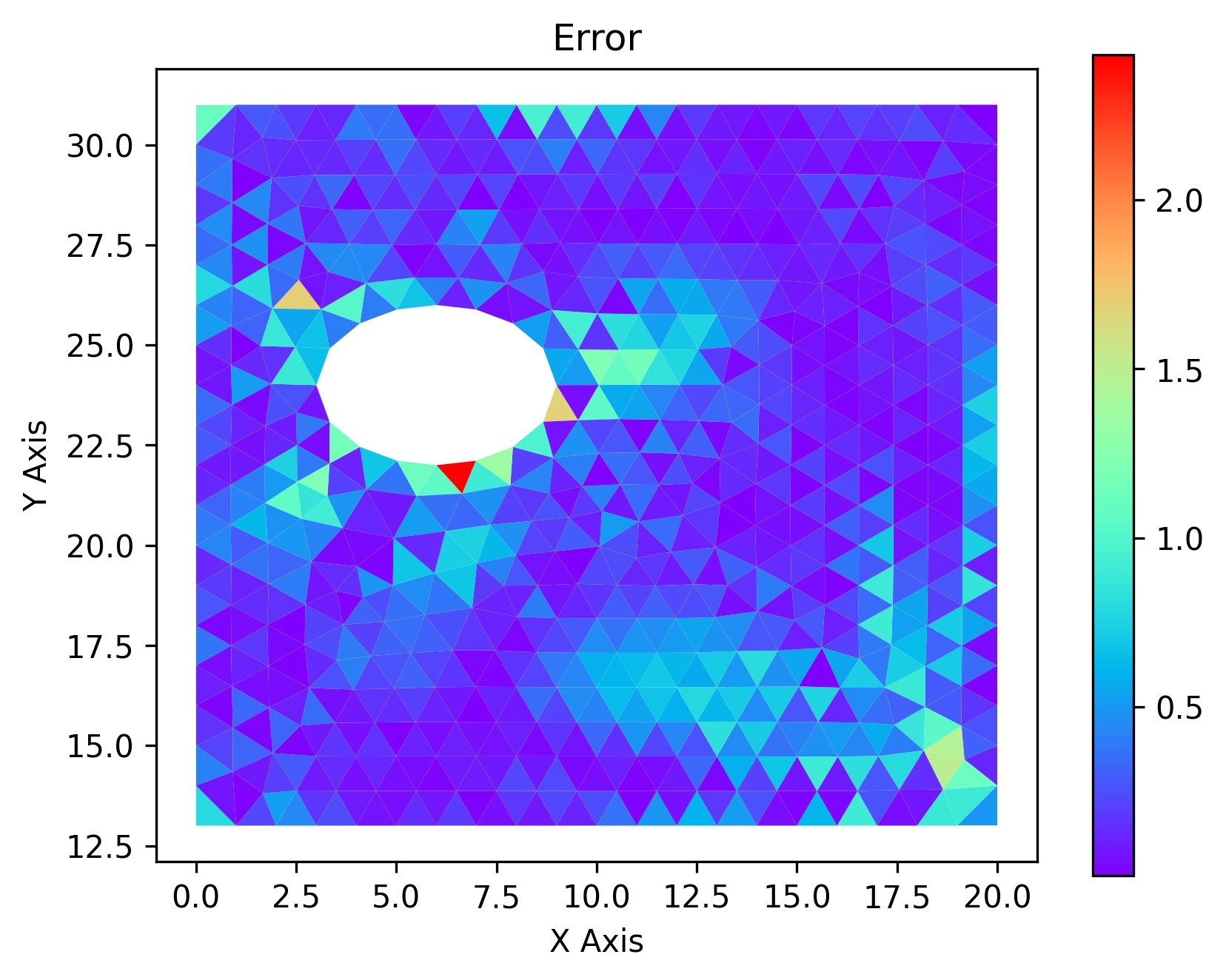}
  \caption{Shape \#2 error}
\end{subfigure}
\end{figure}

\begin{figure}[!htbp]\ContinuedFloat
\begin{subfigure}{\figWidth\textwidth}
  \centering
  \includegraphics[width=\linWidthRatio\linewidth]{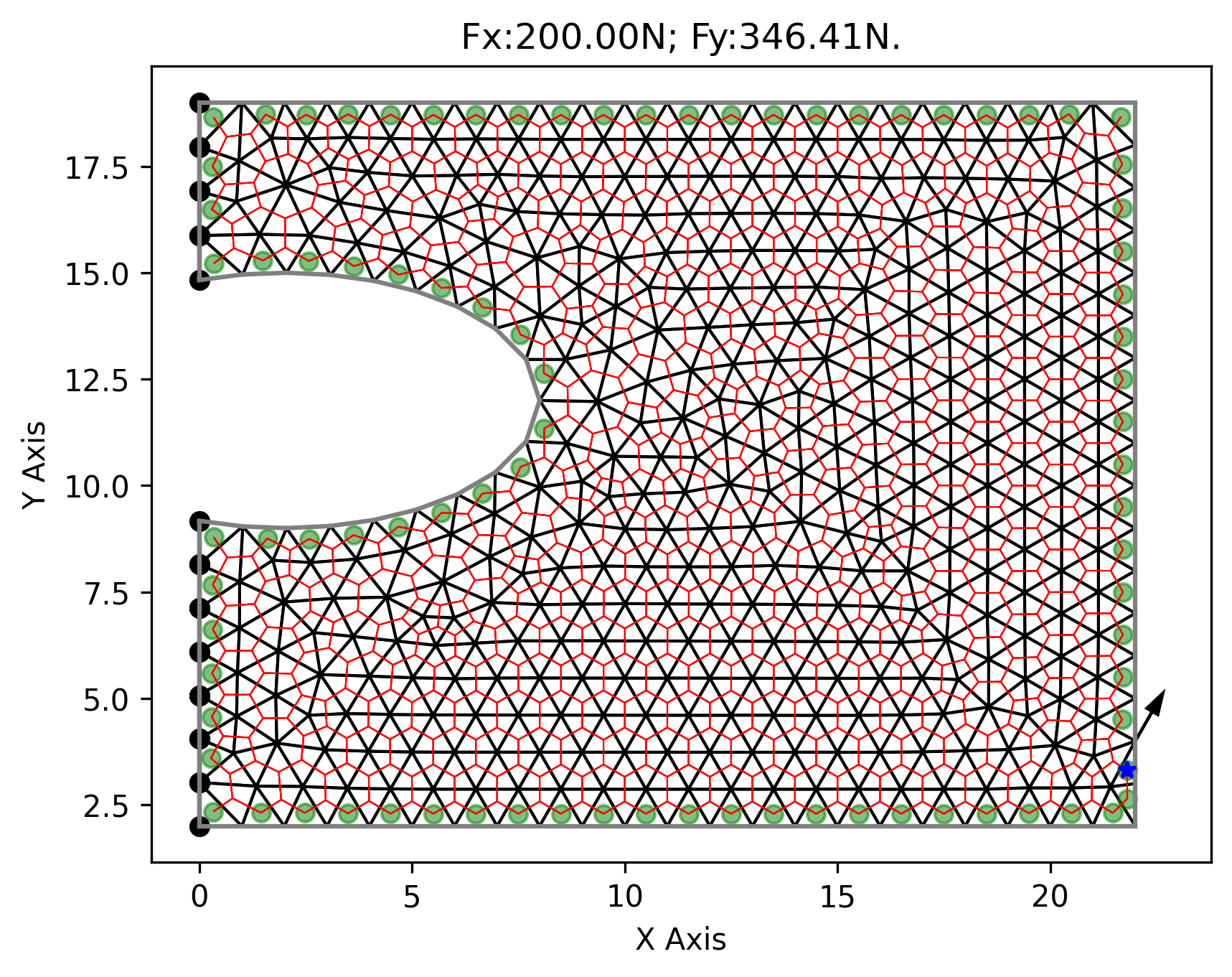}
  \caption{Shape \#3 simulation settings}
\end{subfigure}%
\begin{subfigure}{\figWidth\textwidth}
  \centering
  \includegraphics[width=\linWidthRatio\linewidth]{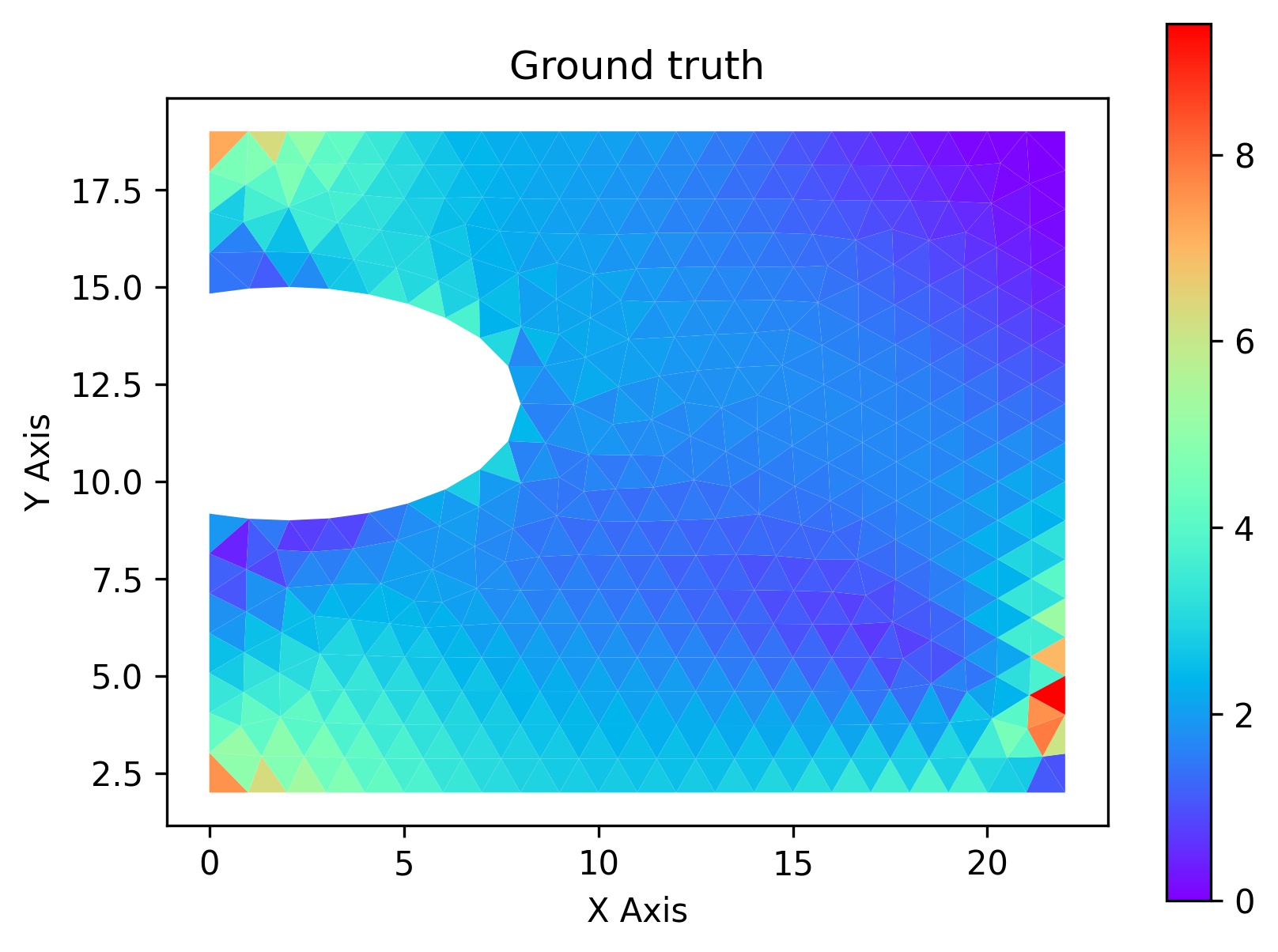}
  \caption{Shape \#3 ground truth}
\end{subfigure}
\begin{subfigure}{\figWidth\textwidth}
  \centering
  \includegraphics[width=\linWidthRatio\linewidth]{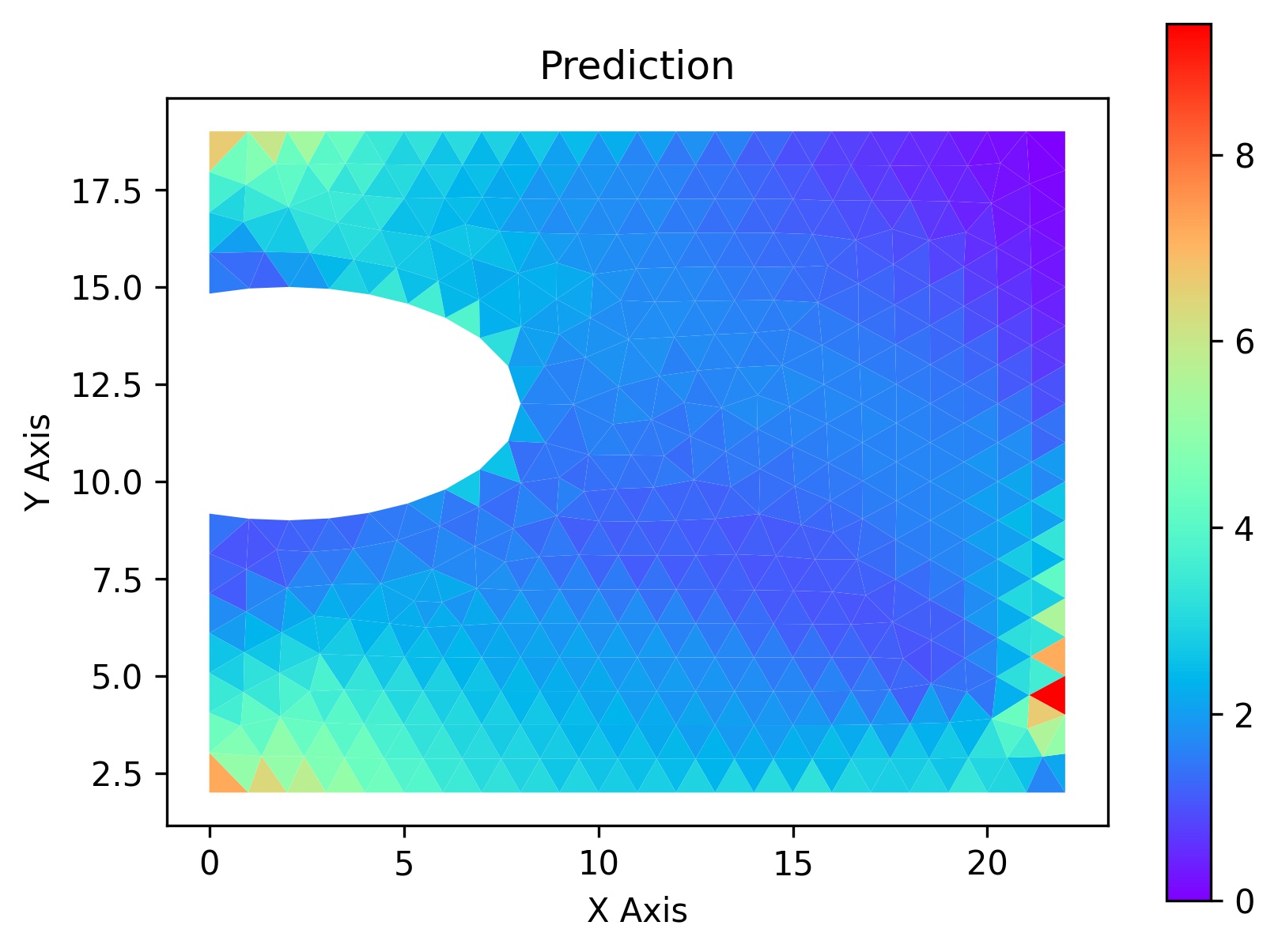}
  \caption{Shape \#3 prediction}
\end{subfigure}
\begin{subfigure}{\figWidth\textwidth}
  \centering
  \includegraphics[width=\linWidthRatio\linewidth]{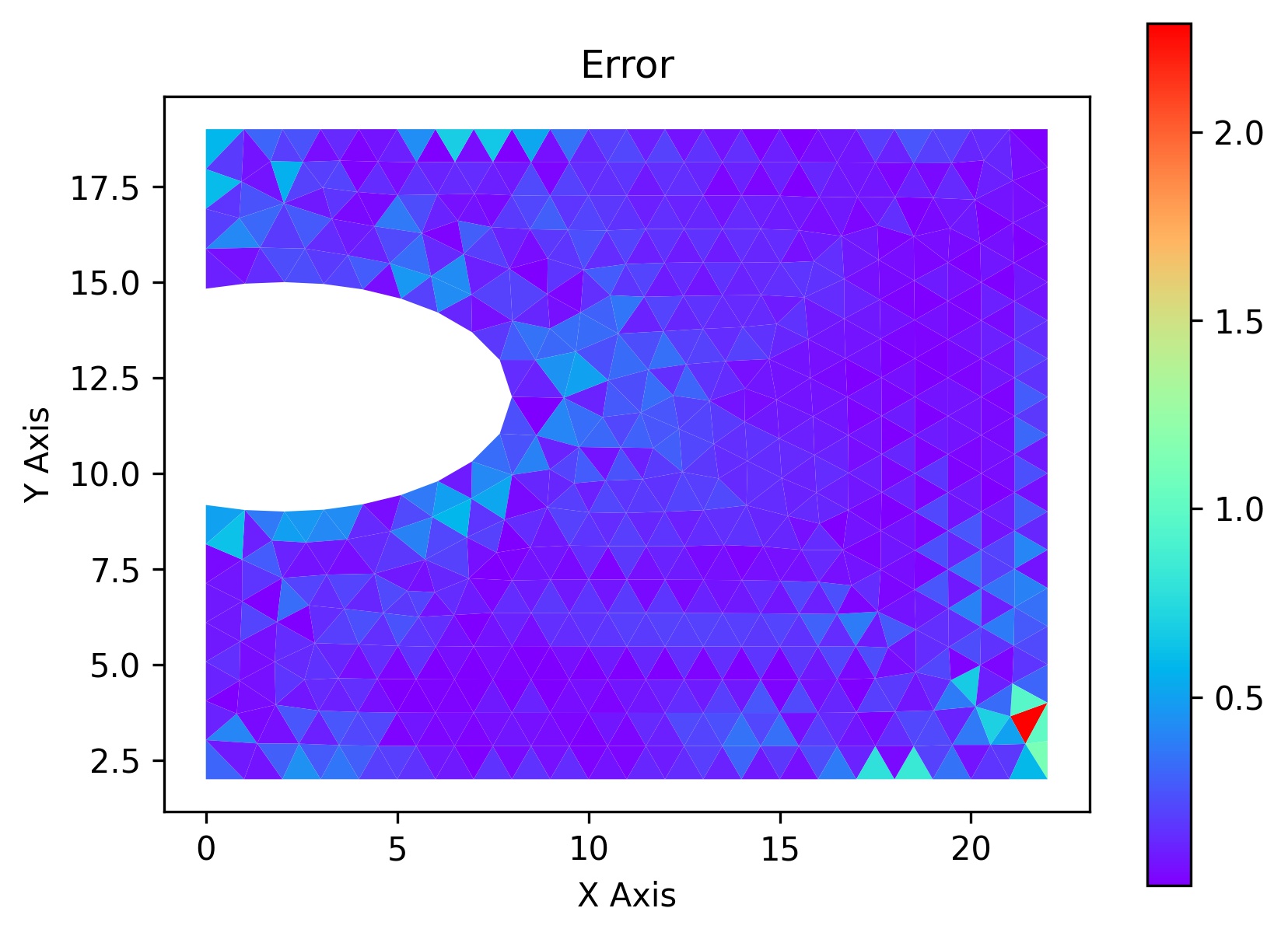}
  \caption{Shape \#3 error}
\end{subfigure}
\end{figure}

\begin{figure}[!htbp]\ContinuedFloat
\begin{subfigure}{\figWidth\textwidth}
  \centering
  \includegraphics[width=\linWidthRatio\linewidth]{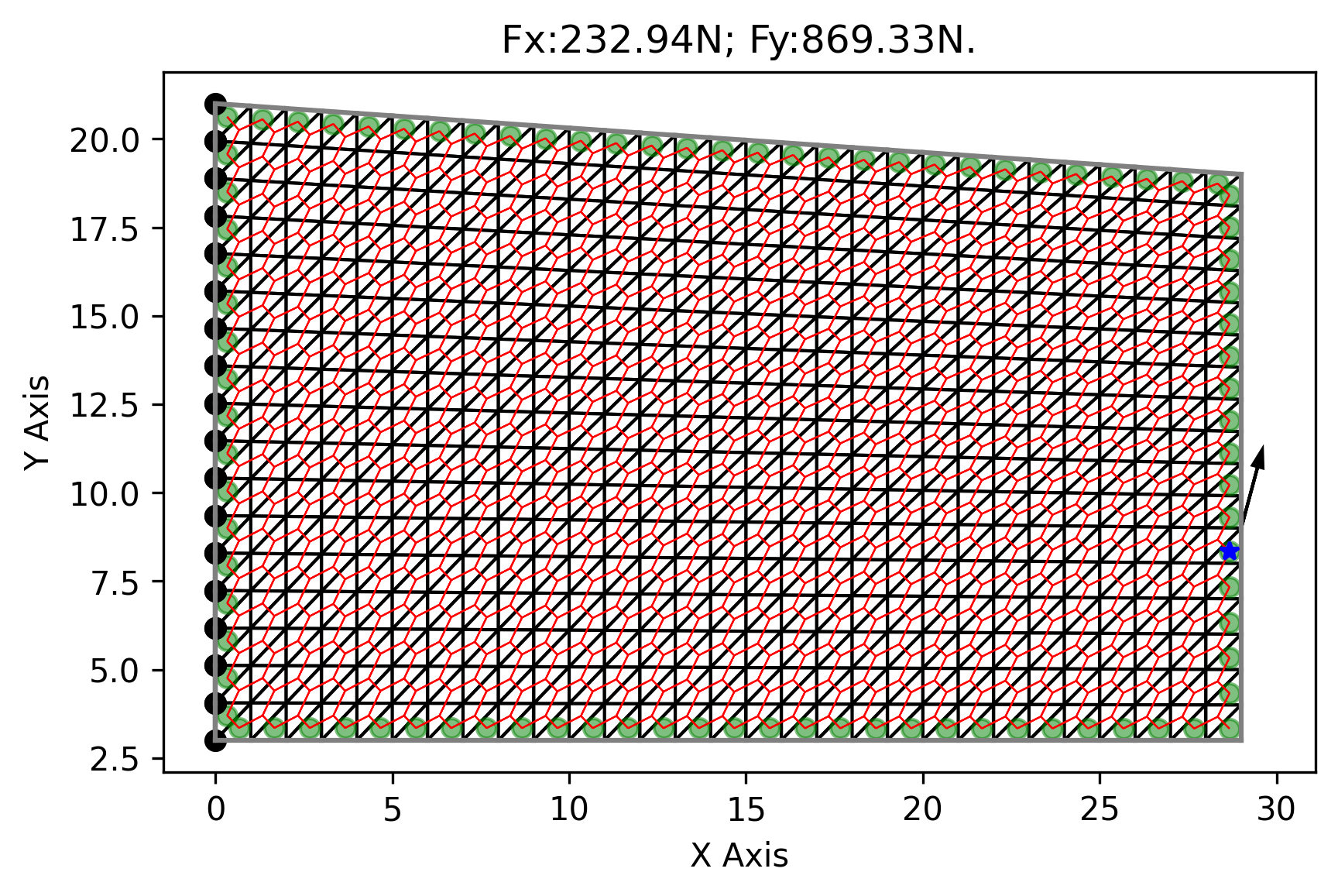}
  \caption{Shape \#4 simulation settings}
  \label{fig:ConvSim}
\end{subfigure}%
\begin{subfigure}{\figWidth\textwidth}
  \centering
  \includegraphics[width=\linWidthRatio\linewidth]{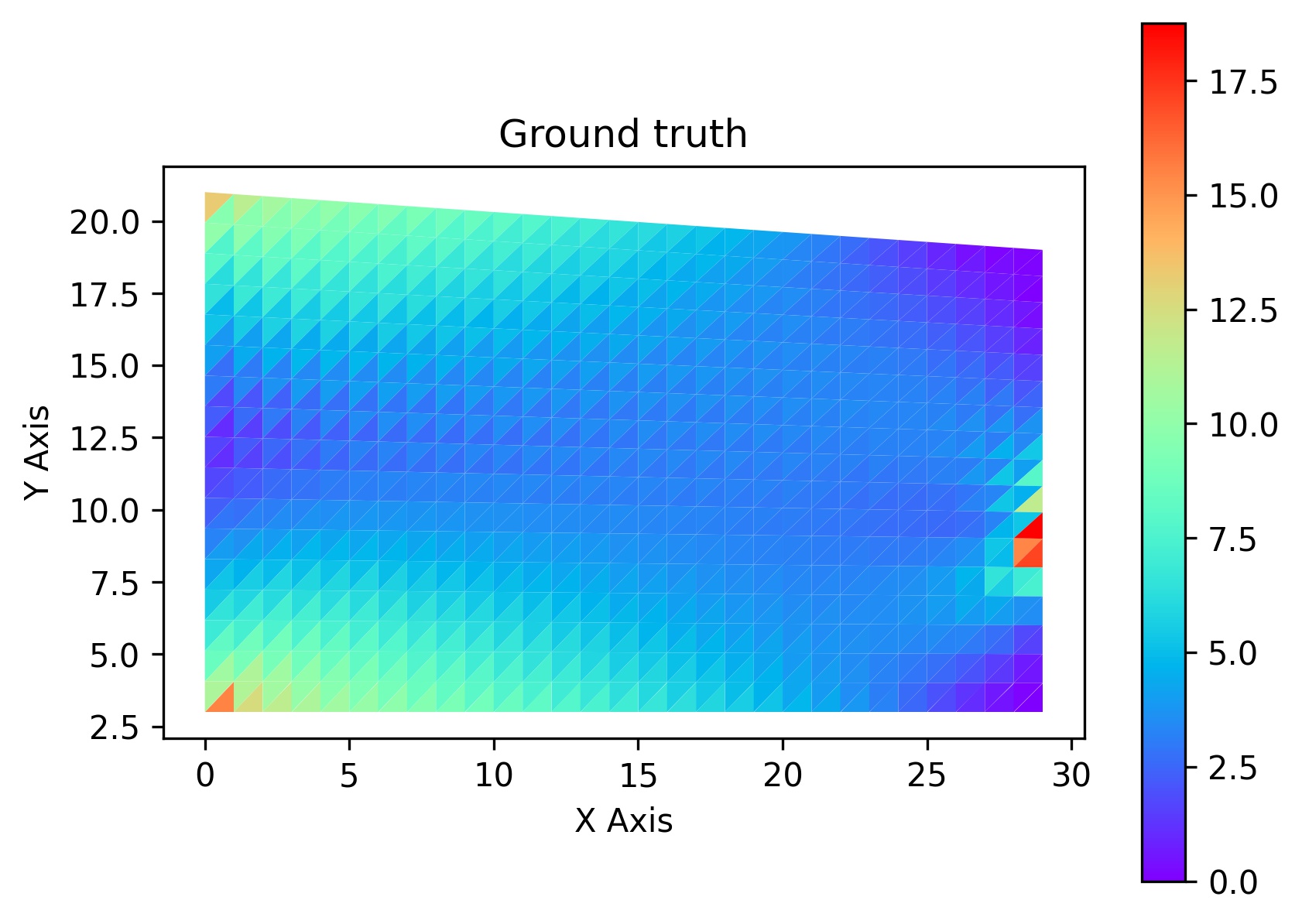}
  \caption{Shape \#4 ground truth}
\end{subfigure}
\begin{subfigure}{\figWidth\textwidth}
  \centering
  \includegraphics[width=\linWidthRatio\linewidth]{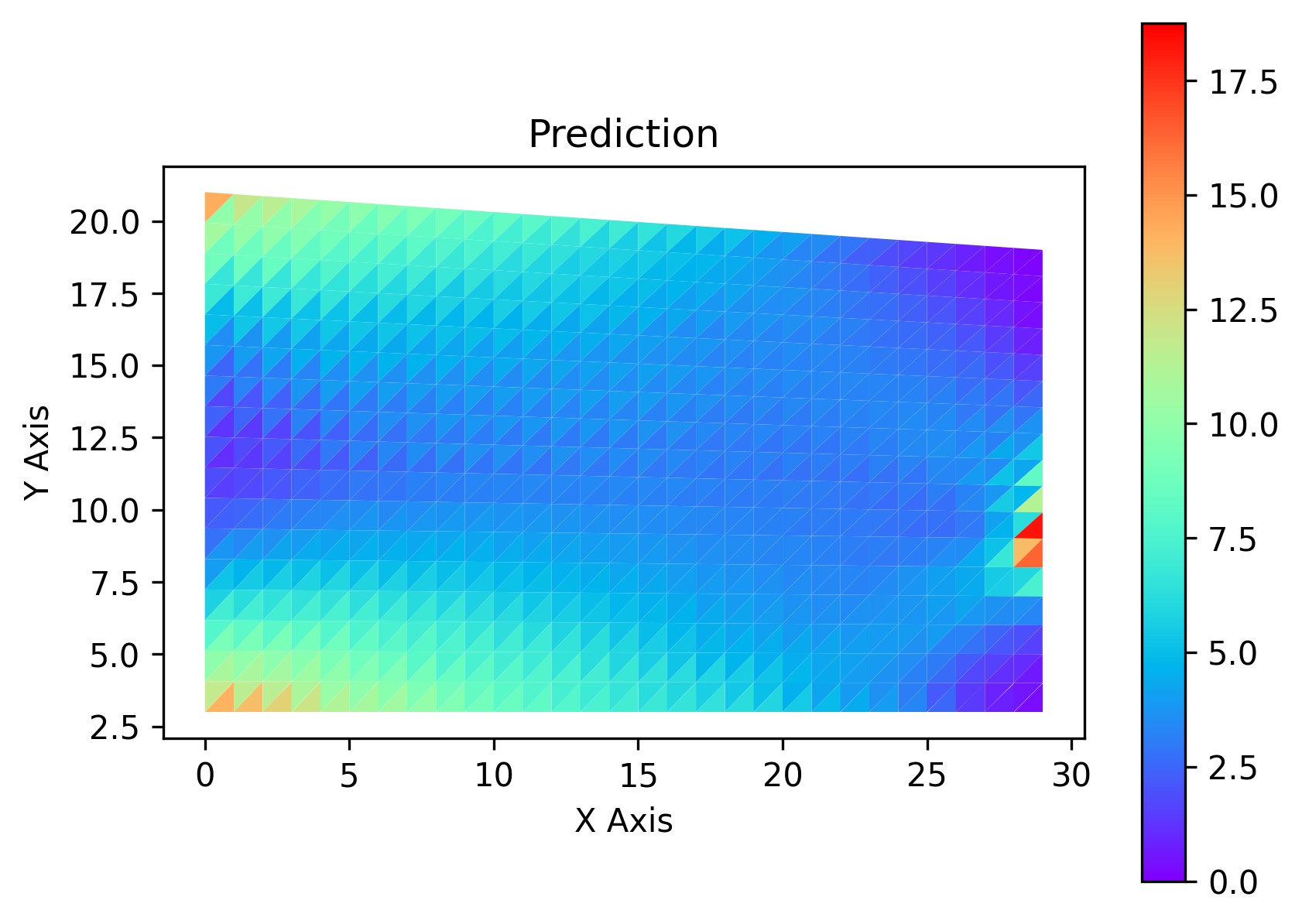}
  \caption{Shape \#4 prediction}
\end{subfigure}
\begin{subfigure}{\figWidth\textwidth}
  \centering
  \includegraphics[width=\linWidthRatio\linewidth]{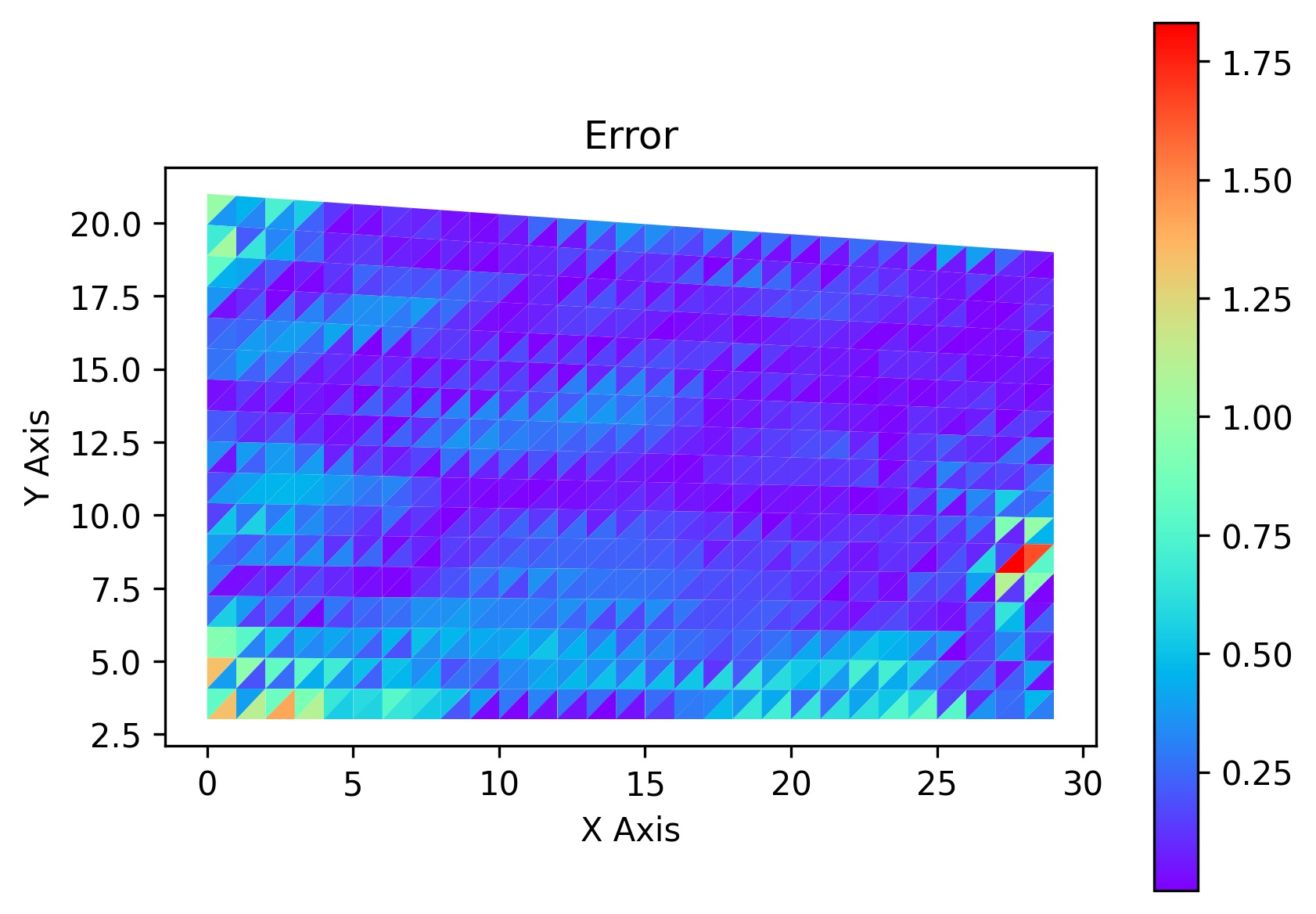}
  \caption{Shape \#4 error}
\end{subfigure}
\end{figure}

\begin{figure}[!htbp]\ContinuedFloat
\begin{subfigure}{\figWidth\textwidth}
  \centering
  \includegraphics[width=\linWidthRatio\linewidth]{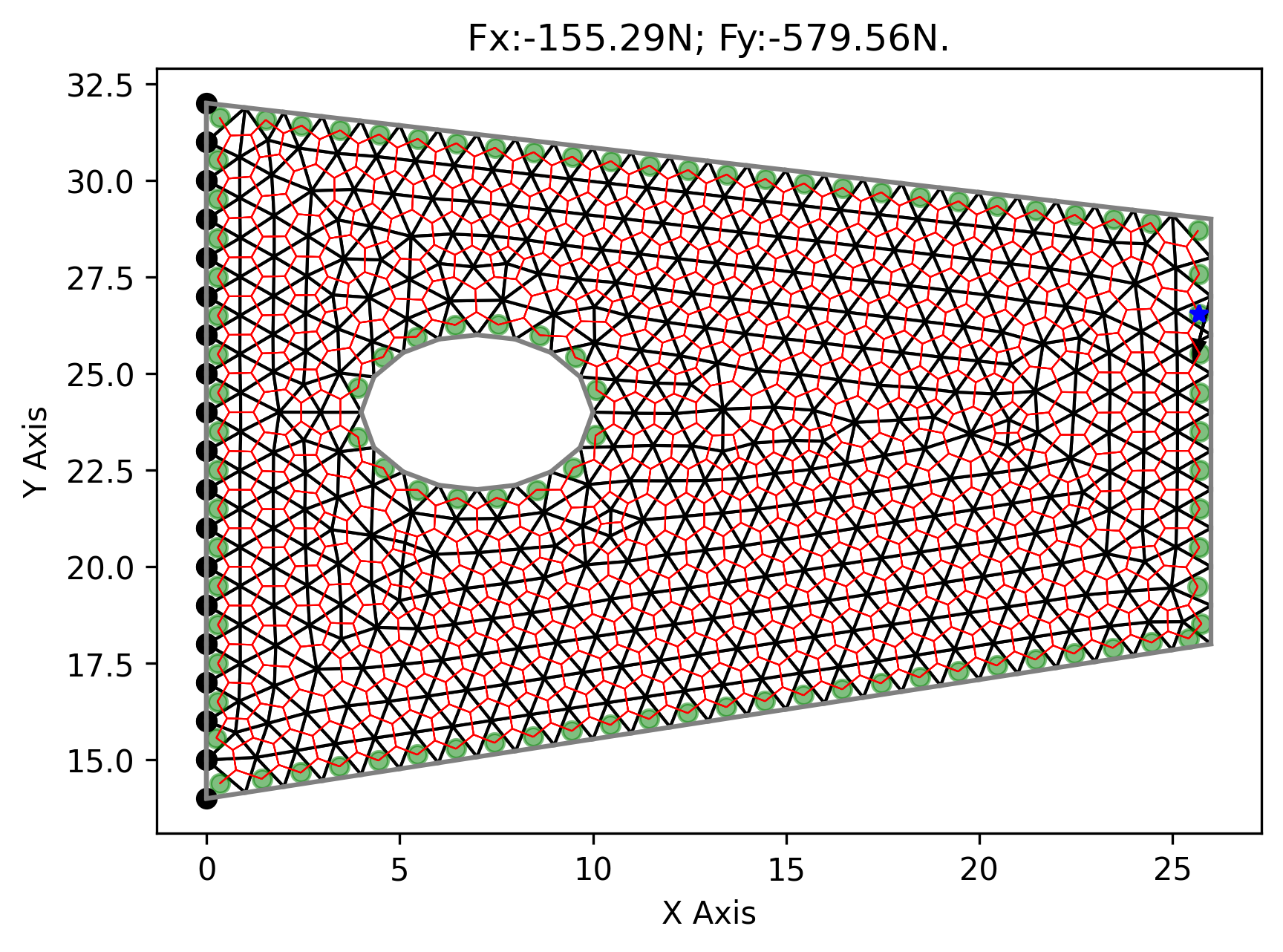}
  \caption{Shape \#5 simulation settings}
\end{subfigure}%
\begin{subfigure}{\figWidth\textwidth}
  \centering
  \includegraphics[width=\linWidthRatio\linewidth]{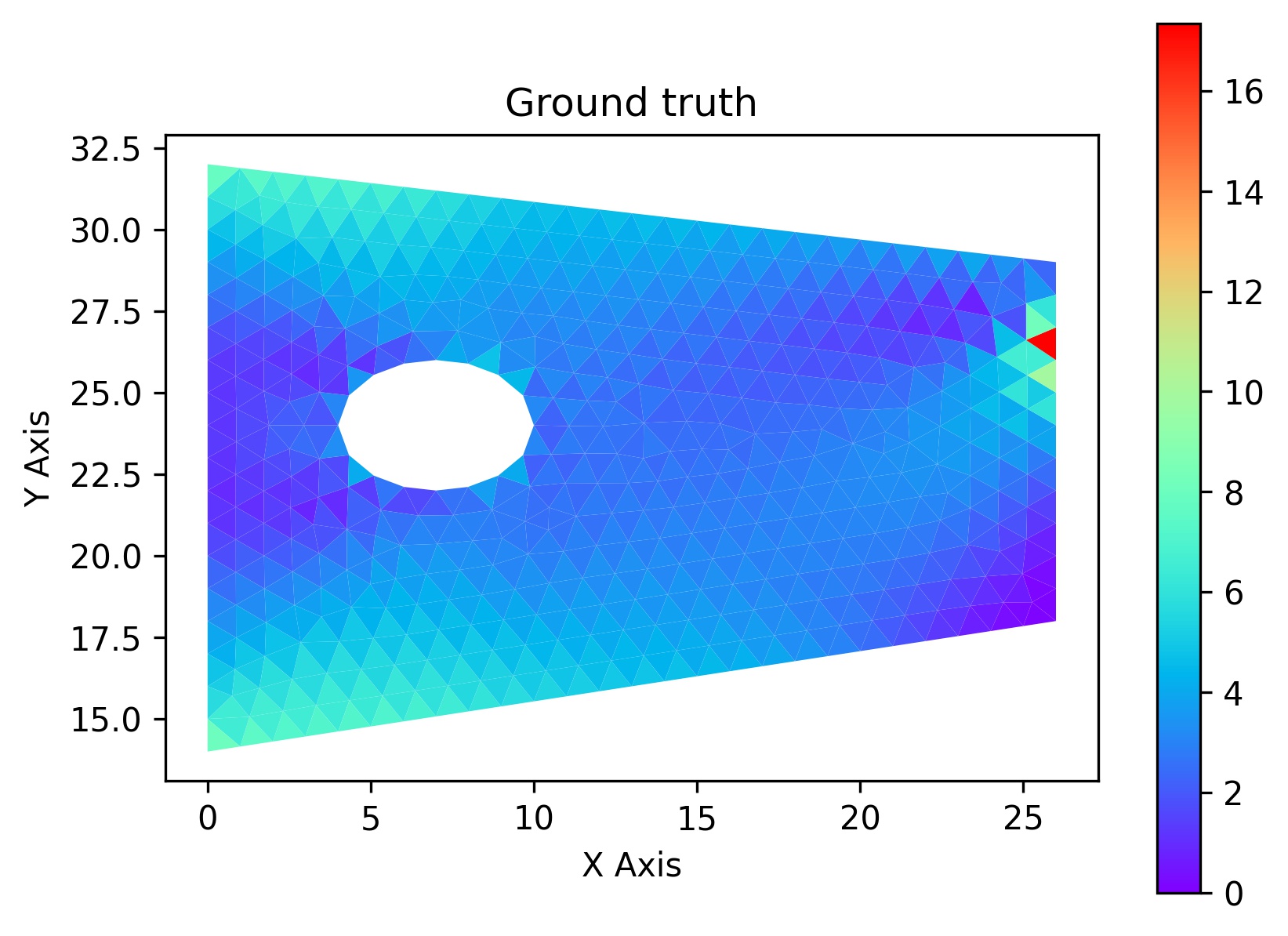}
  \caption{Shape \#5 ground truth}
\end{subfigure}
\begin{subfigure}{\figWidth\textwidth}
  \centering
  \includegraphics[width=\linWidthRatio\linewidth]{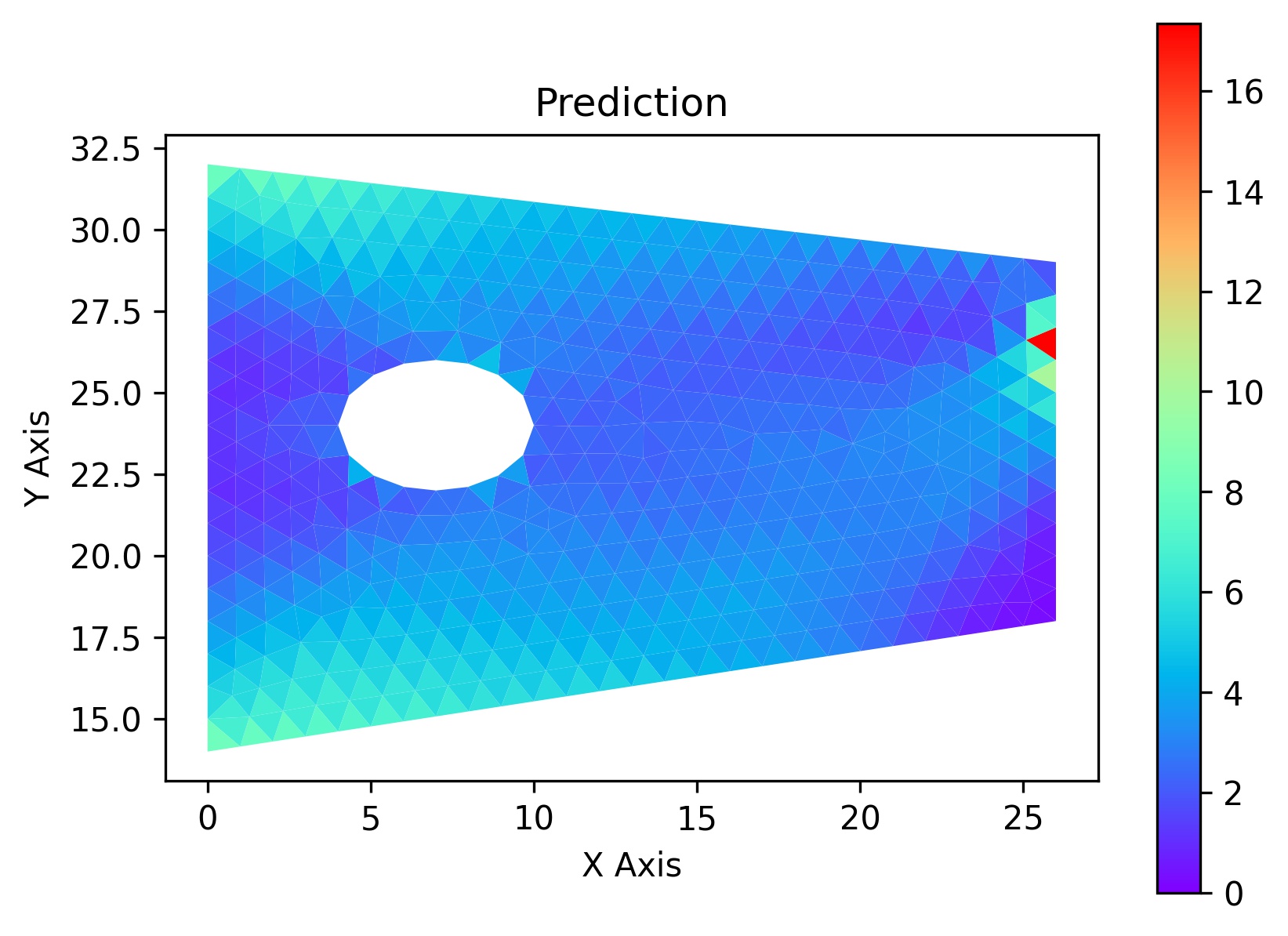}
  \caption{Shape \#5 prediction}
\end{subfigure}
\begin{subfigure}{\figWidth\textwidth}
  \centering
  \includegraphics[width=\linWidthRatio\linewidth]{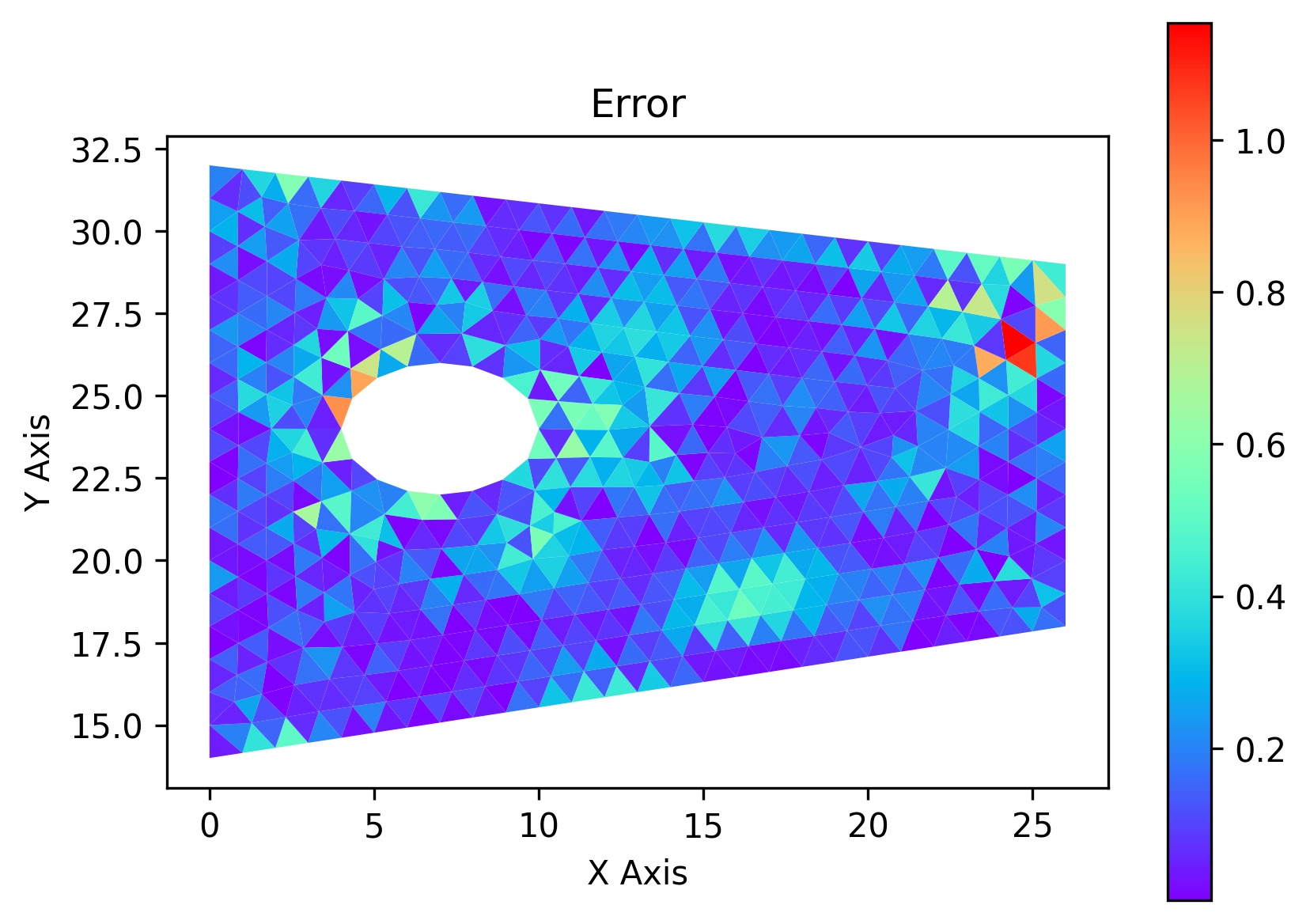}
  \caption{Shape \#5 error}
\end{subfigure}
\end{figure}

\begin{figure}[!htbp]\ContinuedFloat
\begin{subfigure}{\figWidth\textwidth}
  \centering
  \includegraphics[width=\linWidthRatio\linewidth]{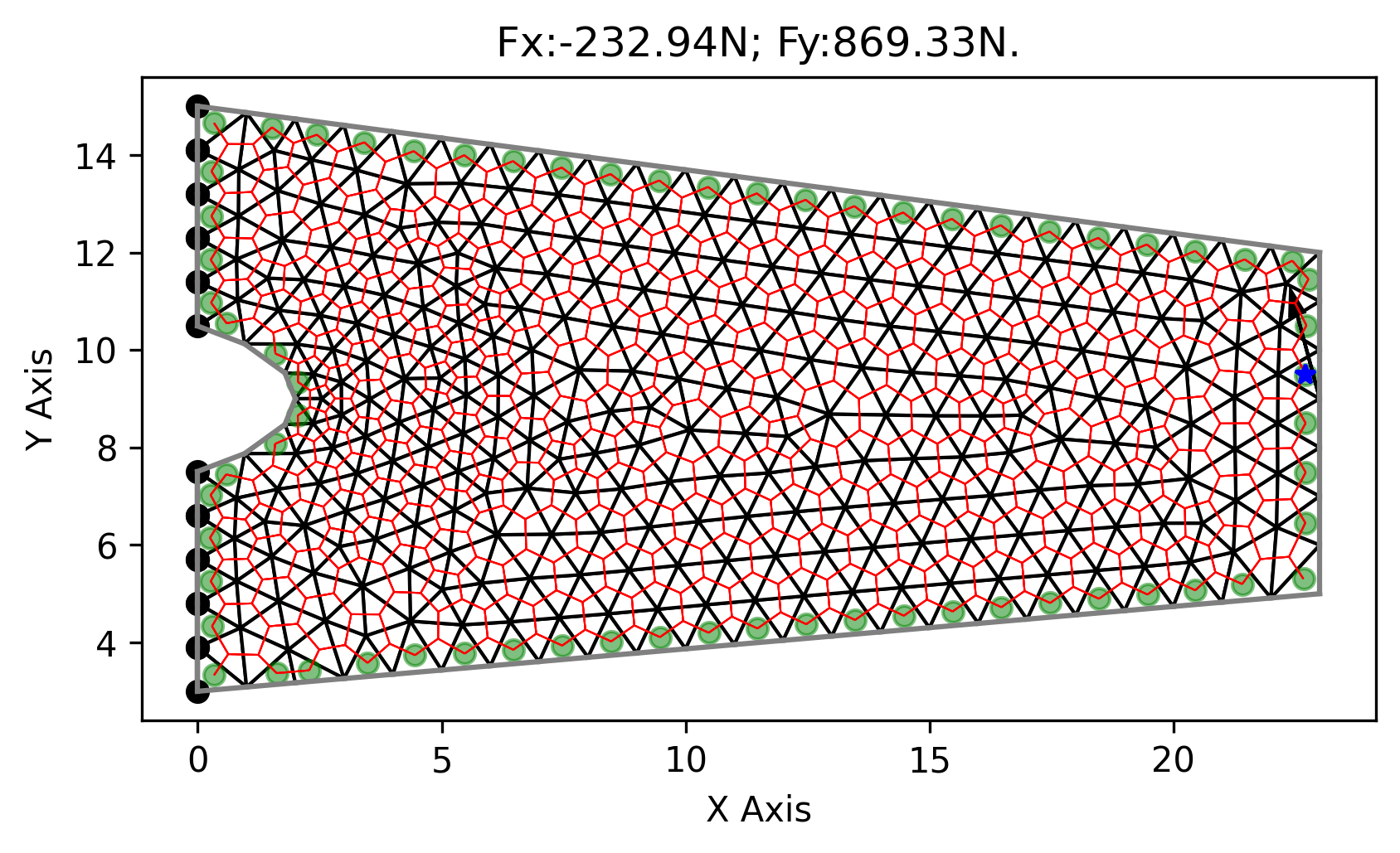}
  \caption{Shape \#6 simulation settings}
\end{subfigure}%
\begin{subfigure}{\figWidth\textwidth}
  \centering
  \includegraphics[width=\linWidthRatio\linewidth]{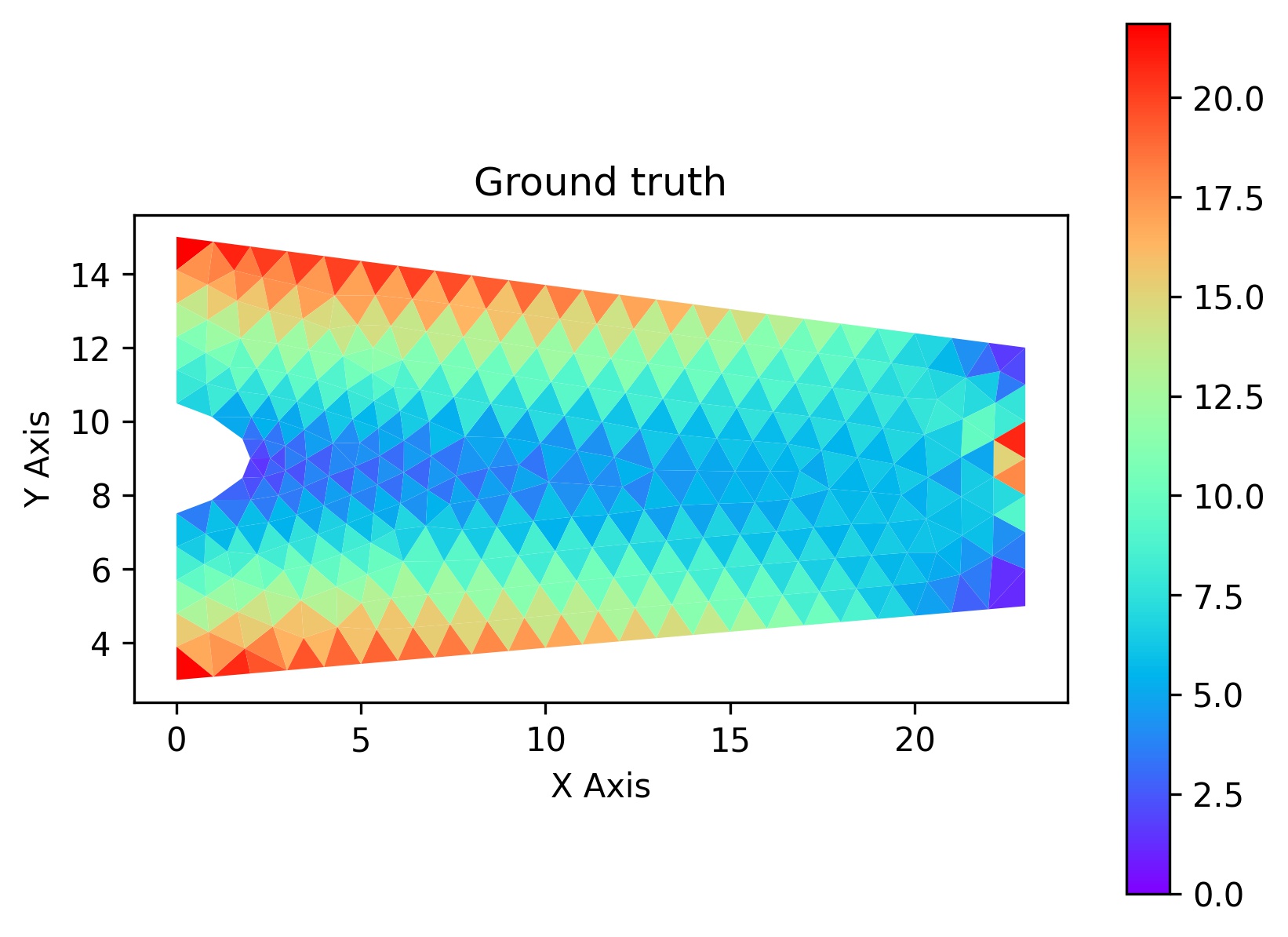}
  \caption{Shape \#6 ground truth}
\end{subfigure}
\begin{subfigure}{\figWidth\textwidth}
  \centering
  \includegraphics[width=\linWidthRatio\linewidth]{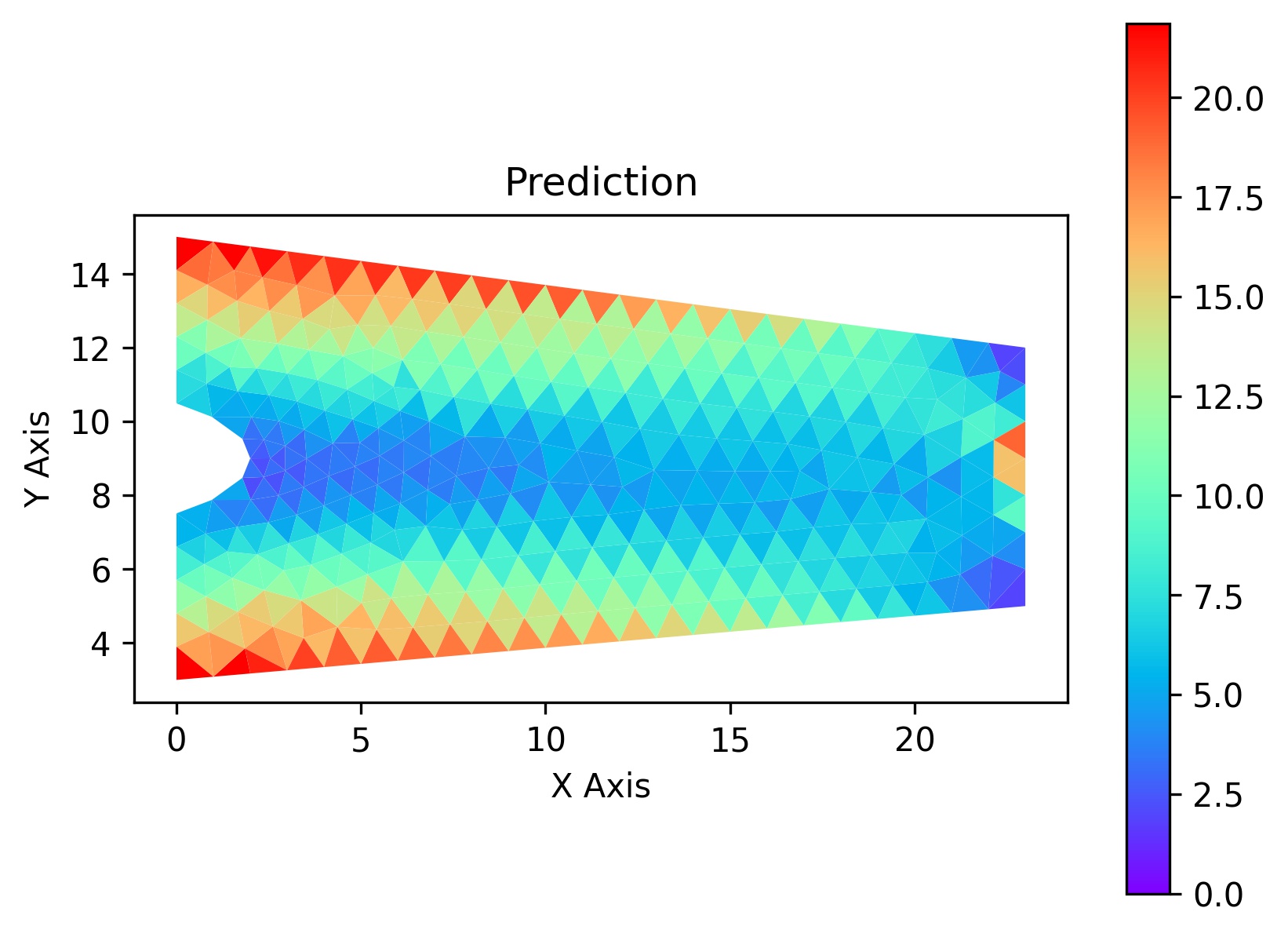}
  \caption{Shape \#6 prediction}
\end{subfigure}
\begin{subfigure}{\figWidth\textwidth}
  \centering
  \includegraphics[width=\linWidthRatio\linewidth]{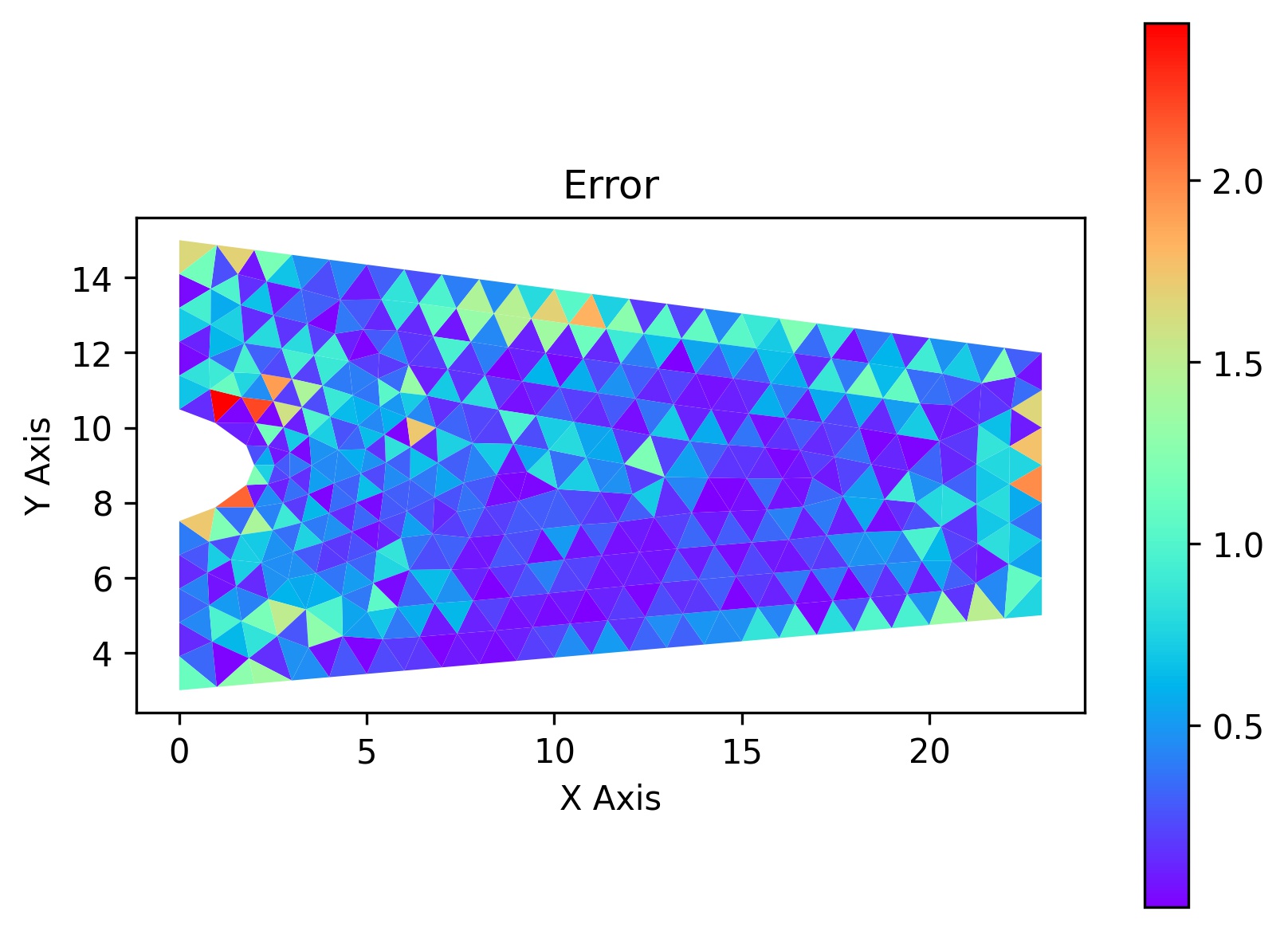}
  \caption{Shape \#6 error}
\end{subfigure}
\end{figure}

\begin{figure}[!htbp]\ContinuedFloat
\begin{subfigure}{\figWidth\textwidth}
  \centering
  \includegraphics[width=\linWidthRatio\linewidth]{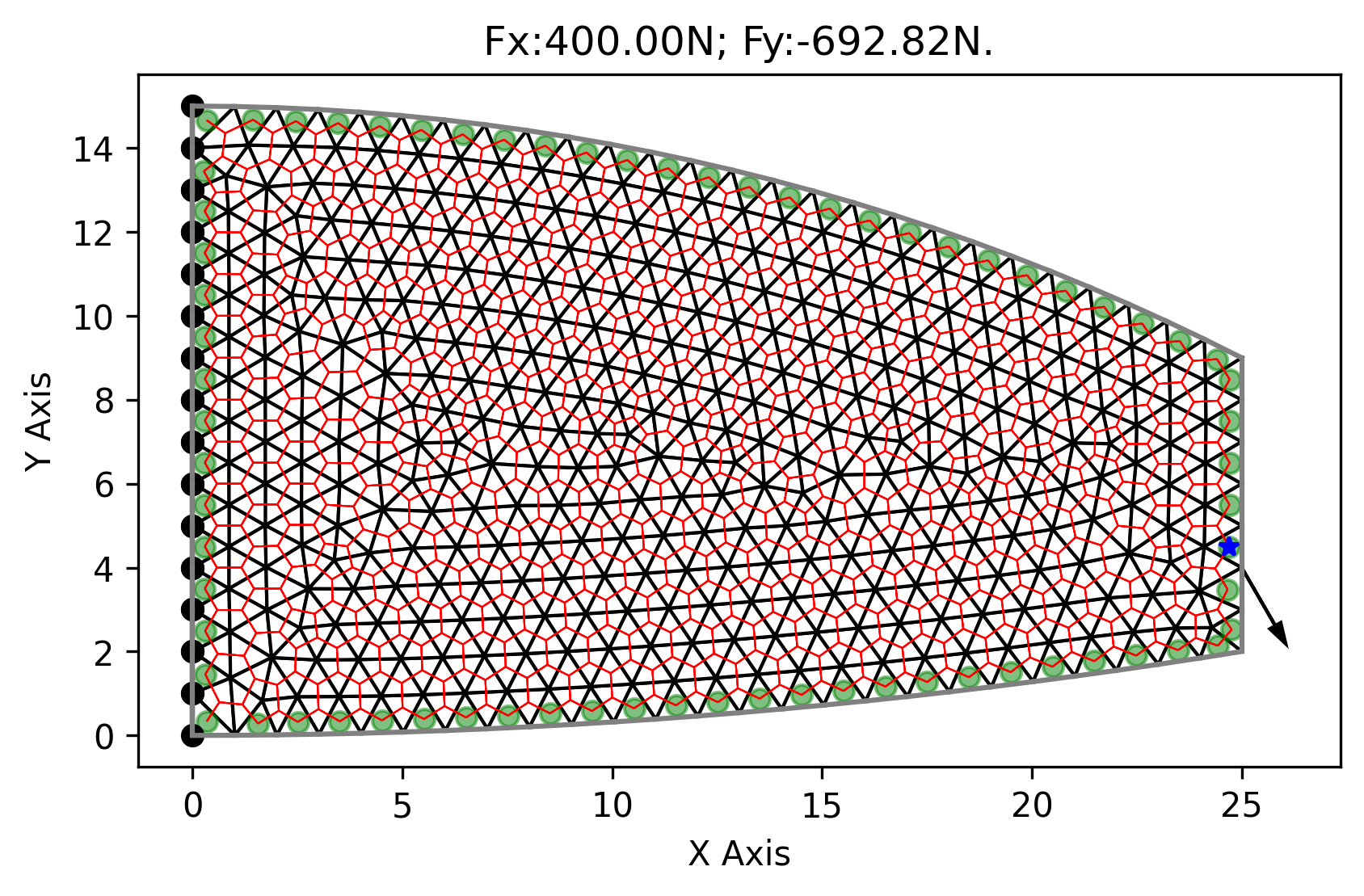}
  \caption{Shape \#7 simulation settings}
\end{subfigure}%
\begin{subfigure}{\figWidth\textwidth}
  \centering
  \includegraphics[width=\linWidthRatio\linewidth]{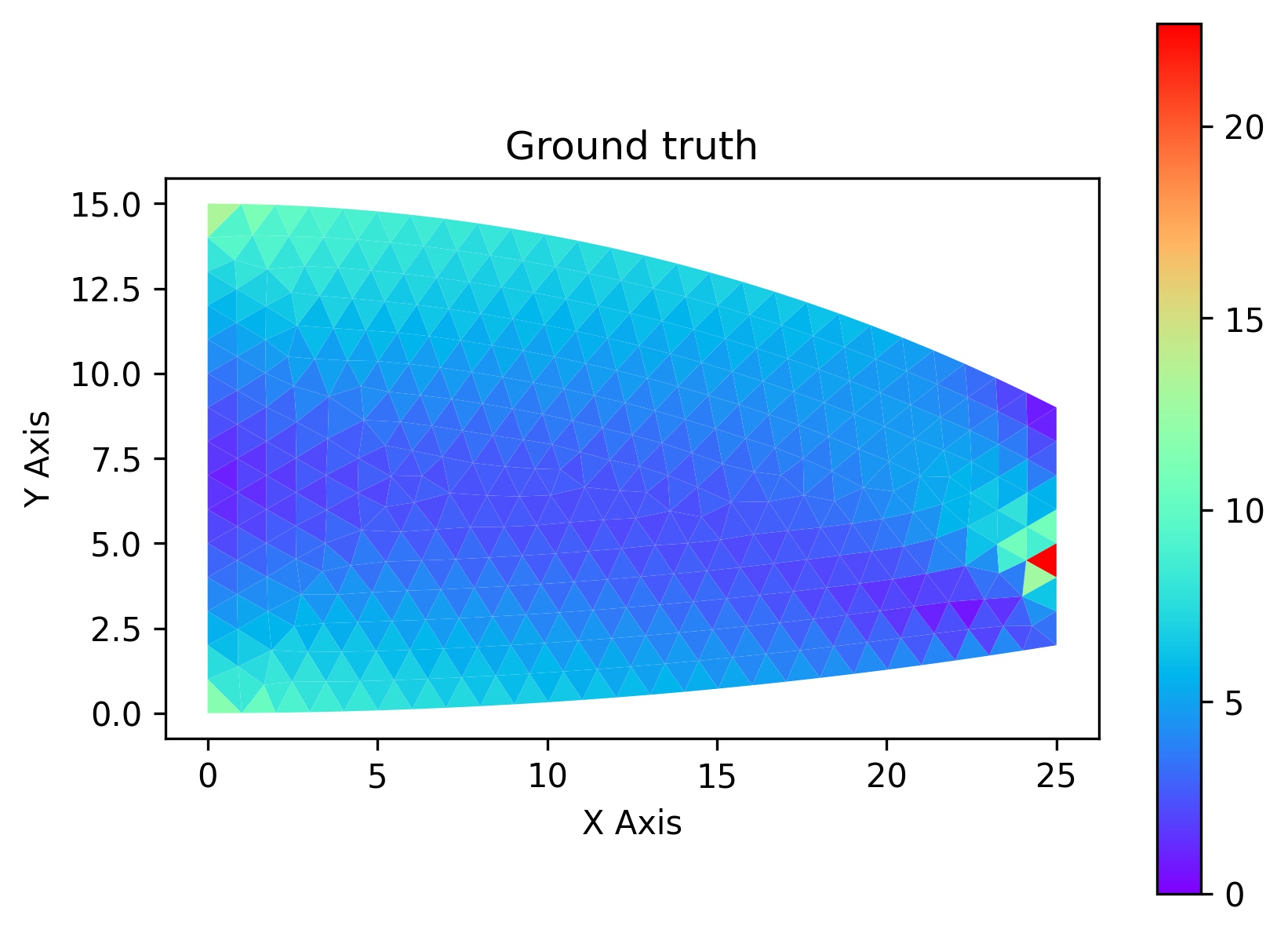}
  \caption{Shape \#7 ground truth}
\end{subfigure}
\begin{subfigure}{\figWidth\textwidth}
  \centering
  \includegraphics[width=\linWidthRatio\linewidth]{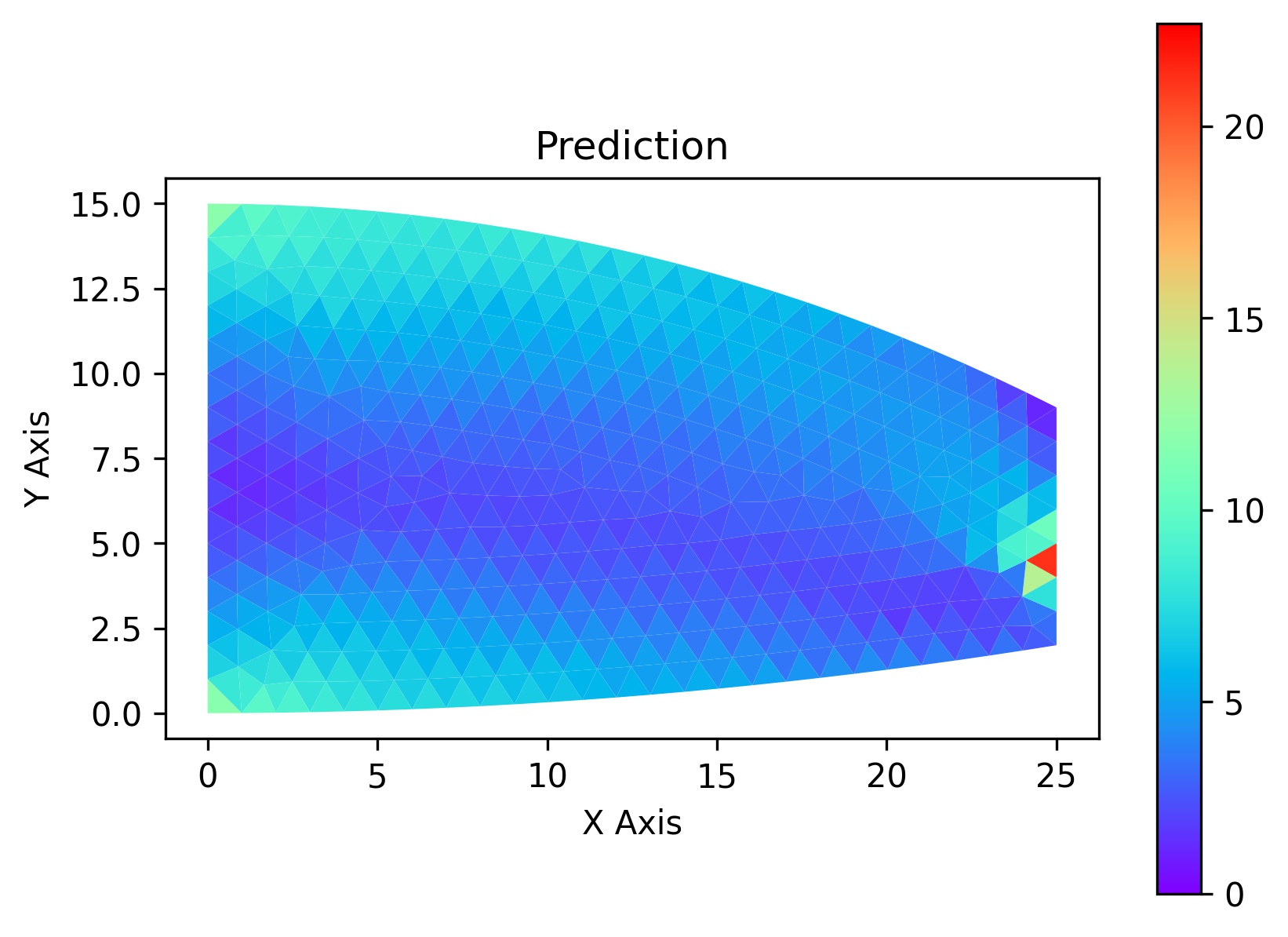}
  \caption{Shape \#7 prediction}
\end{subfigure}
\begin{subfigure}{\figWidth\textwidth}
  \centering
  \includegraphics[width=\linWidthRatio\linewidth]{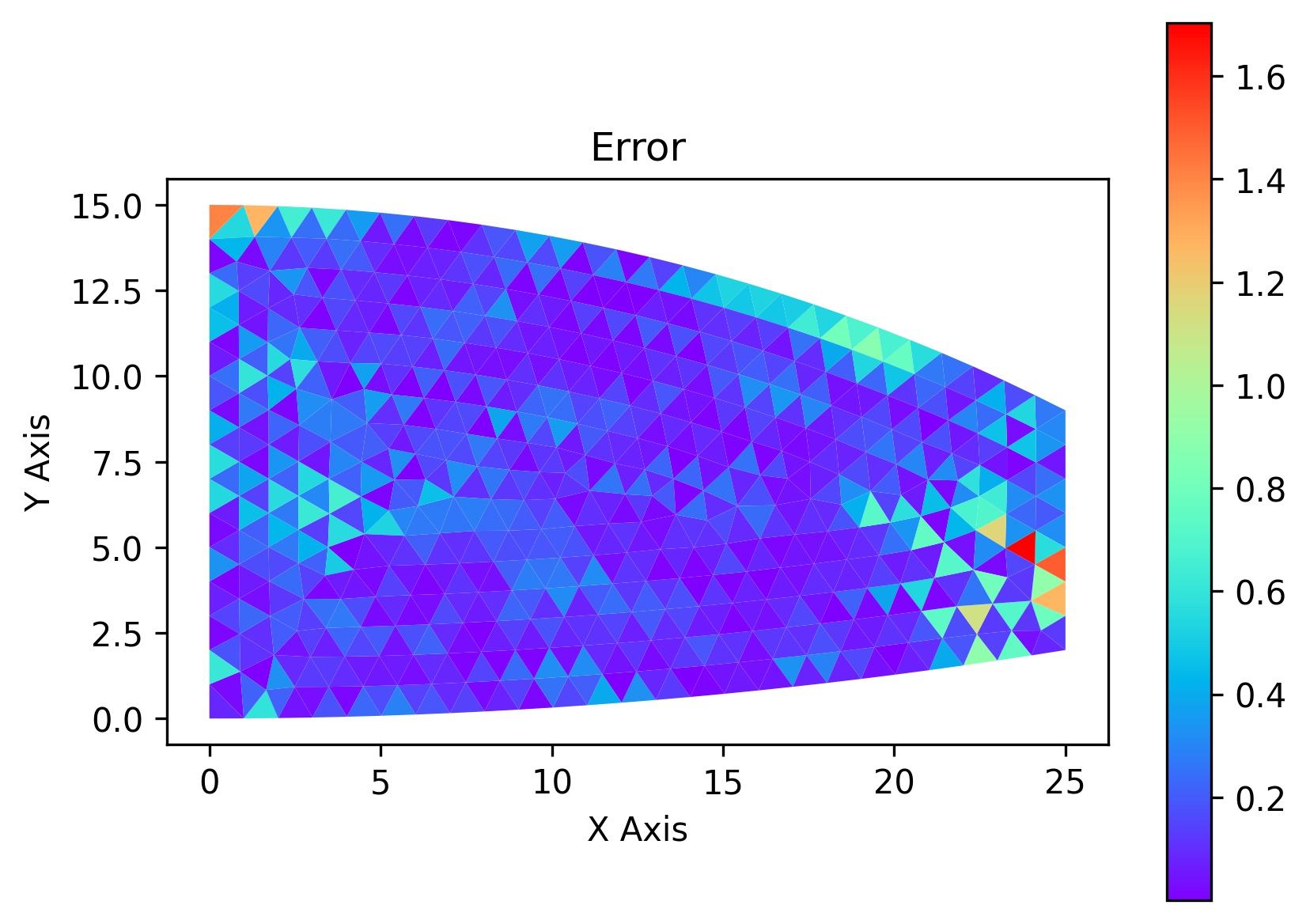}
  \caption{Shape \#7 error}
\end{subfigure}
\end{figure}

\begin{figure}[!htbp]\ContinuedFloat
\begin{subfigure}{\figWidth\textwidth}
  \centering
  \includegraphics[width=\linWidthRatio\linewidth]{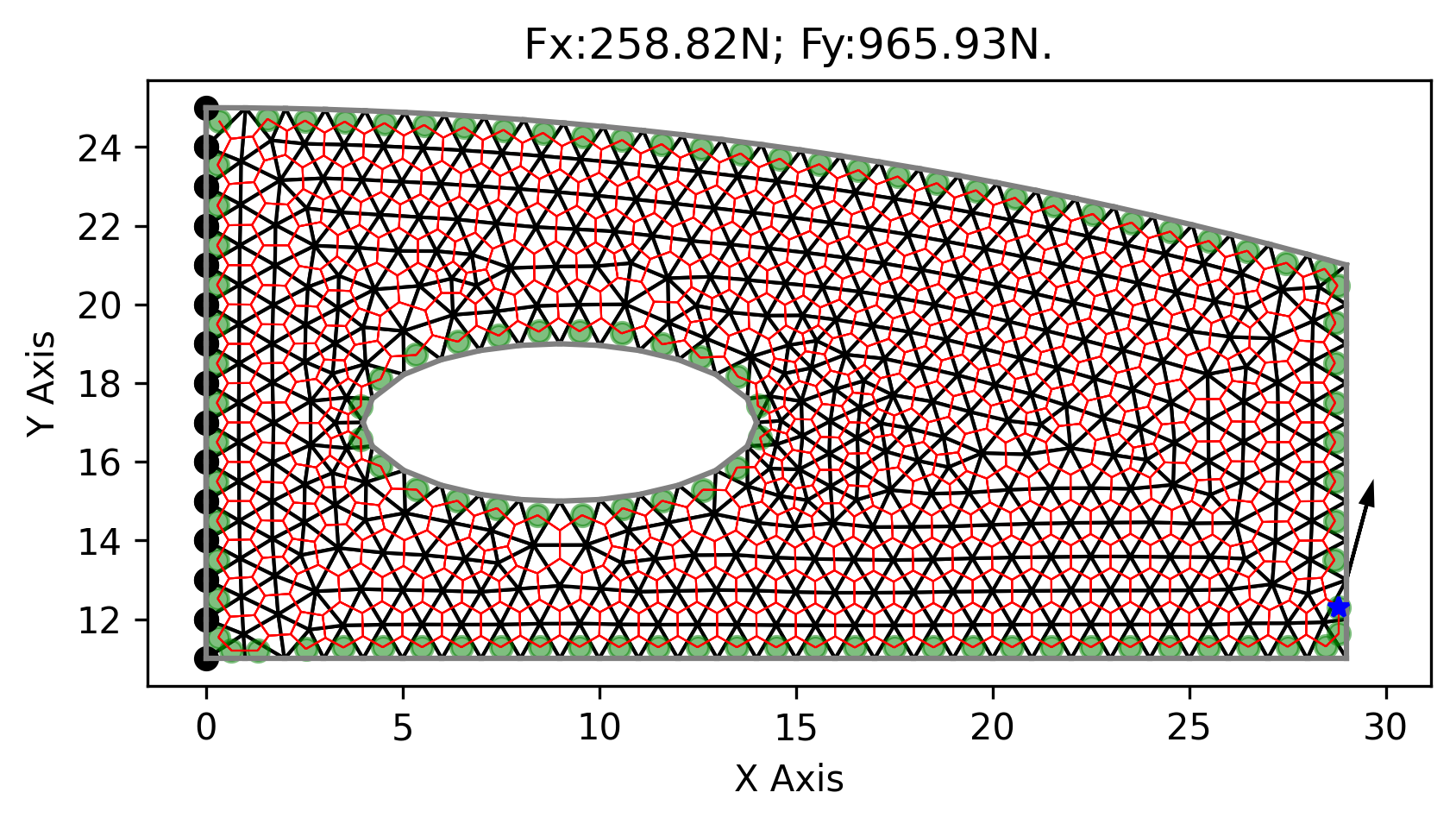}
  \caption{Shape \#8 simulation settings}
\end{subfigure}%
\begin{subfigure}{\figWidth\textwidth}
  \centering
  \includegraphics[width=\linWidthRatio\linewidth]{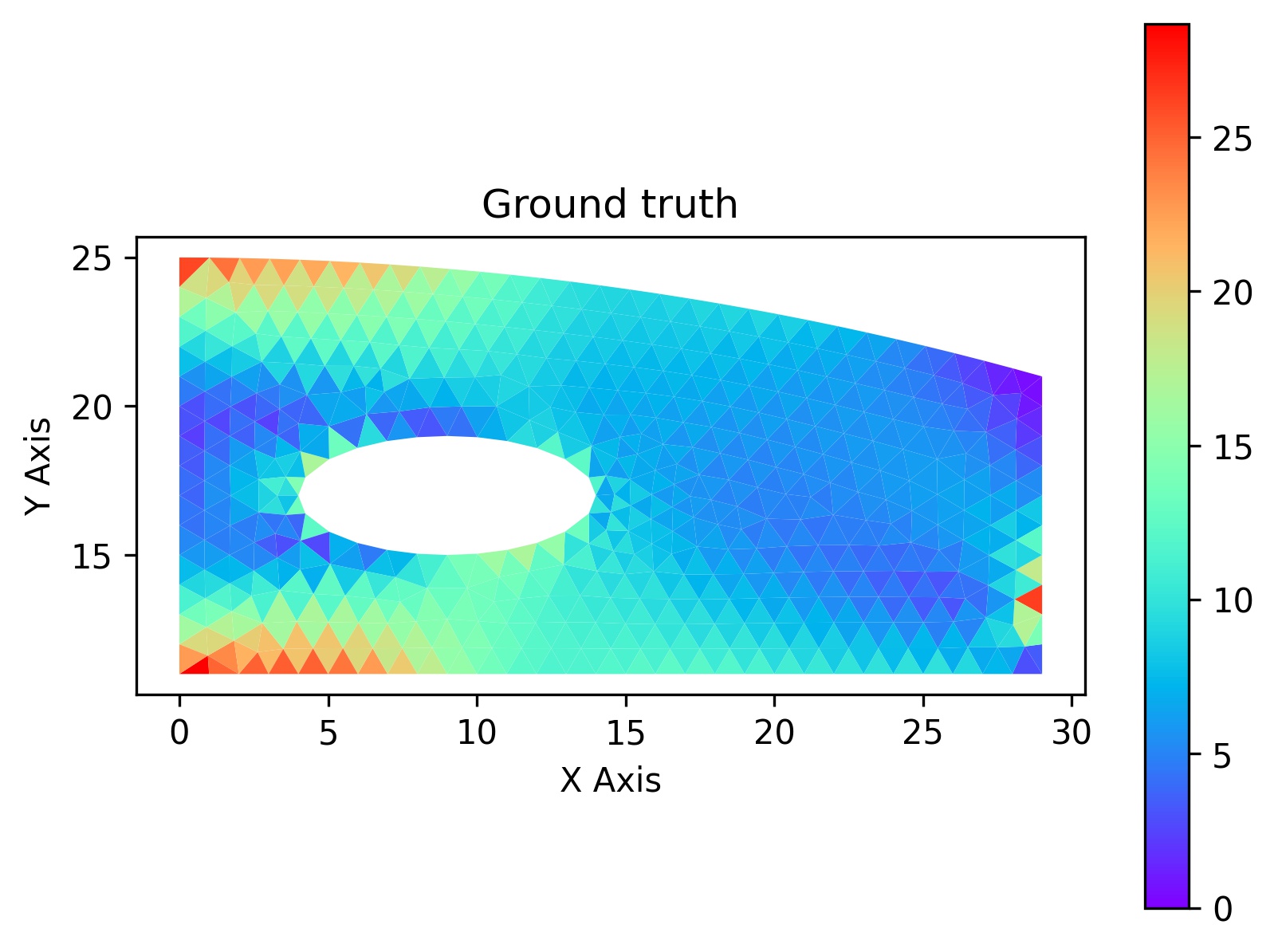}
  \caption{Shape \#8 ground truth}
\end{subfigure}
\begin{subfigure}{\figWidth\textwidth}
  \centering
  \includegraphics[width=\linWidthRatio\linewidth]{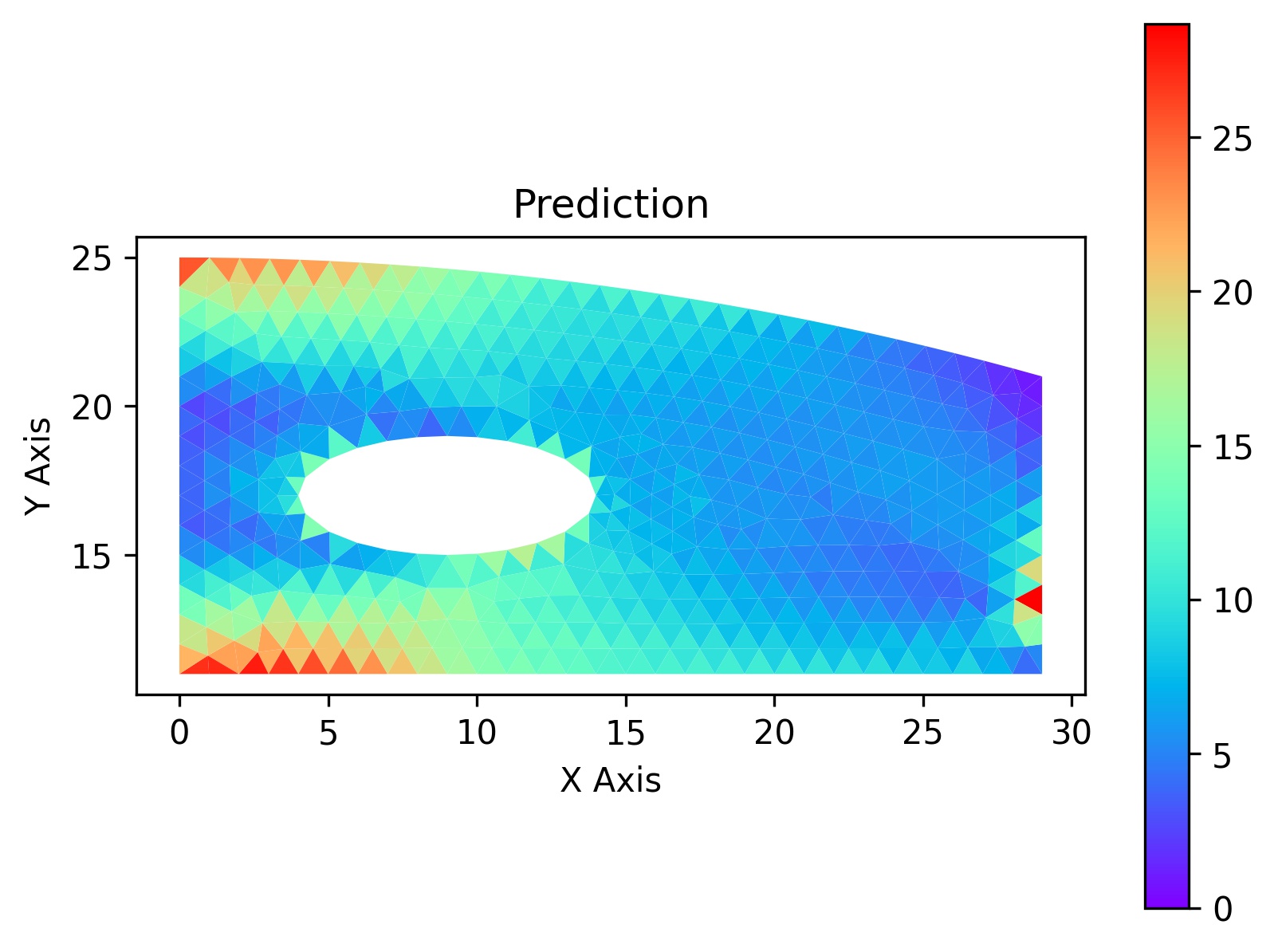}
  \caption{Shape \#8 prediction}
\end{subfigure}
\begin{subfigure}{\figWidth\textwidth}
  \centering
  \includegraphics[width=\linWidthRatio\linewidth]{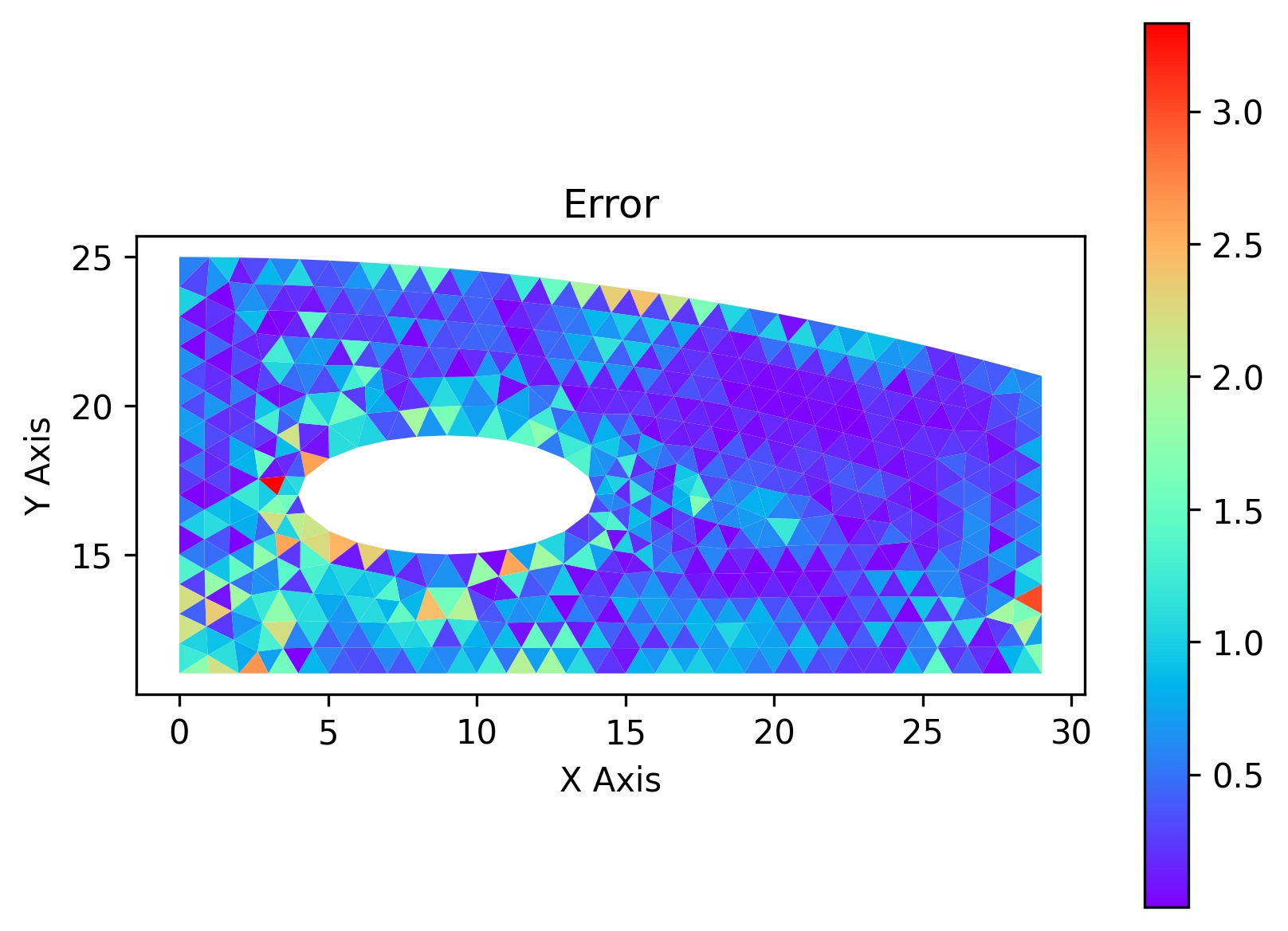}
  \caption{Shape \#8 error}
\end{subfigure}
\end{figure}

\begin{figure}[!htbp]\ContinuedFloat
\begin{subfigure}{\figWidth\textwidth}
  \centering
  \includegraphics[width=\linWidthRatio\linewidth]{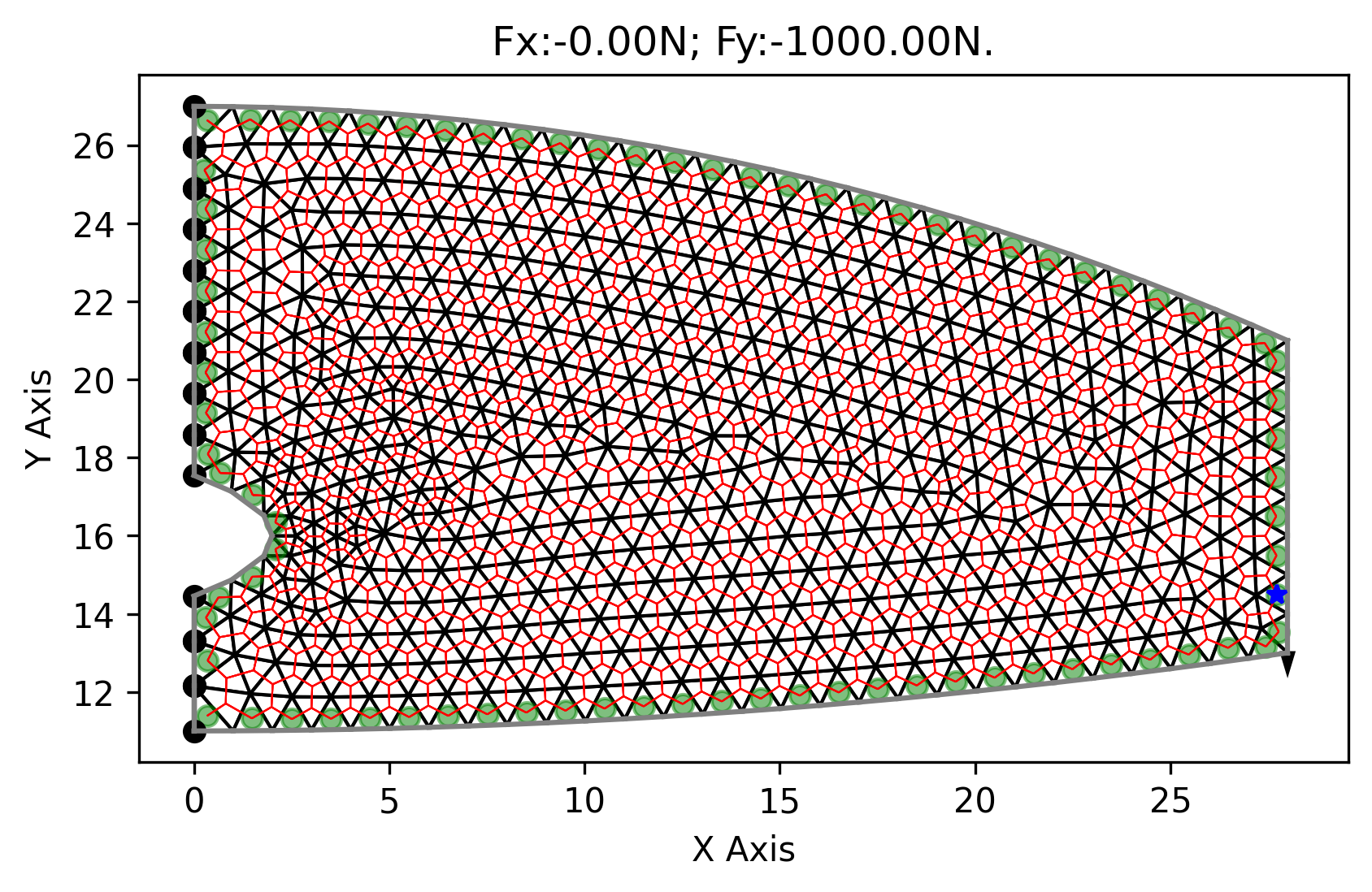}
  \caption{Shape \#9 simulation settings}
\end{subfigure}%
\begin{subfigure}{\figWidth\textwidth}
  \centering
  \includegraphics[width=\linWidthRatio\linewidth]{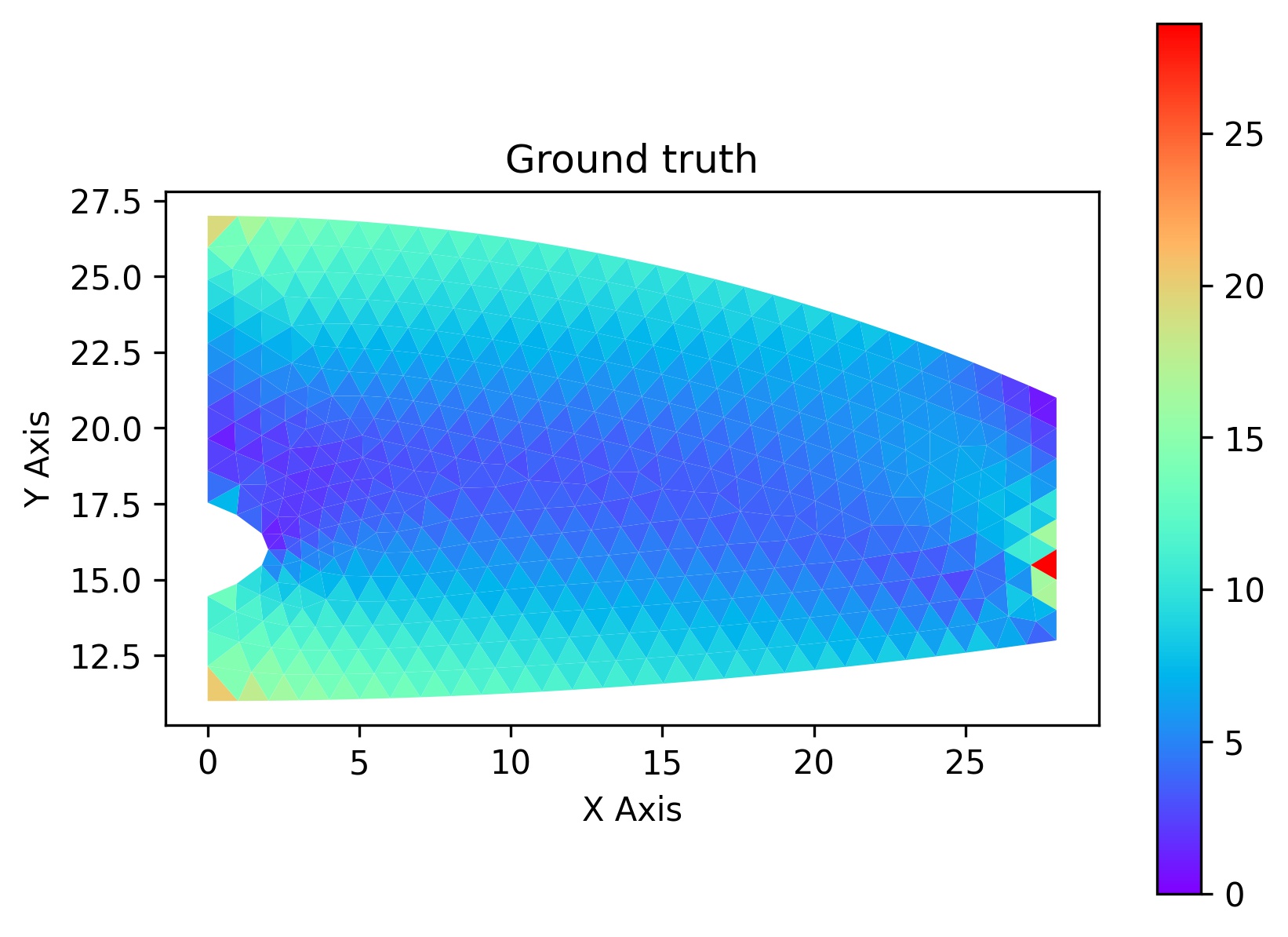}
  \caption{Shape \#9 ground truth}
\end{subfigure}
\begin{subfigure}{\figWidth\textwidth}
  \centering
  \includegraphics[width=\linWidthRatio\linewidth]{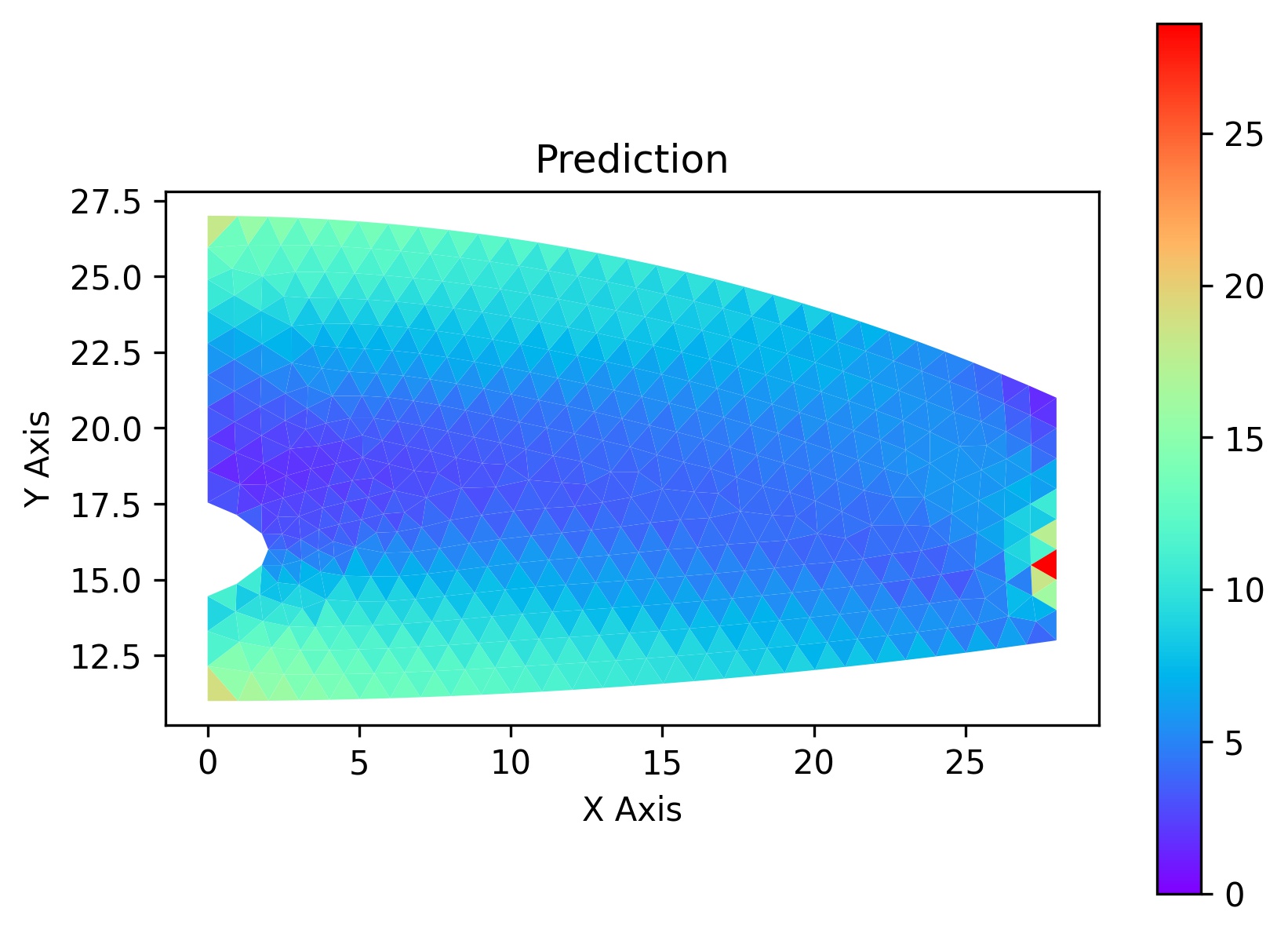}
  \caption{Shape \#9 prediction}
\end{subfigure}
\begin{subfigure}{\figWidth\textwidth}
  \centering
  \includegraphics[width=\linWidthRatio\linewidth]{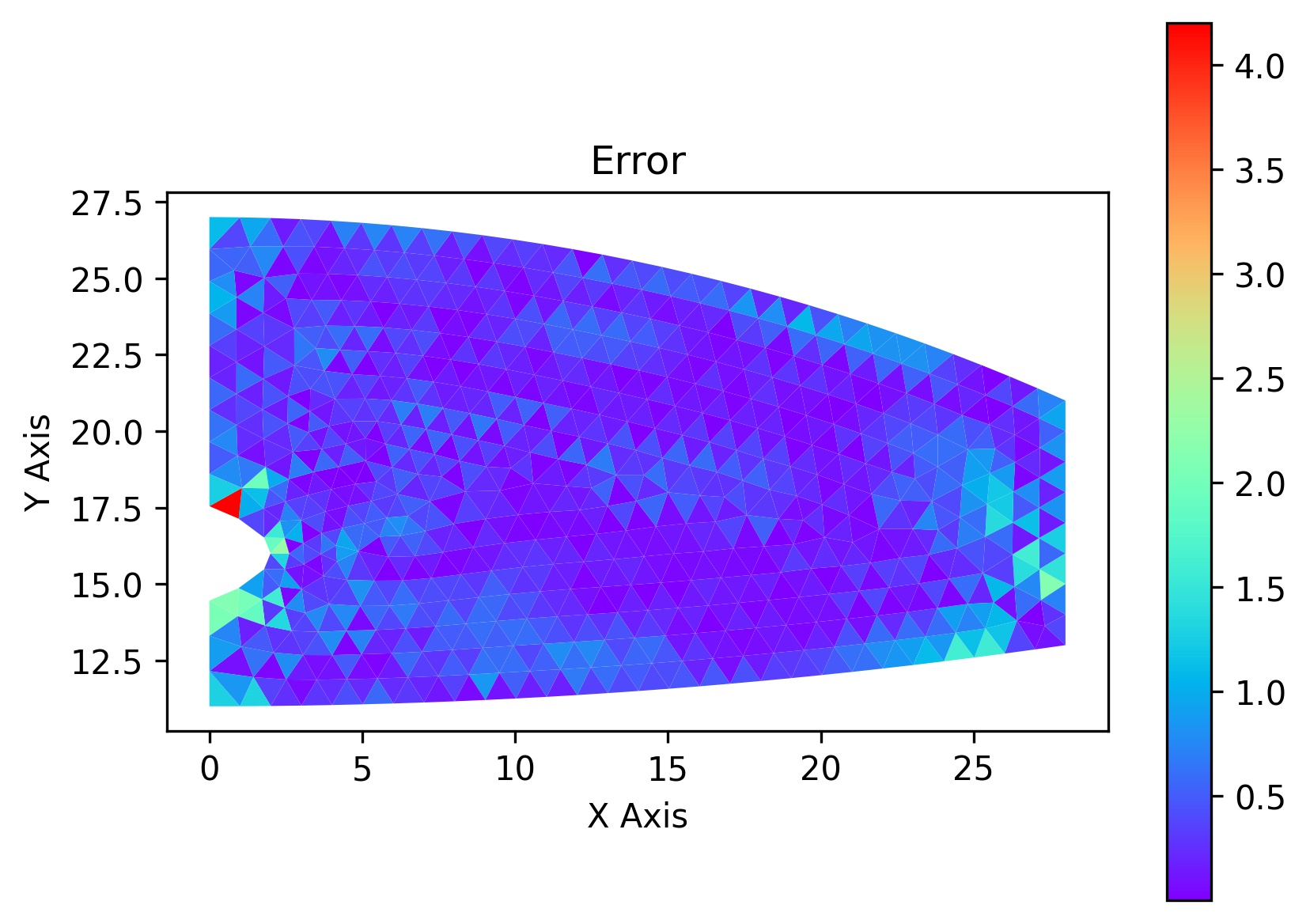}
  \caption{Shape \#9 error}
\end{subfigure}
\caption{Prediction results for von Mises stress field}
\label{fig:StressResults}
\end{figure}

\subsection{Topology optimization prediction}\label{stress}

Using the same method, we regress the topology optimization results to investigate the BOGE’s potential on solving abstract decision-making problems. Using the same DeepGCN structure, we employ the same graph embedding inputs to regress the volume density $z_e$. Since topology optimization usually happens after the stress field simulation, we also consider the condition when the stress field information is added to the features of graph vertices. In this condition, the stress field information $\mathbf{\sigma} = [\sigma_x, \sigma_y, \tau_{xy}]^T$ is attached to the graph-vertex features in BOGE approach. The training data employs BOGE with $l_e=1$ on 45k topology optimization dataset ("Topology optimization (45k)" in Table~\ref{tab:2_1_Dataset}), and its training results are shown in Table~\ref{tab:3_4_topo}. The dropout layers have been added after each layer with the dropout rate of 0.1 to provide a better prediction result.

\begin{table}[!h]
\centering
\includegraphics[keepaspectratio, width=\columnwidth]{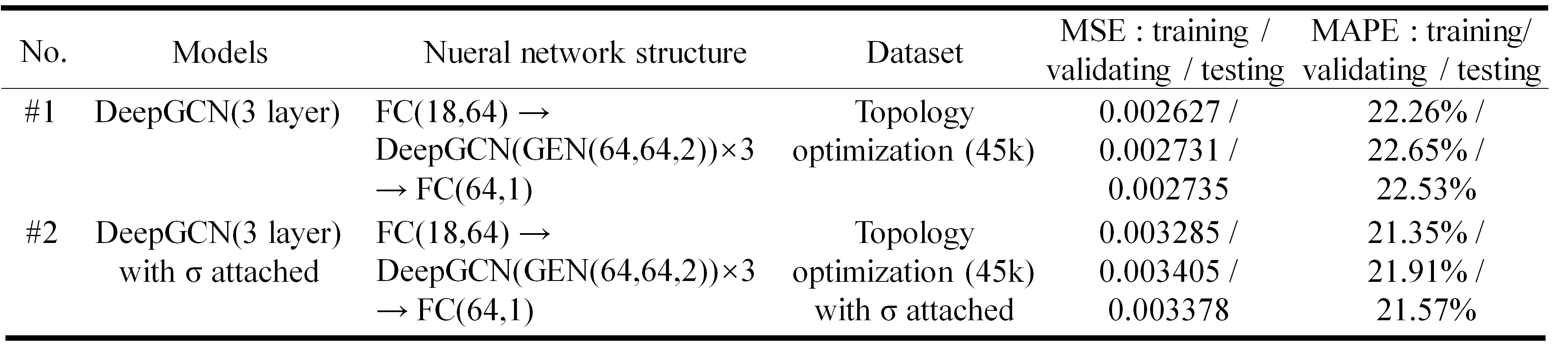}
\caption{Topology optimization prediction results (notations for neural network layers are explained in Table~\ref{tab:3_1_results}; training/validating/testing accuracy are presented without regularization)}
\label{tab:3_4_topo}
\end{table}

\begin{table}[!h]
\centering
\includegraphics[keepaspectratio, width=\columnwidth]{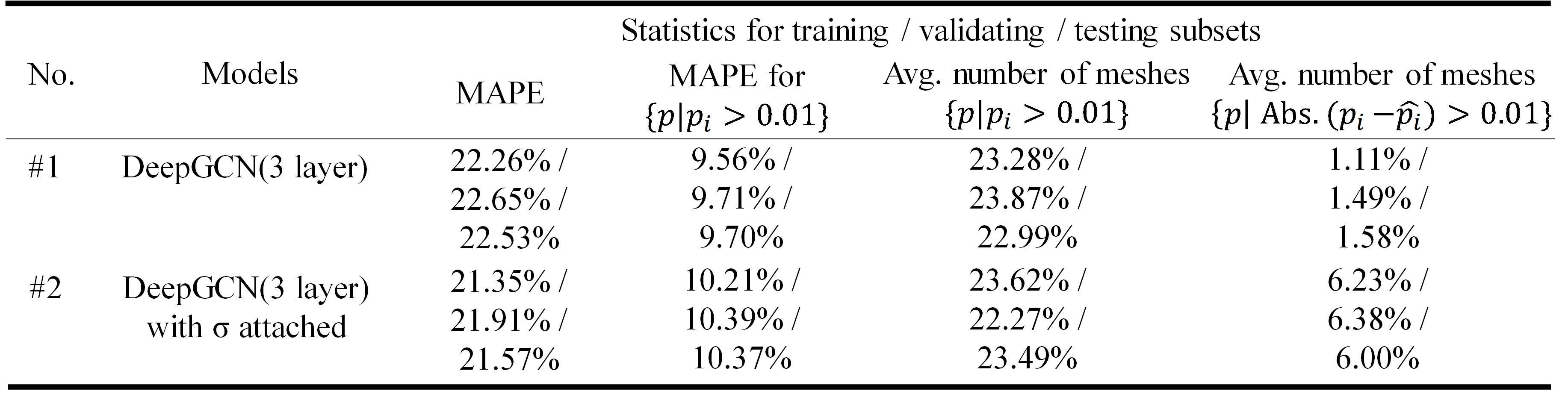}
\caption{Analyses for topology optimization prediction results (training/validating/testing accuracy are presented without regularization)}
\label{tab:3_5_acc}
\end{table}

From the training results, it can be seen that BOGE without encoding the stress information (model \#1 in Table~\ref{tab:3_4_topo}) reaches 0.002735 testing accuracy, which is better than the best testing MSE - 0.059943 in \cite{nie2021topologygan} with similar cantilever beam settings and the comparable number of meshes. However, the MAPE is relatively large and reaches 22.53\% in testing accuracy. It is due to two effects. First, the range or the average value for the volume density in the topological optimization is relatively small. For stress prediction shown in Table~\ref{tab:2_1_Dataset}, the average von Mises stress is around 2MPa while the maximum stress magnitude can reach 74.09MPa, which means the error divided by its ground truth value can be relatively small for the MAPE. Assume one element with the stress of 2MPa (which is the average stress) has the prediction error of 0.1MPa, the MAPE for this element is only 5\%. However, from Table~\ref{tab:2_1_Dataset}, the average volume density is around 0.29 whereas the maximum volume density is 1.0, which means that the same error can contribute largely to the final MAPE. For example, if the element with 0.29 volume density has the prediction error of 0.1, the MAPE for this element can be 34\%. This explains one of the reasons that some researchers haven't employed the MAPE for evaluating topology optimization results \cite{nie2021topologygan, zhang2019deep, banga20183d}. Further, the large MAPE is generated because some outliers around the edge of the optimized beam contribute to most of the errors. Assume that the resolution of the volume density is 0.01, which is also the $\varepsilon$ in Eq.~(\ref{mape}). From model \#1 in Table~\ref{tab:3_5_acc}, when only meshes with its ground-truth volume density larger than $\varepsilon$ are considered (${p|p_i>0.01}$), the MAPE for those meshes is only 9.70\%. Also, for each testing graph (or simulation model), around 22\% - 23\% of its meshes have the volume density larger than 0.01 for its ground truth value, but only 1.58\% of the total meshes contribute to the errors larger than $\varepsilon$. Some prediction results from model \#1 in Table~\ref{tab:3_4_topo} are shown in Fig.~\ref{fig:TopoResults}. It can be seen that those 1.58\% of meshes mostly are located primarily on the edges of the optimized beams. Since the boundary of the optimized beams usually has a volume density around $\varepsilon=0.01$ as its ground truth value, any small prediction error can be largely amplified by the MAPE. However, regardless of the MAPE, the predicted volume density shown in Fig.~\ref{fig:TopoResults} and its testing accuracy are sufficient to validate the BOGE's advantage in regressing complex decision-making problems. 

With the stress field information, BOGE shows a worse prediction result, which reaches the MSE of 0.003378 for testing accuracy (model \#2 in Table~\ref{tab:3_4_topo}). Also, more outliers (6.00\% from model \#2 in Table~\ref{tab:3_5_acc}) are found in this condition. Adding additional features to the graph vertex can lead to message congestion\cite{loukas2019graph} and oversquashing\cite{alon2020bottleneck}. Therefore, unlike the CNN-based surrogate model\cite{nie2021topologygan}, GNN models are sensitive to the size of the input tensor and cannot be simply enhanced by adding extra input features. Balancing the input feature, hidden layer features, and the GNN layers can provide better results according to the GNN's property.

The topology optimization prediction shows BOGE’s capability of making abstract decisions which benefits the smart design technology. The average prediction time for the topology-optimization GNN model (model \#1 in Table~\ref{tab:3_4_topo}) on GPU is around 0.011ms, similar to that of GNN stress field prediction, but far less than the ABAQUS computation time (around 489s in Table~\ref{tab:2_1_Dataset}). This is because the SIMP method requires several iterations of computations that take much longer than the stress prediction. However, the GNN surrogate model can directly output the optimized result which largely saves the time cost. The overall performance of BOGE demonstrates its irreplaceable advantage in regressing the physical field on structured elements.The training loss and the validation/testing accuracy for the stress predictions is shown in Fig.~\ref{fig:loss}(\subref{fig:TopoLoss})

\def \figWidth{0.24} 
\def \linWidthRatio{0.7} 

\def \idx{1~} 
\begin{figure}[!htbp]
\centering
\includegraphics[width=\linWidthRatio\linewidth]{fig/3_0_Legend.png}
\newline
\begin{subfigure}{\figWidth\textwidth}
  \centering
  \includegraphics[width=\linWidthRatio\linewidth]{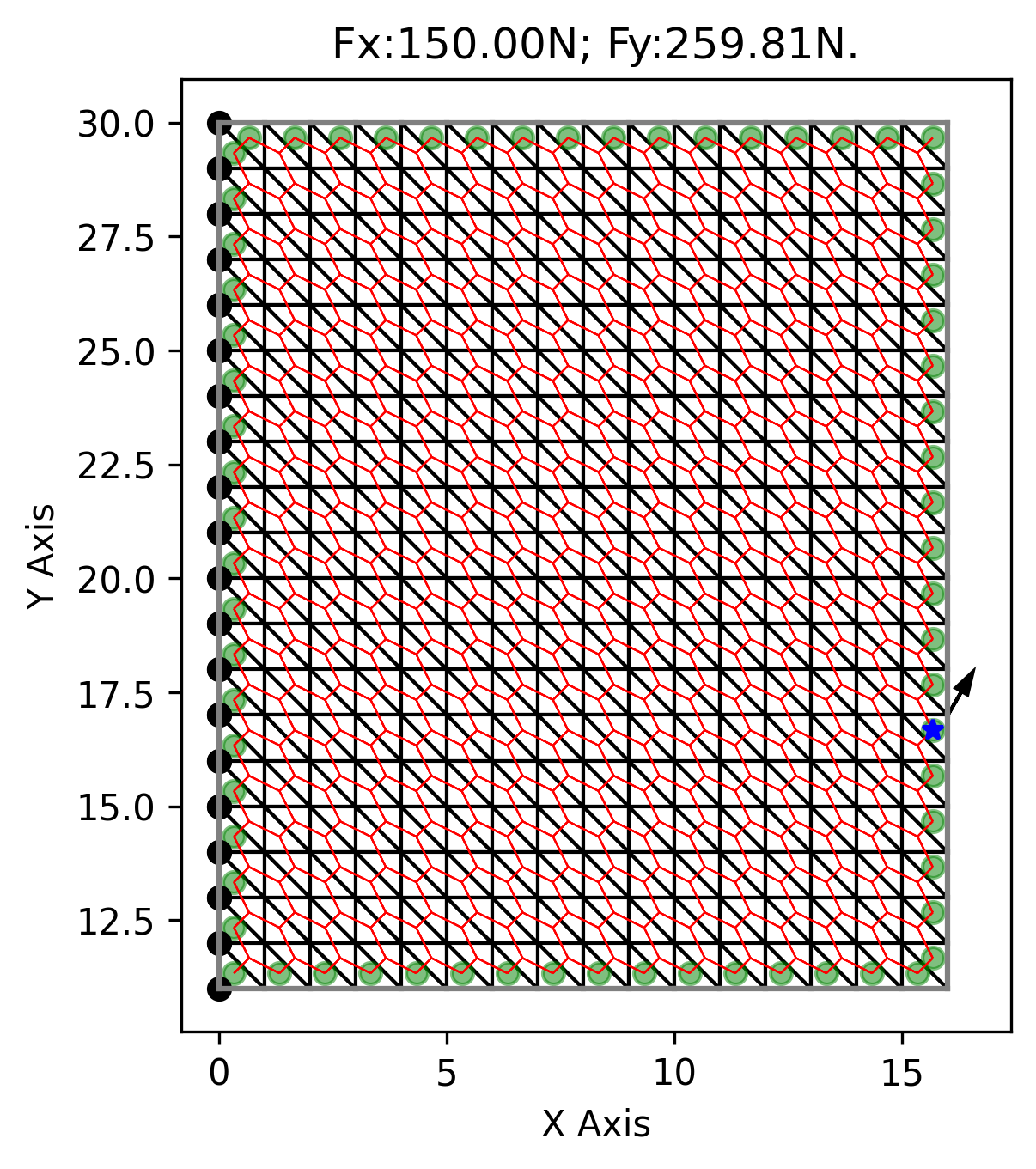}
  \caption{Shape \#\idx simulation settings}
\end{subfigure}%
\begin{subfigure}{\figWidth\textwidth}
  \centering
  \includegraphics[width=\linWidthRatio\linewidth]{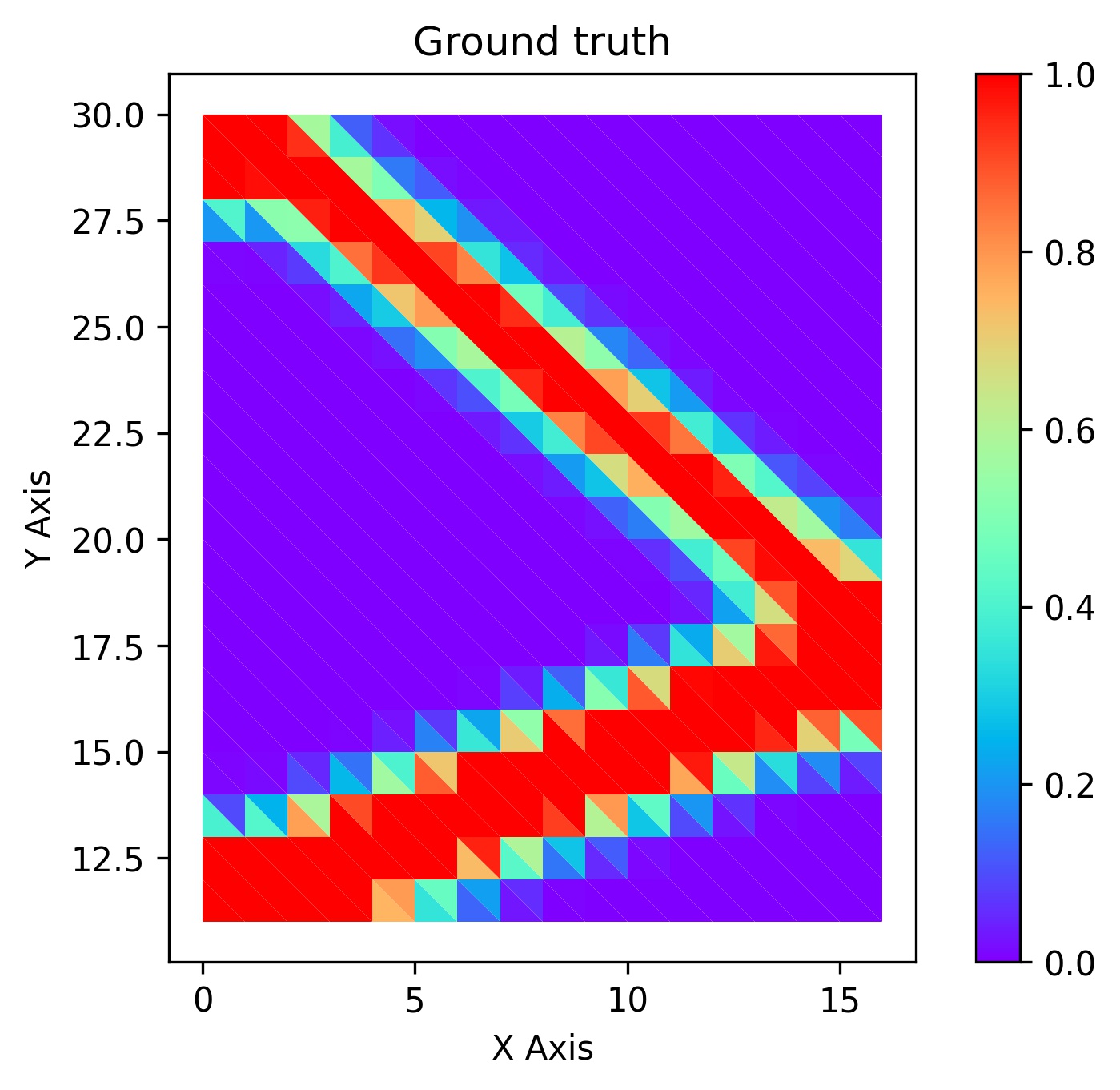}
  \caption{Shape \#\idx ground truth}
\end{subfigure}
\begin{subfigure}{\figWidth\textwidth}
  \centering
  \includegraphics[width=\linWidthRatio\linewidth]{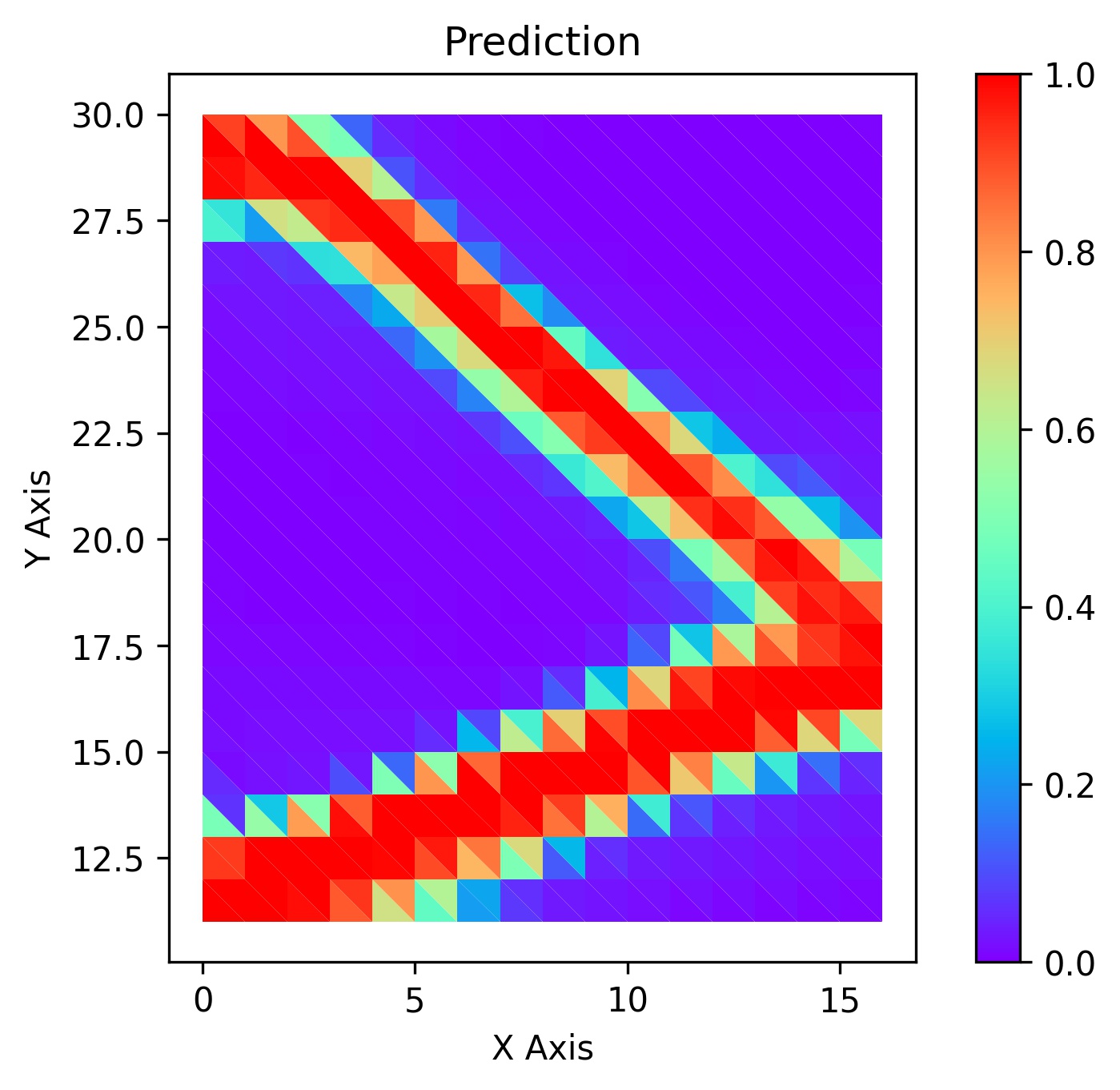}
  \caption{Shape \#\idx prediction}
\end{subfigure}
\begin{subfigure}{\figWidth\textwidth}
  \centering
  \includegraphics[width=\linWidthRatio\linewidth]{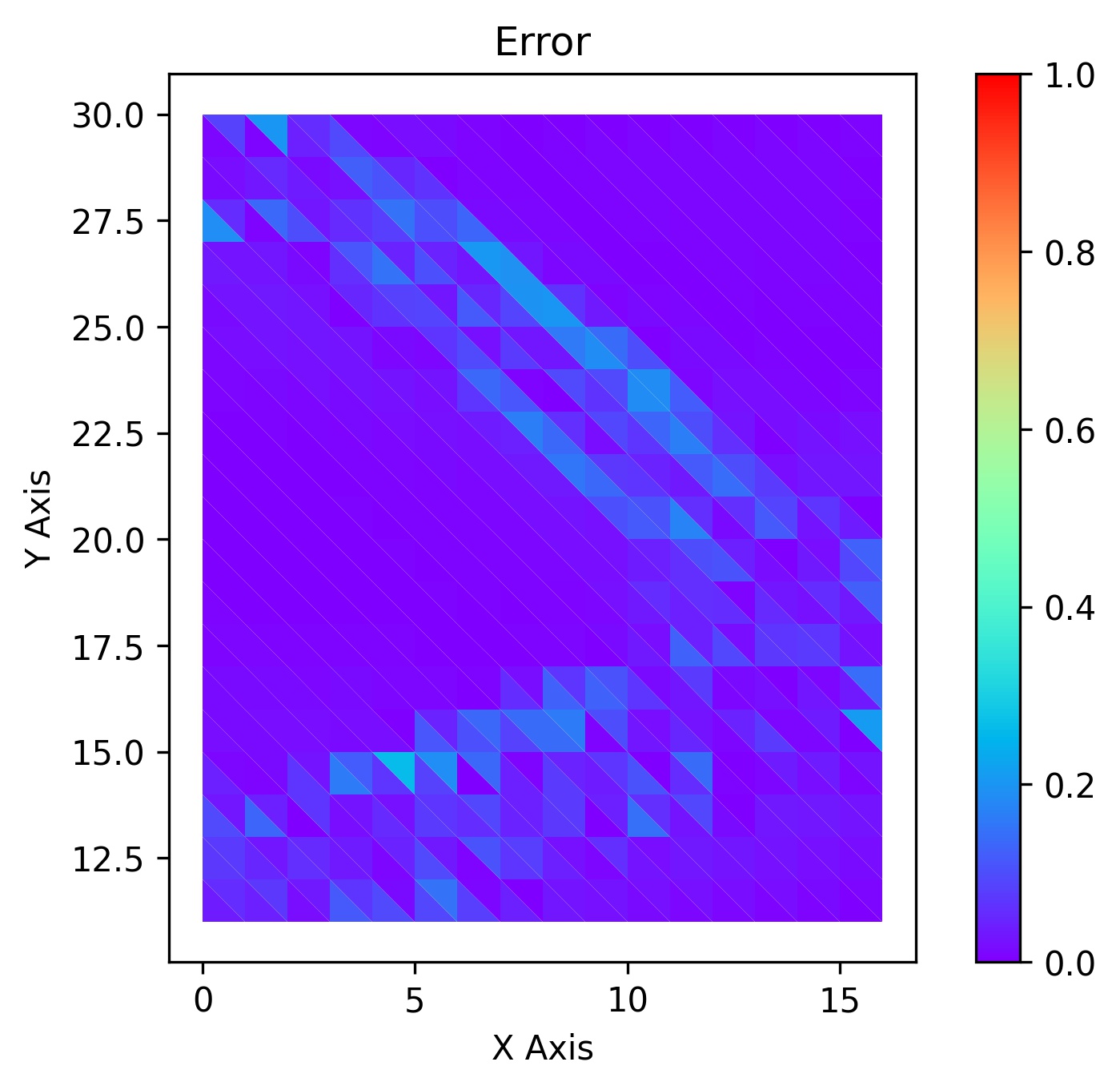}
  \caption{Shape \#\idx error}
\end{subfigure}
\end{figure}

\def \idx{2~} 
\begin{figure}[!htbp]\ContinuedFloat
\begin{subfigure}{\figWidth\textwidth}
  \centering
  \includegraphics[width=\linWidthRatio\linewidth]{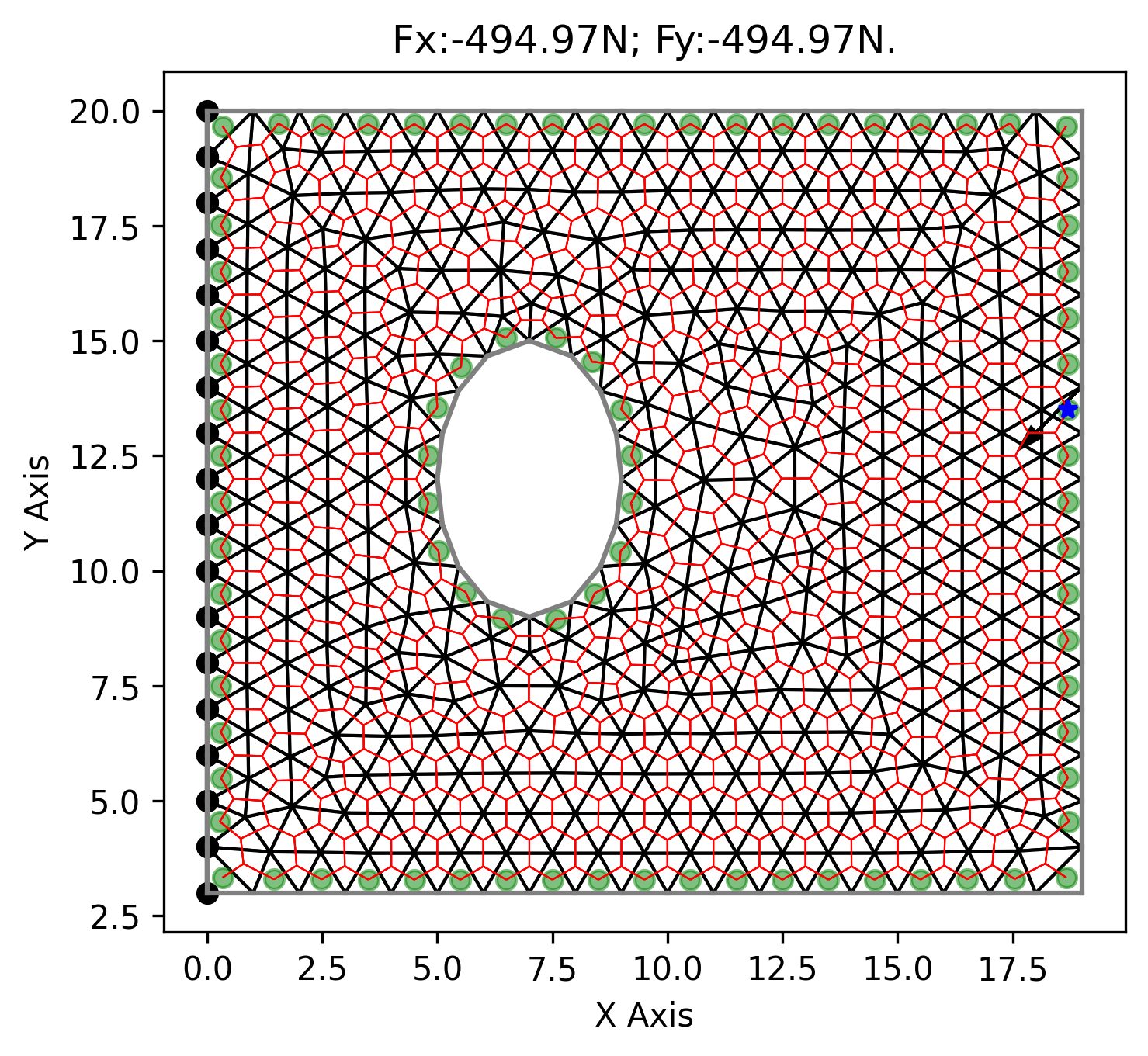}
  \caption{Shape \#\idx simulation settings}
\end{subfigure}%
\begin{subfigure}{\figWidth\textwidth}
  \centering
  \includegraphics[width=\linWidthRatio\linewidth]{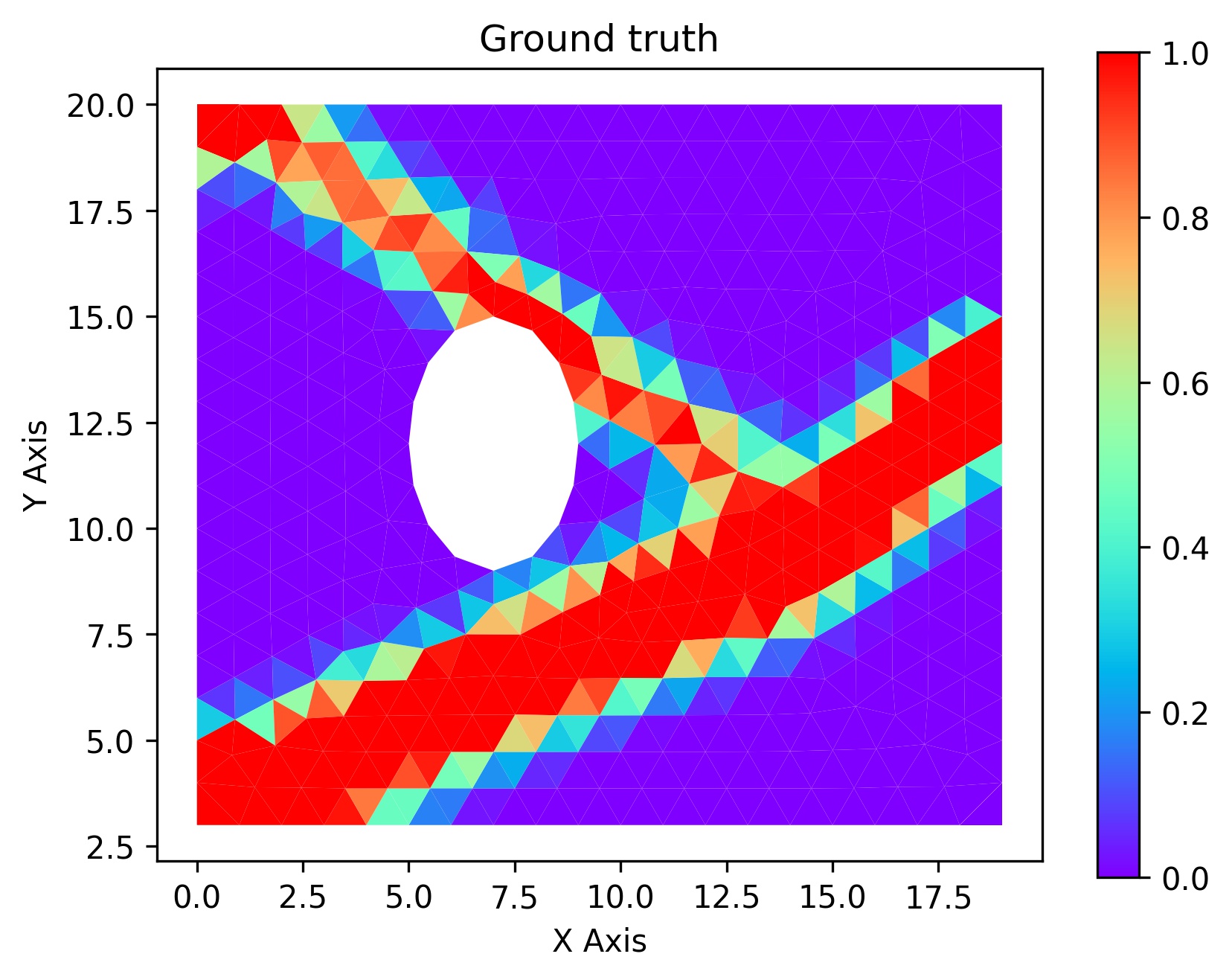}
  \caption{Shape \#\idx ground truth}
\end{subfigure}
\begin{subfigure}{\figWidth\textwidth}
  \centering
  \includegraphics[width=\linWidthRatio\linewidth]{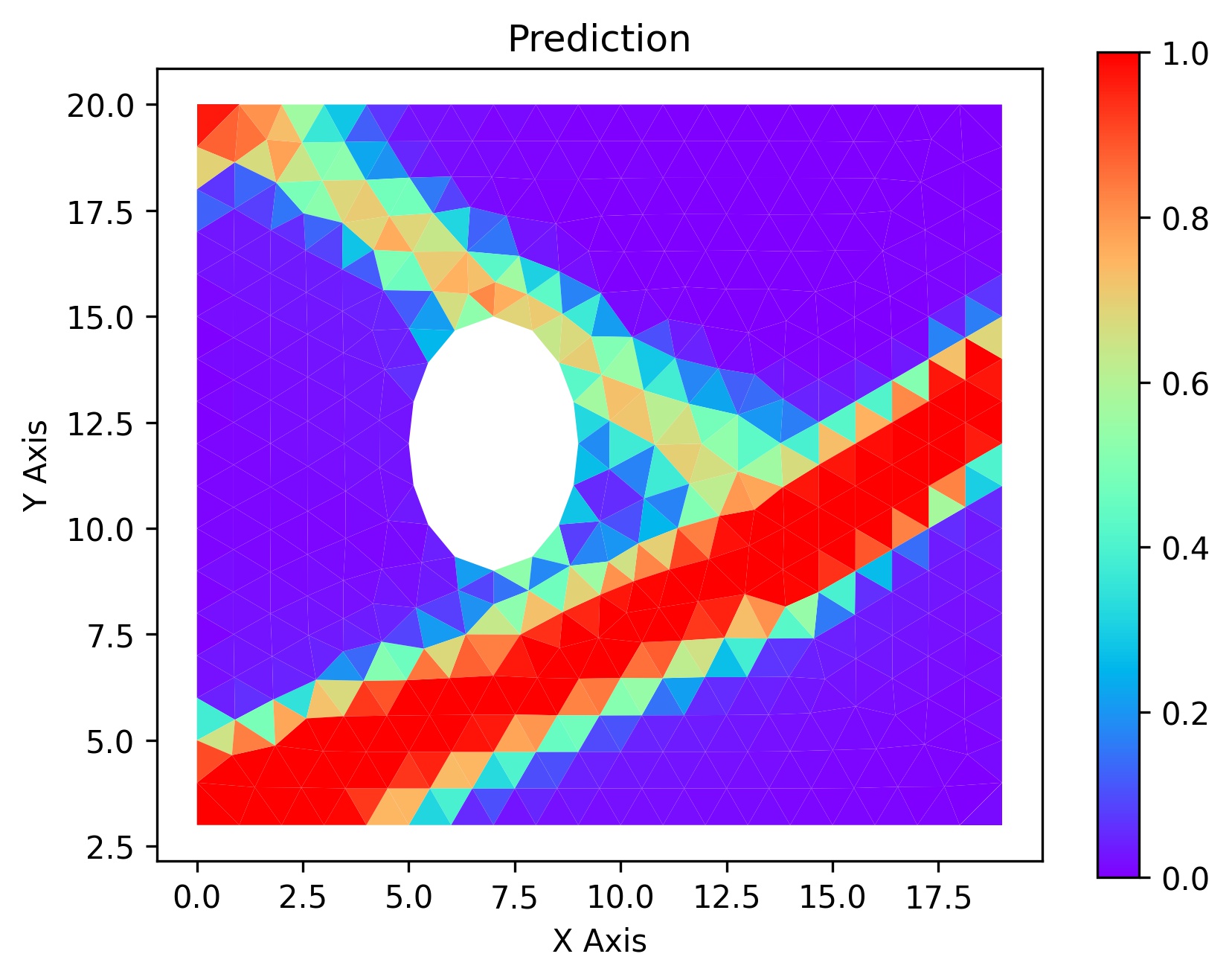}
  \caption{Shape \#\idx prediction}
\end{subfigure}
\begin{subfigure}{\figWidth\textwidth}
  \centering
  \includegraphics[width=\linWidthRatio\linewidth]{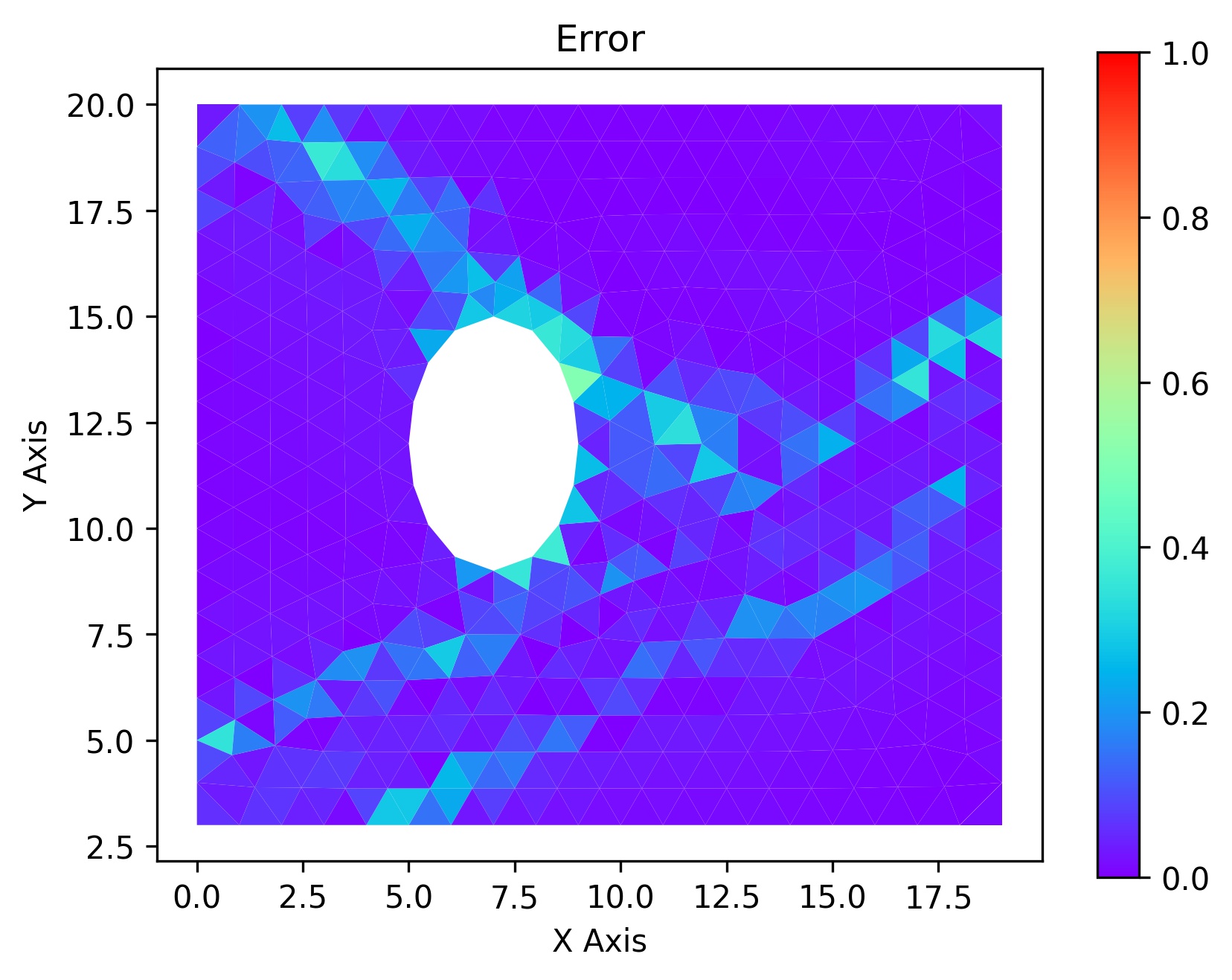}
  \caption{Shape \#\idx error}
\end{subfigure}
\end{figure}

\def \idx{3~} 
\begin{figure}[!htbp]\ContinuedFloat
\begin{subfigure}{\figWidth\textwidth}
  \centering
  \includegraphics[width=\linWidthRatio\linewidth]{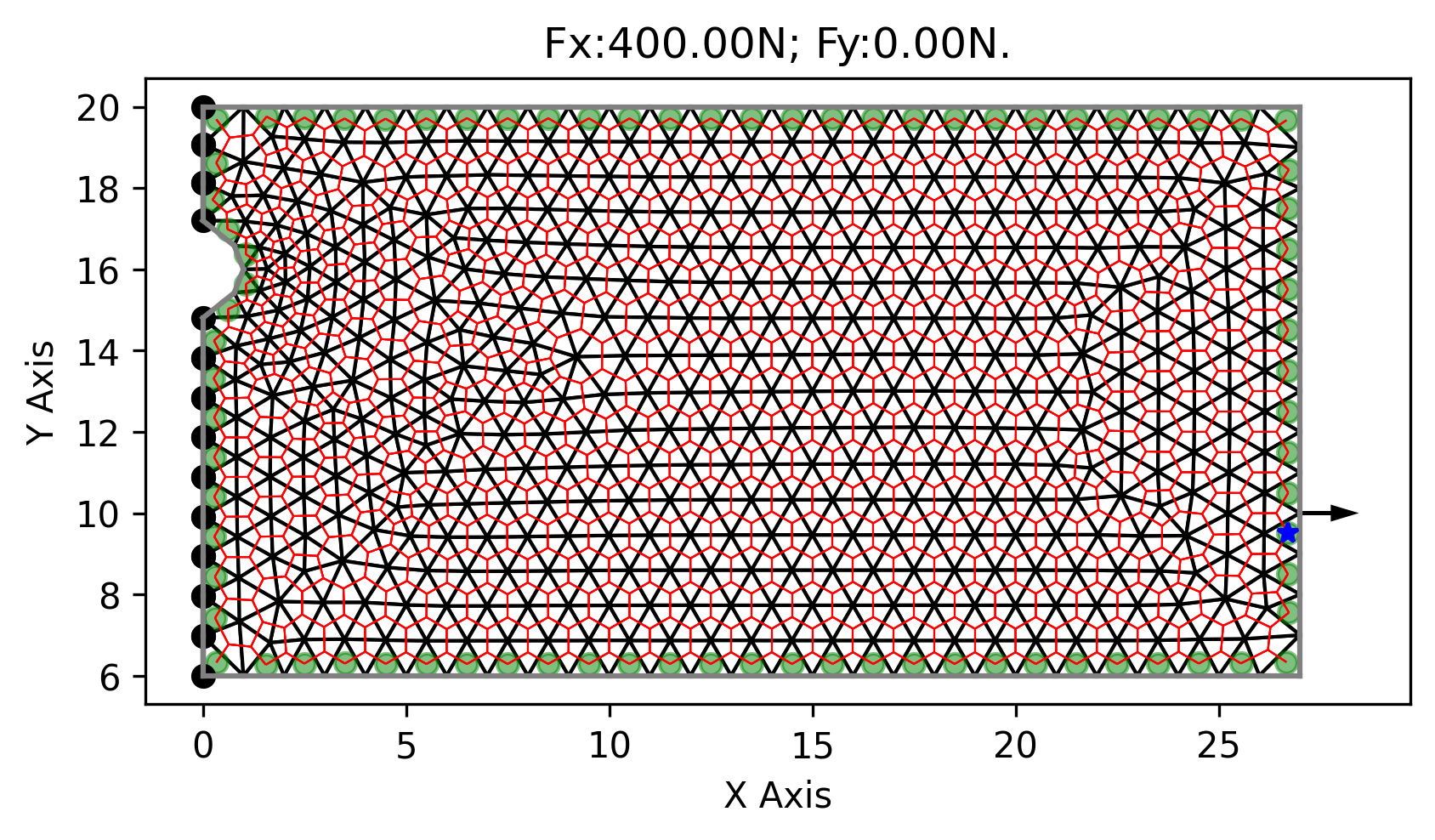}
  \caption{Shape \#\idx simulation settings}
\end{subfigure}%
\begin{subfigure}{\figWidth\textwidth}
  \centering
  \includegraphics[width=\linWidthRatio\linewidth]{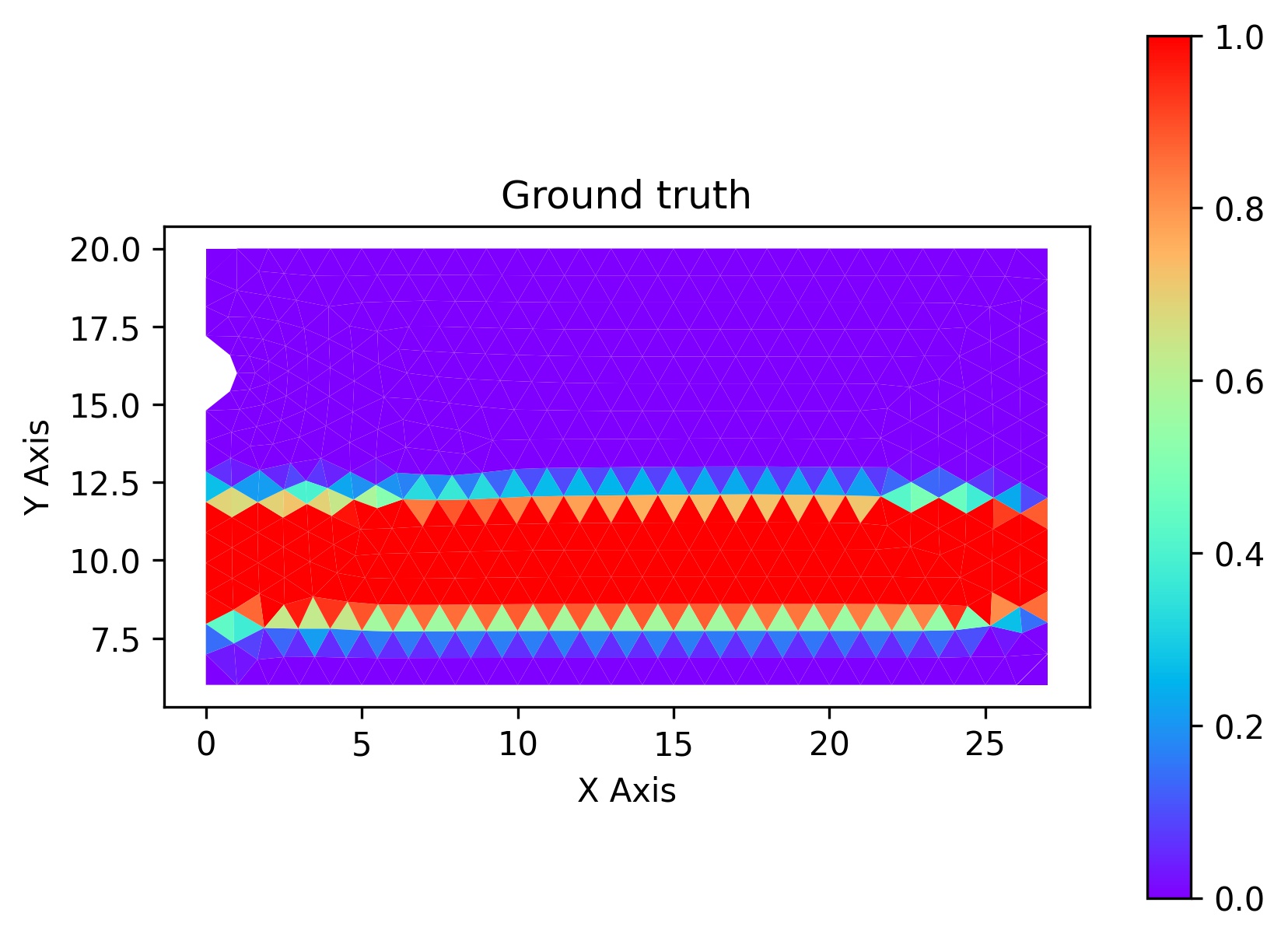}
  \caption{Shape \#\idx ground truth}
\end{subfigure}
\begin{subfigure}{\figWidth\textwidth}
  \centering
  \includegraphics[width=\linWidthRatio\linewidth]{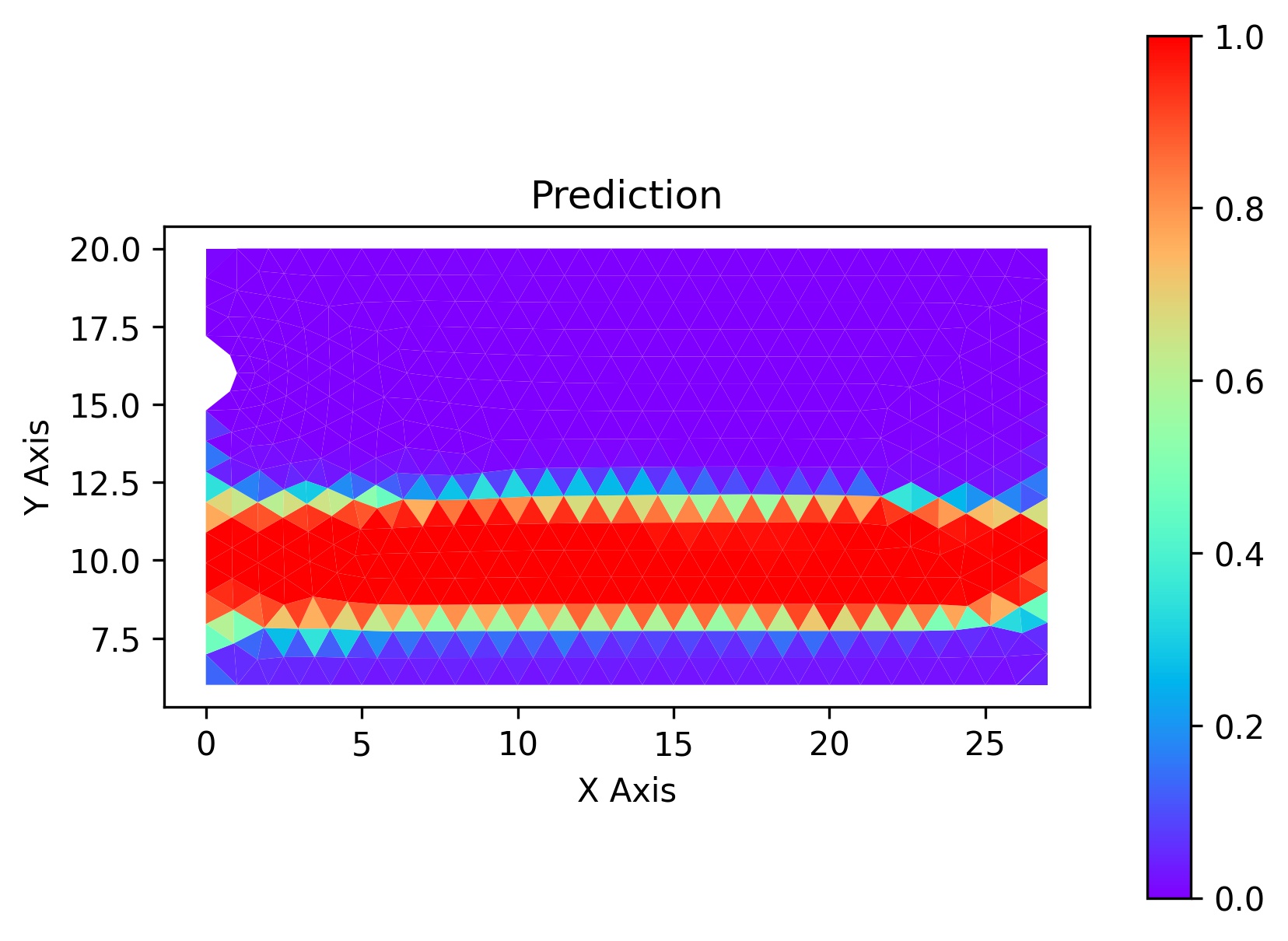}
  \caption{Shape \#\idx prediction}
\end{subfigure}
\begin{subfigure}{\figWidth\textwidth}
  \centering
  \includegraphics[width=\linWidthRatio\linewidth]{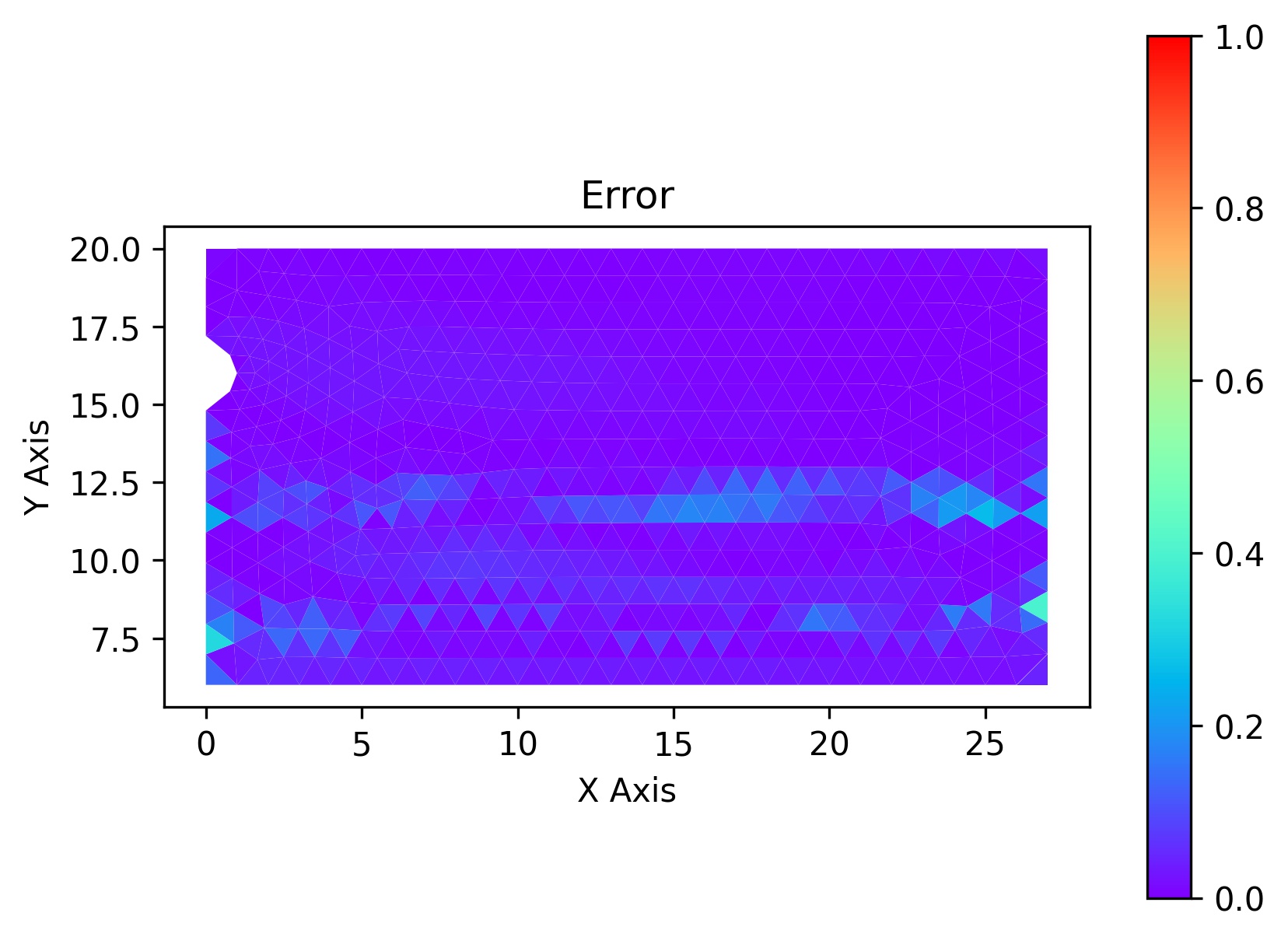}
  \caption{Shape \#\idx error}
\end{subfigure}
\end{figure}

\def \idx{4~} 
\begin{figure}[!htbp]\ContinuedFloat
\begin{subfigure}{\figWidth\textwidth}
  \centering
  \includegraphics[width=\linWidthRatio\linewidth]{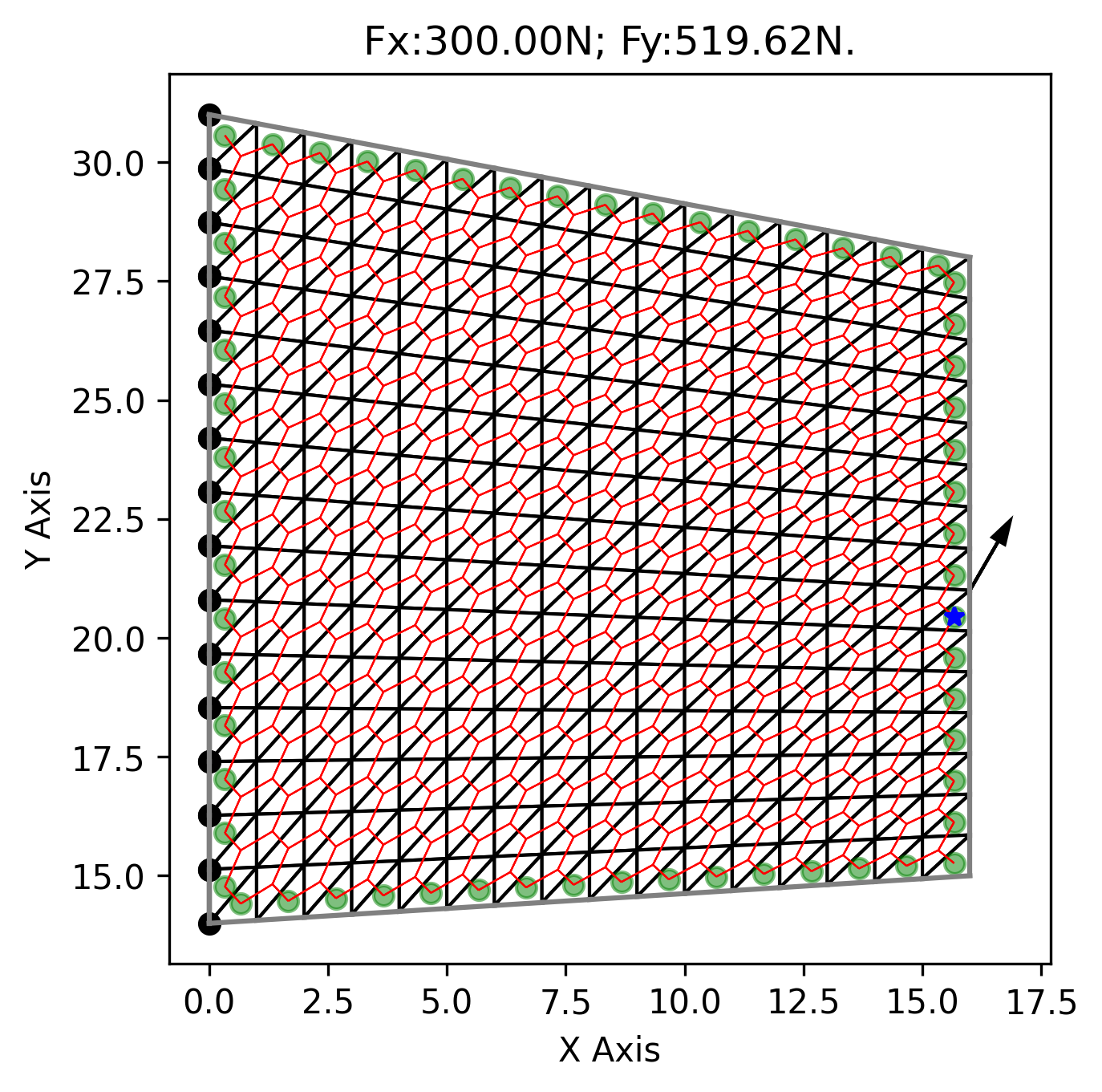}
  \caption{Shape \#\idx simulation settings}
\end{subfigure}%
\begin{subfigure}{\figWidth\textwidth}
  \centering
  \includegraphics[width=\linWidthRatio\linewidth]{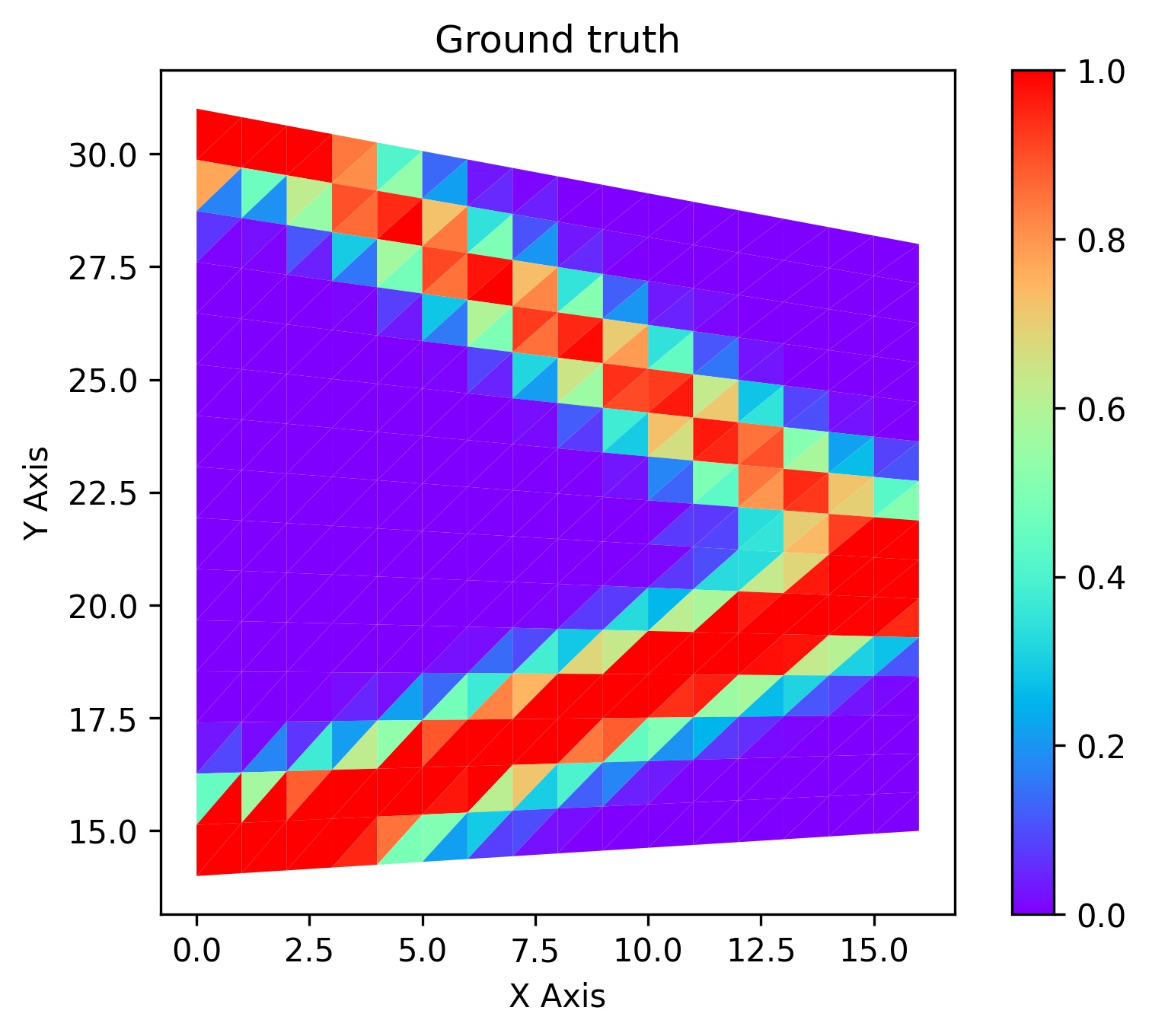}
  \caption{Shape \#\idx ground truth}
\end{subfigure}
\begin{subfigure}{\figWidth\textwidth}
  \centering
  \includegraphics[width=\linWidthRatio\linewidth]{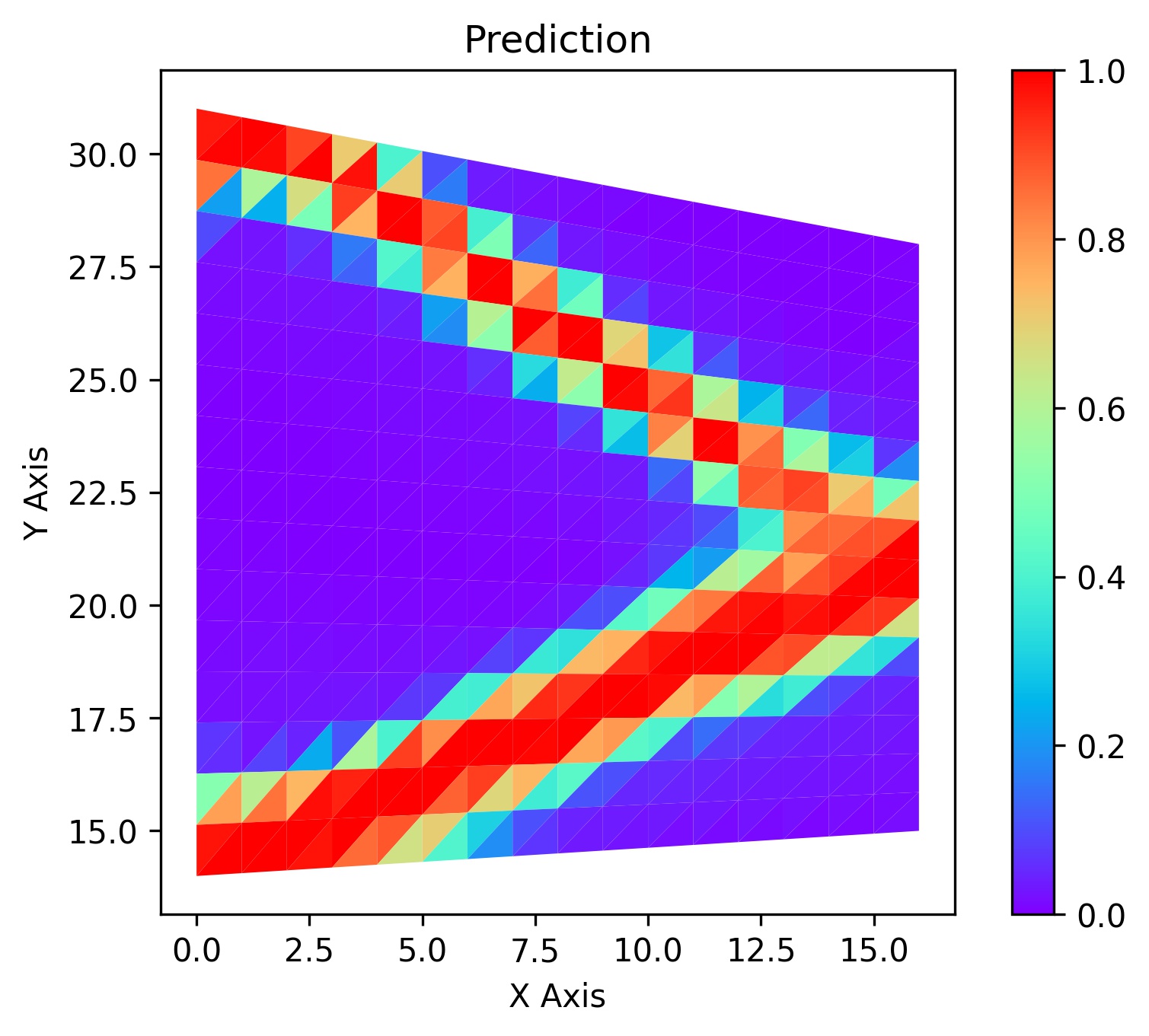}
  \caption{Shape \#\idx prediction}
\end{subfigure}
\begin{subfigure}{\figWidth\textwidth}
  \centering
  \includegraphics[width=\linWidthRatio\linewidth]{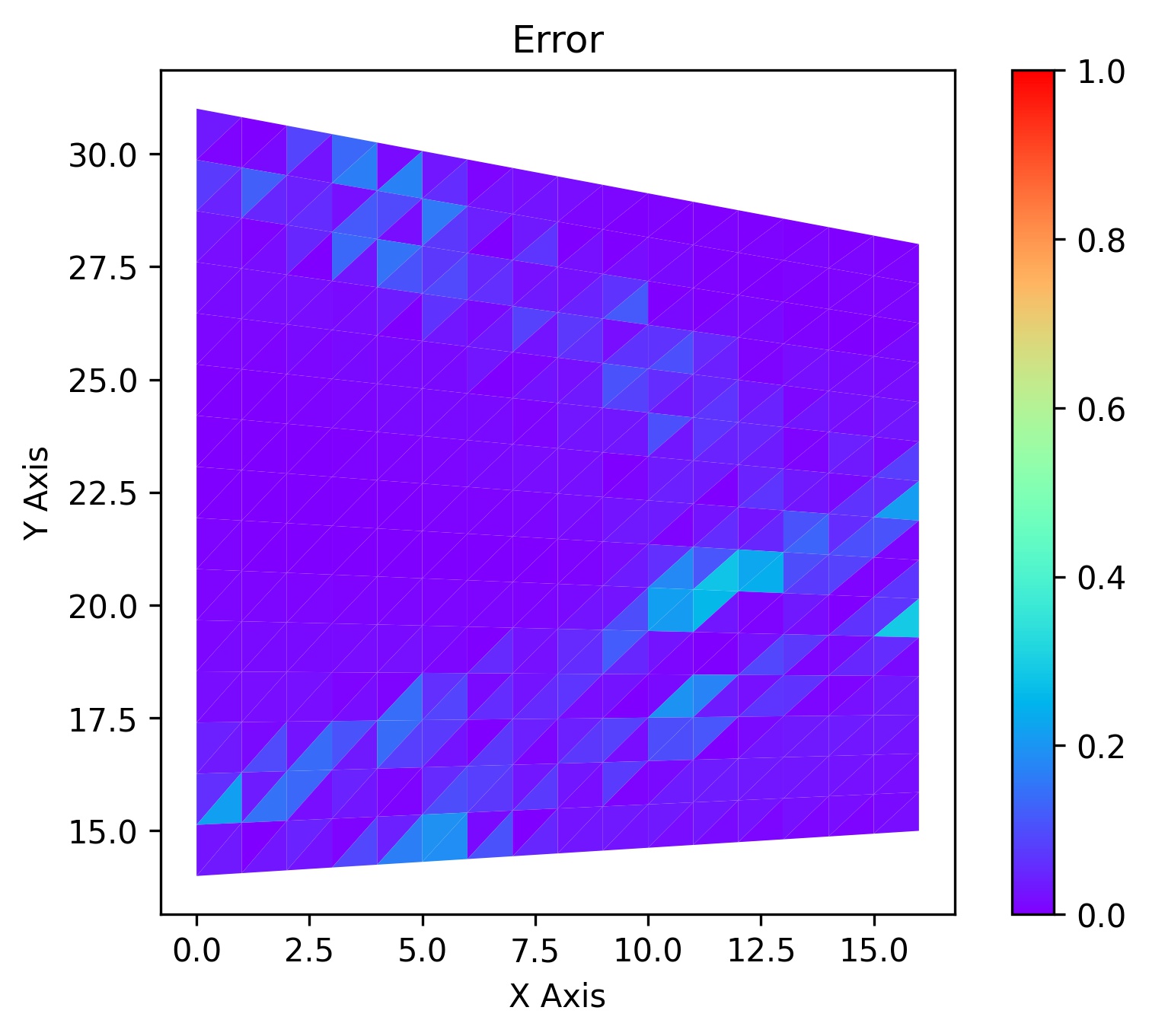}
  \caption{Shape \#\idx error}
\end{subfigure}
\end{figure}

\def \idx{5~} 
\begin{figure}[!htbp]\ContinuedFloat
\begin{subfigure}{\figWidth\textwidth}
  \centering
  \includegraphics[width=\linWidthRatio\linewidth]{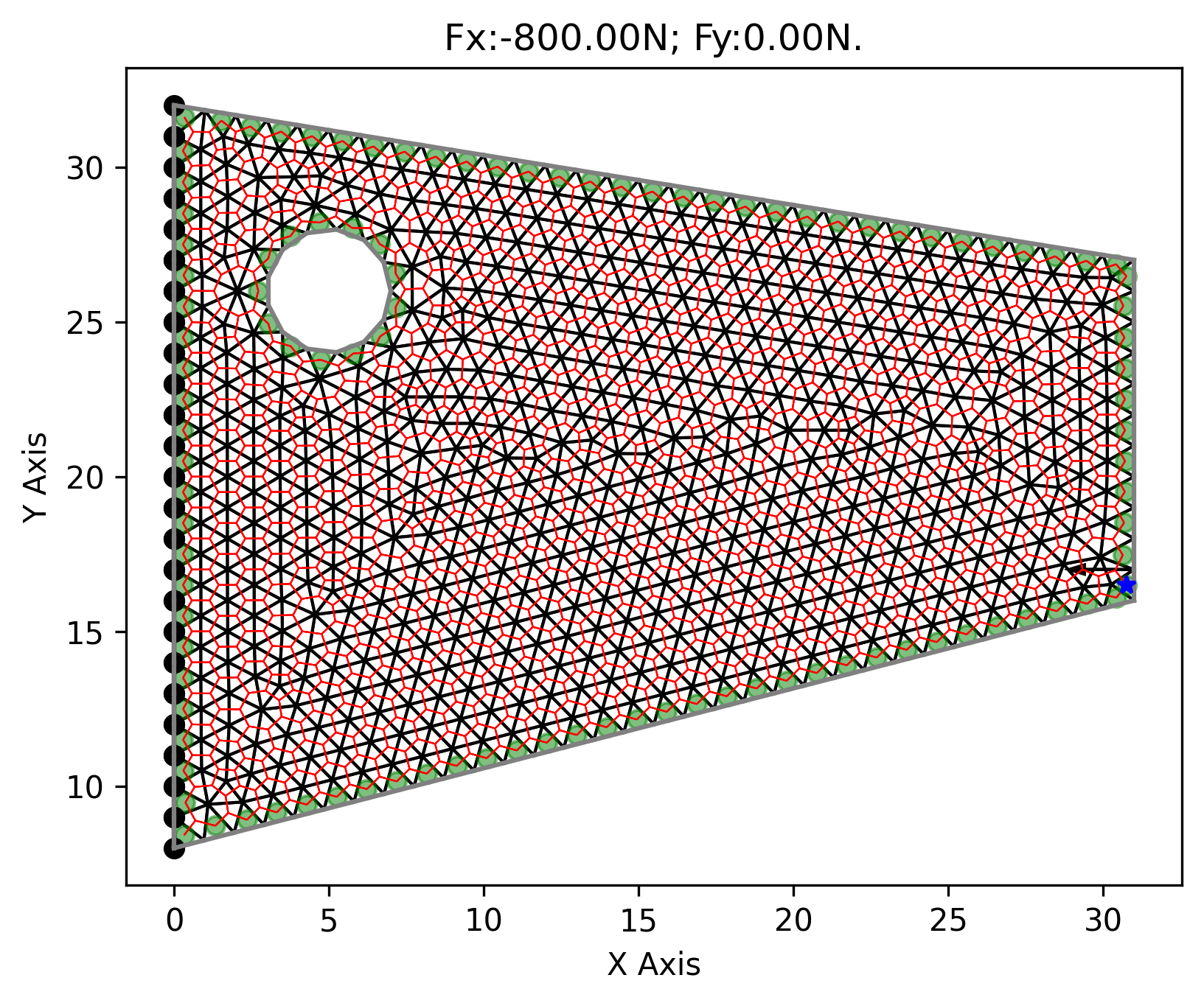}
  \caption{Shape \#\idx simulation settings}
\end{subfigure}%
\begin{subfigure}{\figWidth\textwidth}
  \centering
  \includegraphics[width=\linWidthRatio\linewidth]{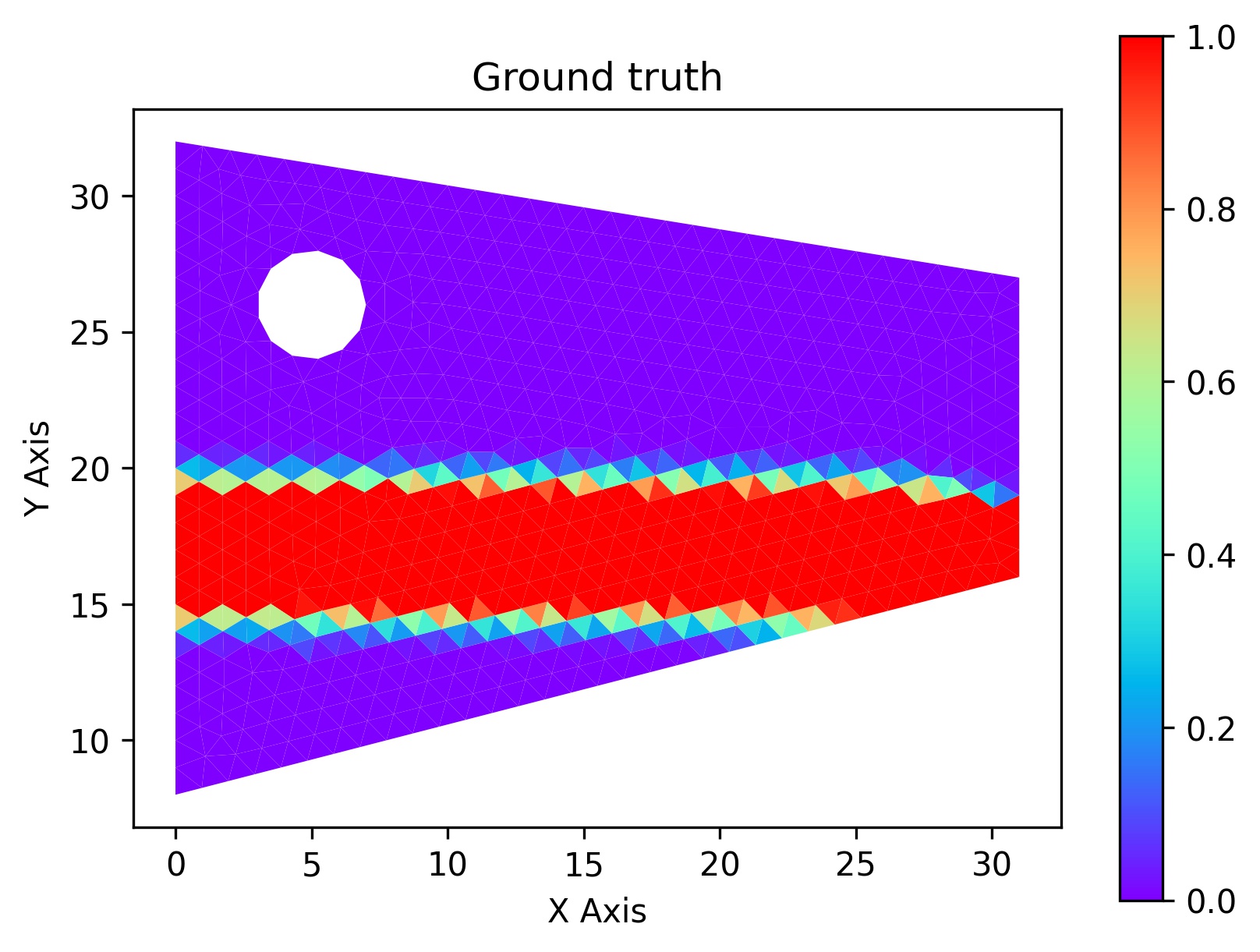}
  \caption{Shape \#\idx ground truth}
\end{subfigure}
\begin{subfigure}{\figWidth\textwidth}
  \centering
  \includegraphics[width=\linWidthRatio\linewidth]{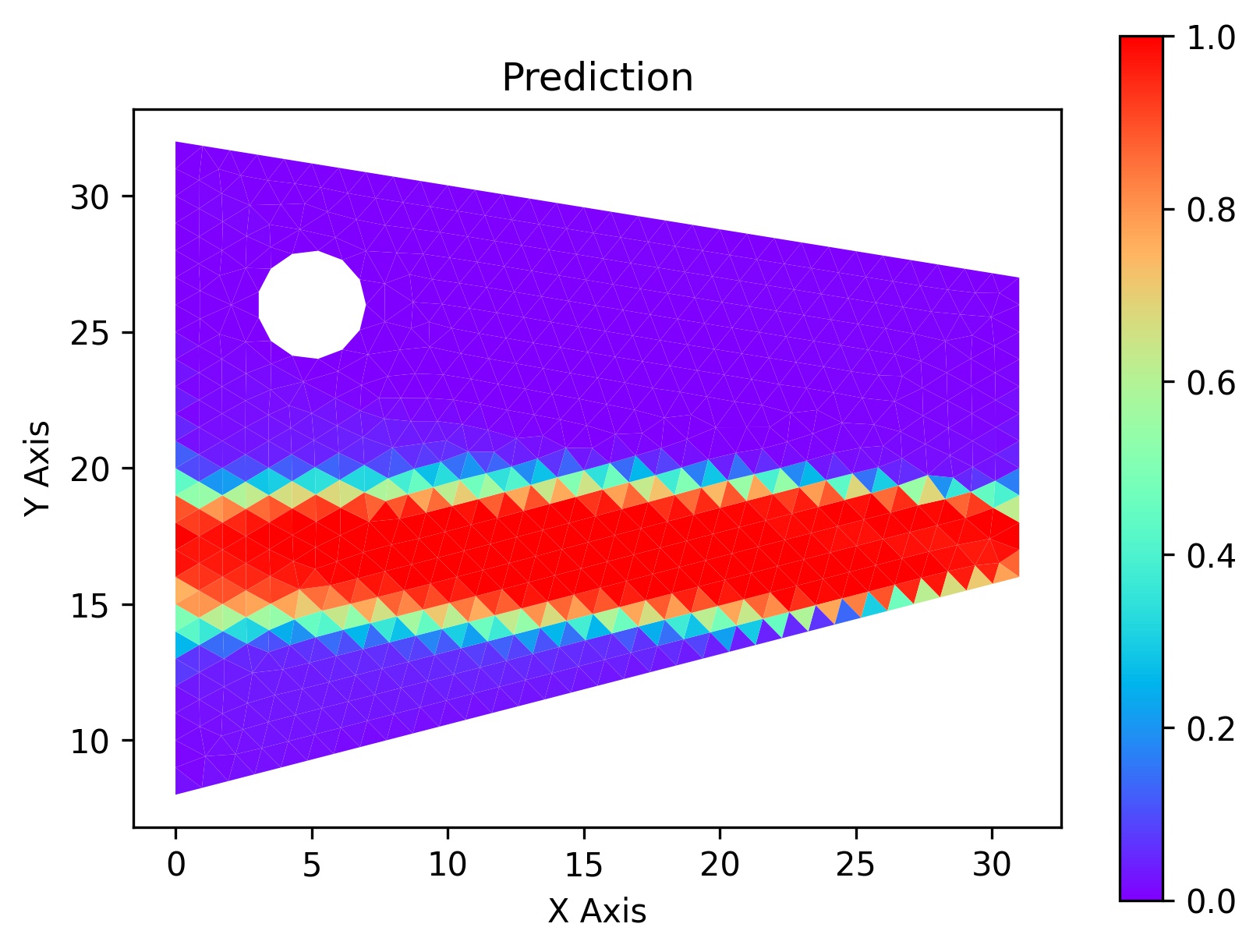}
  \caption{Shape \#\idx prediction}
\end{subfigure}
\begin{subfigure}{\figWidth\textwidth}
  \centering
  \includegraphics[width=\linWidthRatio\linewidth]{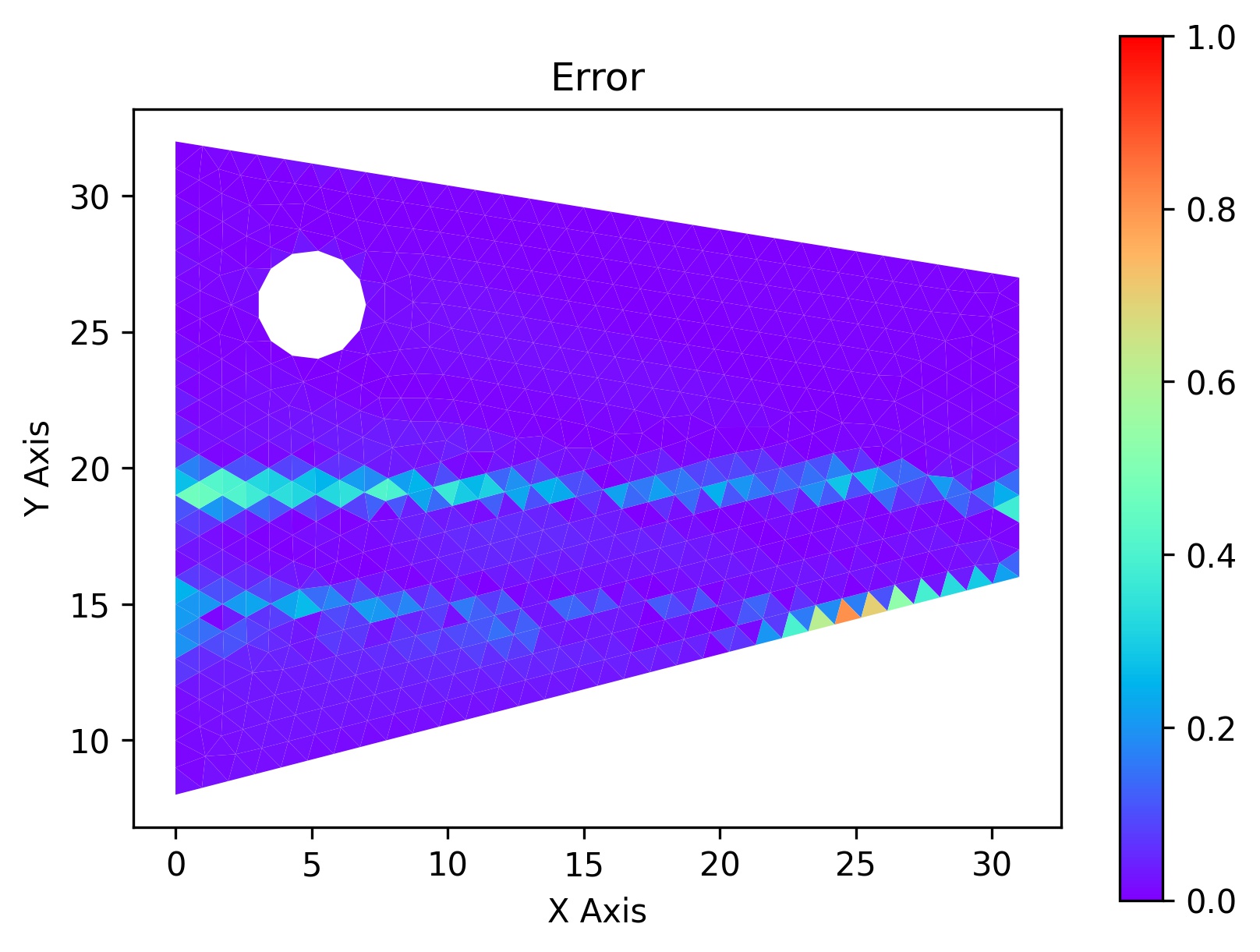}
  \caption{Shape \#\idx error}
\end{subfigure}
\end{figure}

\def \idx{6~} 
\begin{figure}[!htbp]\ContinuedFloat
\begin{subfigure}{\figWidth\textwidth}
  \centering
  \includegraphics[width=\linWidthRatio\linewidth]{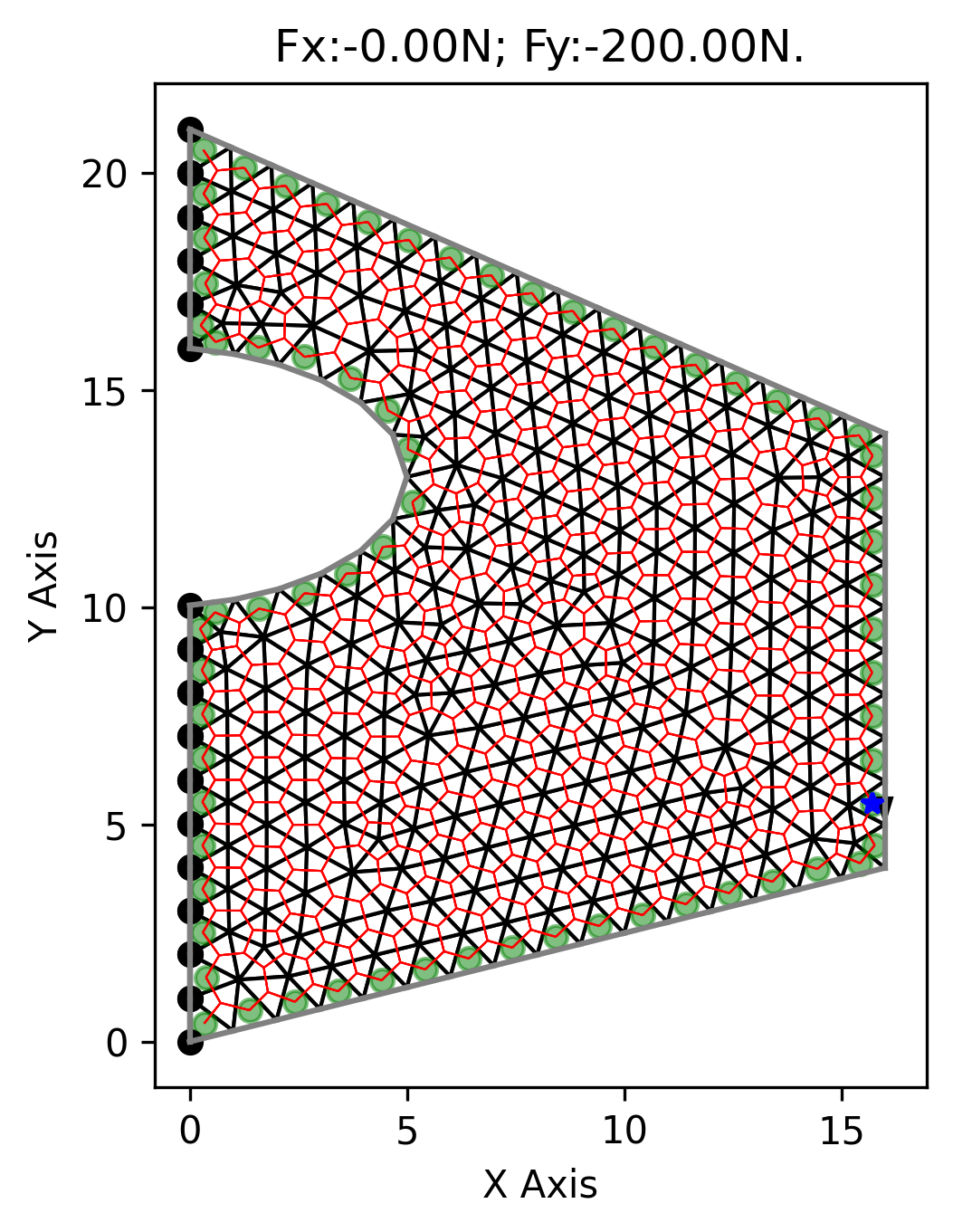}
  \caption{Shape \#\idx simulation settings}
\end{subfigure}%
\begin{subfigure}{\figWidth\textwidth}
  \centering
  \includegraphics[width=\linWidthRatio\linewidth]{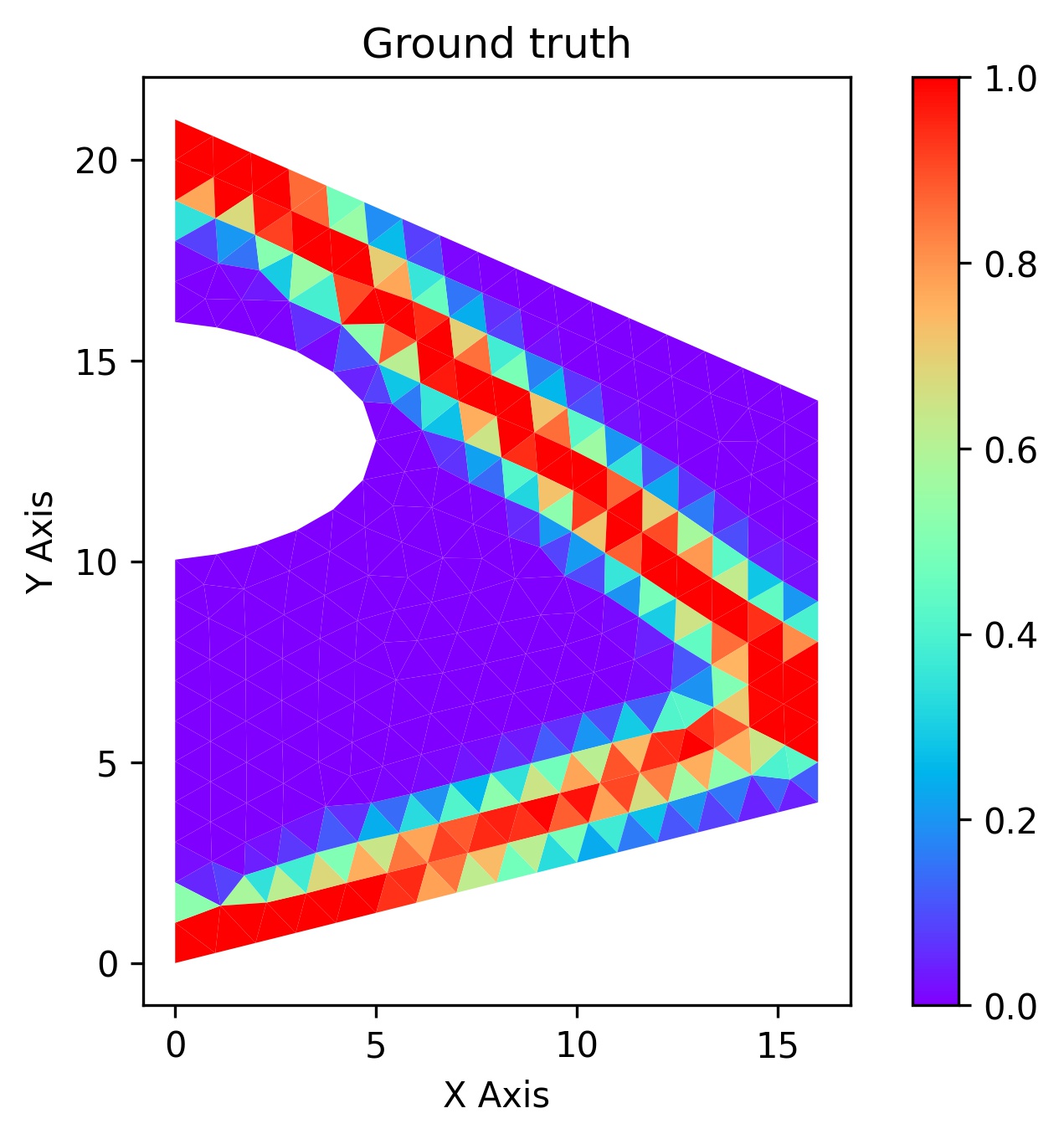}
  \caption{Shape \#\idx ground truth}
\end{subfigure}
\begin{subfigure}{\figWidth\textwidth}
  \centering
  \includegraphics[width=\linWidthRatio\linewidth]{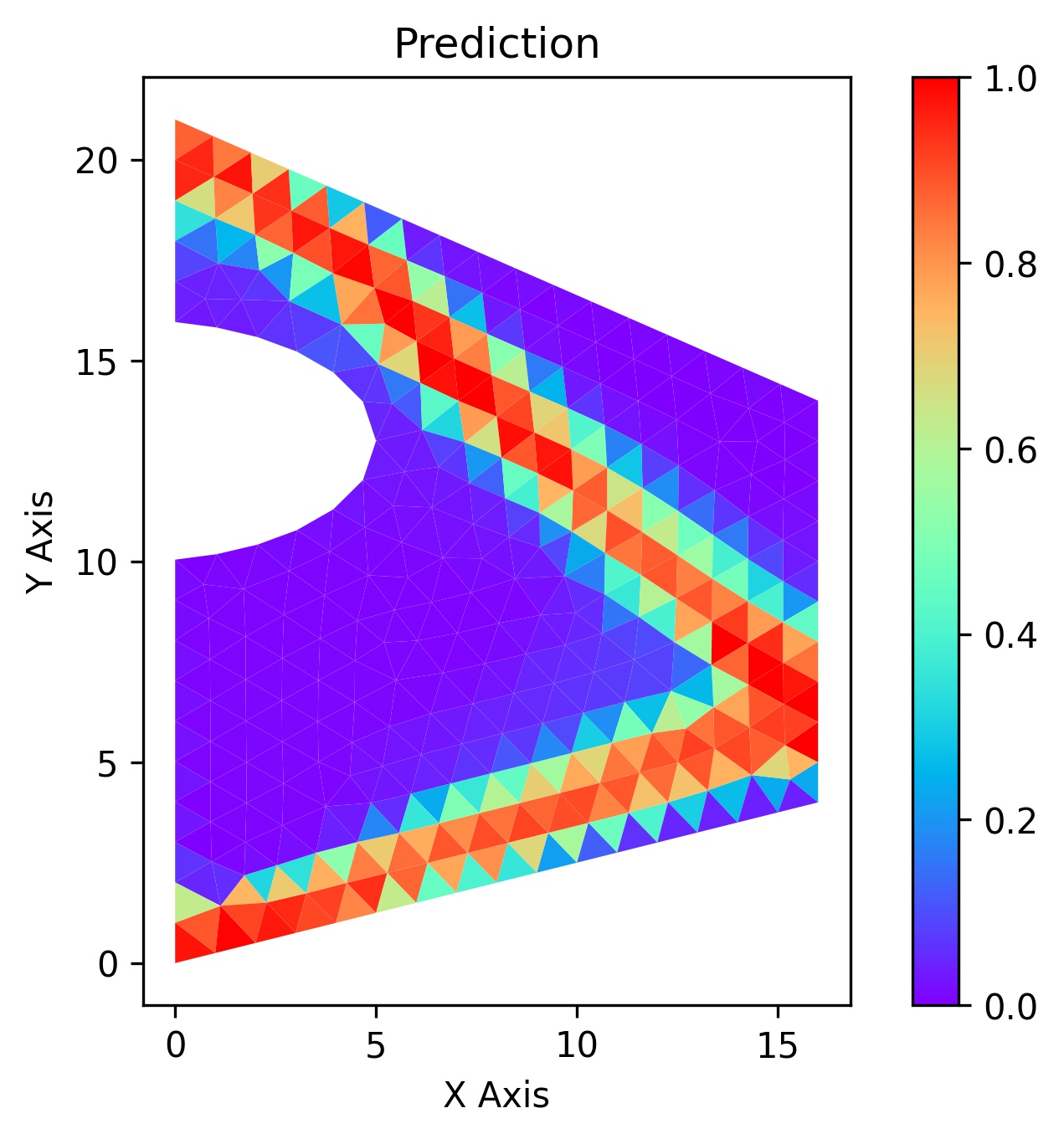}
  \caption{Shape \#\idx prediction}
\end{subfigure}
\begin{subfigure}{\figWidth\textwidth}
  \centering
  \includegraphics[width=\linWidthRatio\linewidth]{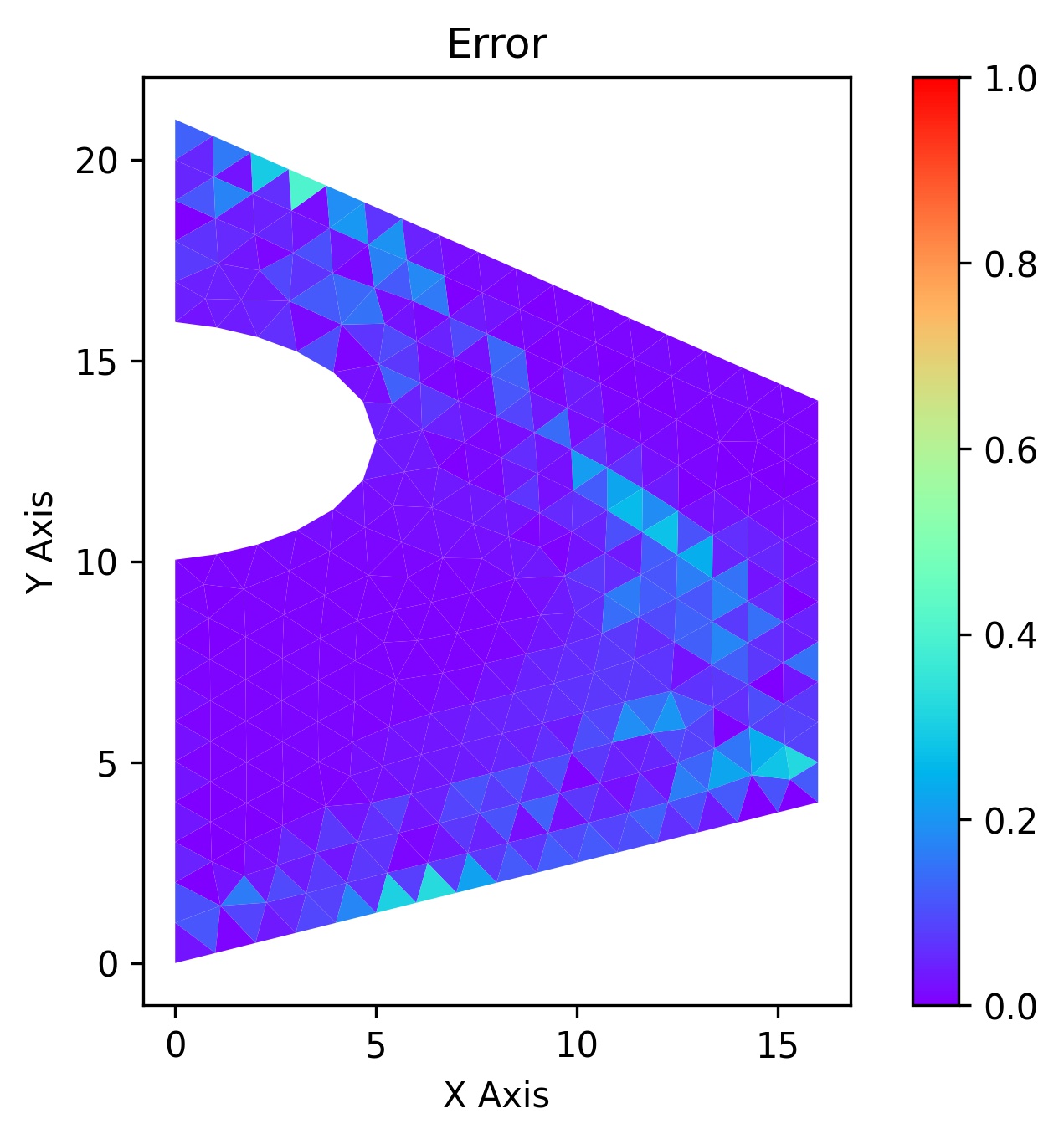}
  \caption{Shape \#\idx error}
\end{subfigure}
\end{figure}

\def \idx{7~} 
\begin{figure}[!htbp]\ContinuedFloat
\begin{subfigure}{\figWidth\textwidth}
  \centering
  \includegraphics[width=\linWidthRatio\linewidth]{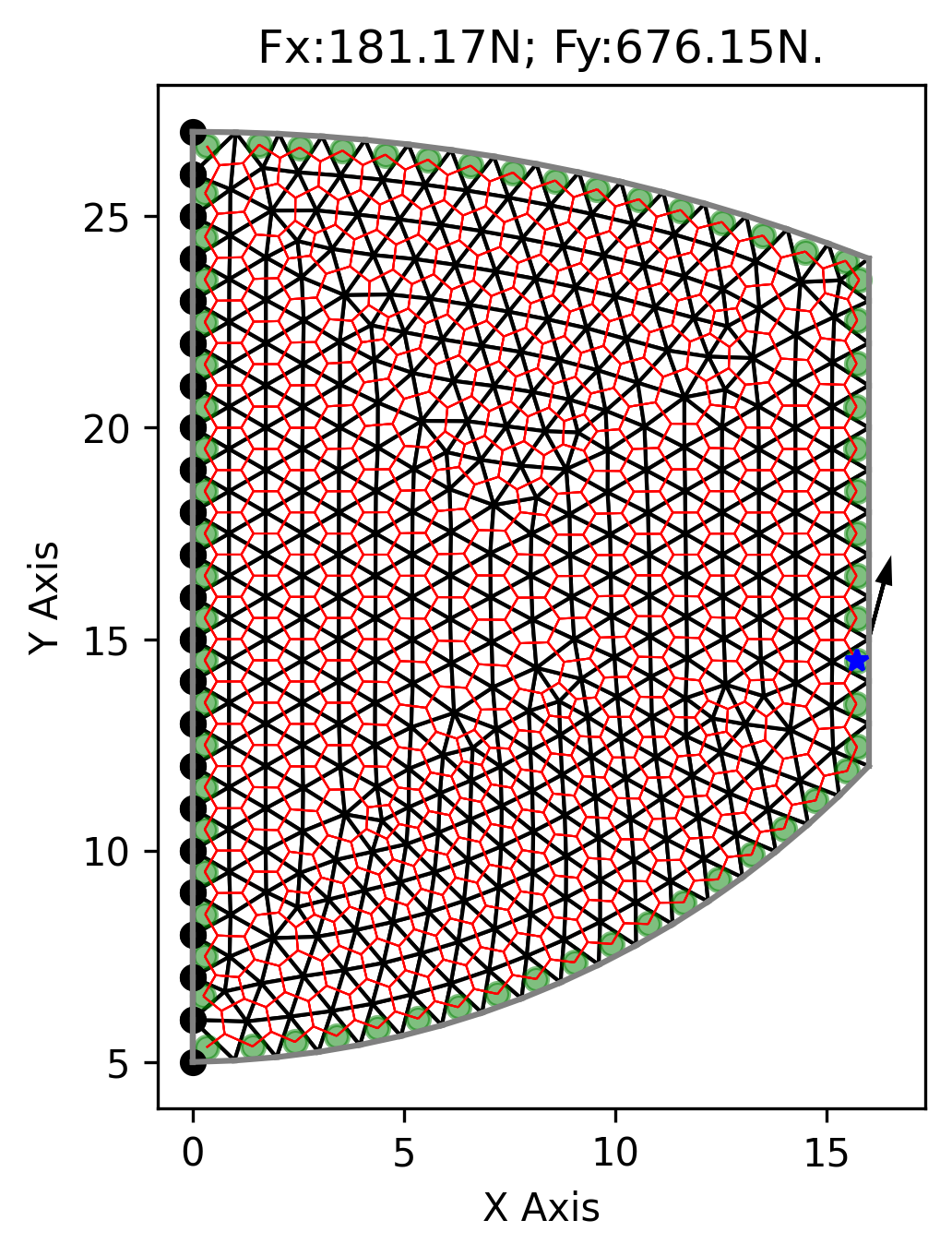}
  \caption{Shape \#\idx simulation settings}
\end{subfigure}%
\begin{subfigure}{\figWidth\textwidth}
  \centering
  \includegraphics[width=\linWidthRatio\linewidth]{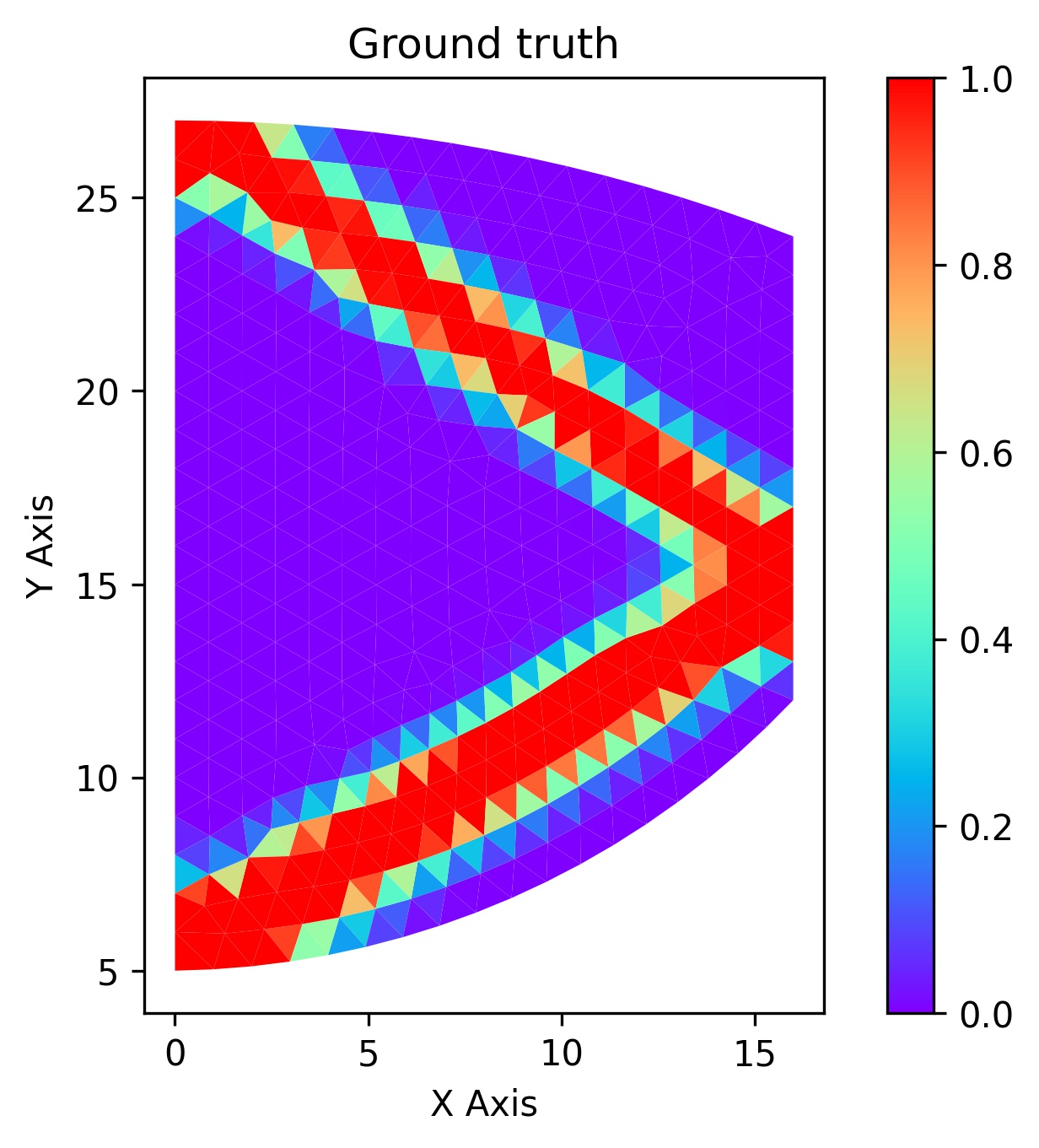}
  \caption{Shape \#\idx ground truth}
\end{subfigure}
\begin{subfigure}{\figWidth\textwidth}
  \centering
  \includegraphics[width=\linWidthRatio\linewidth]{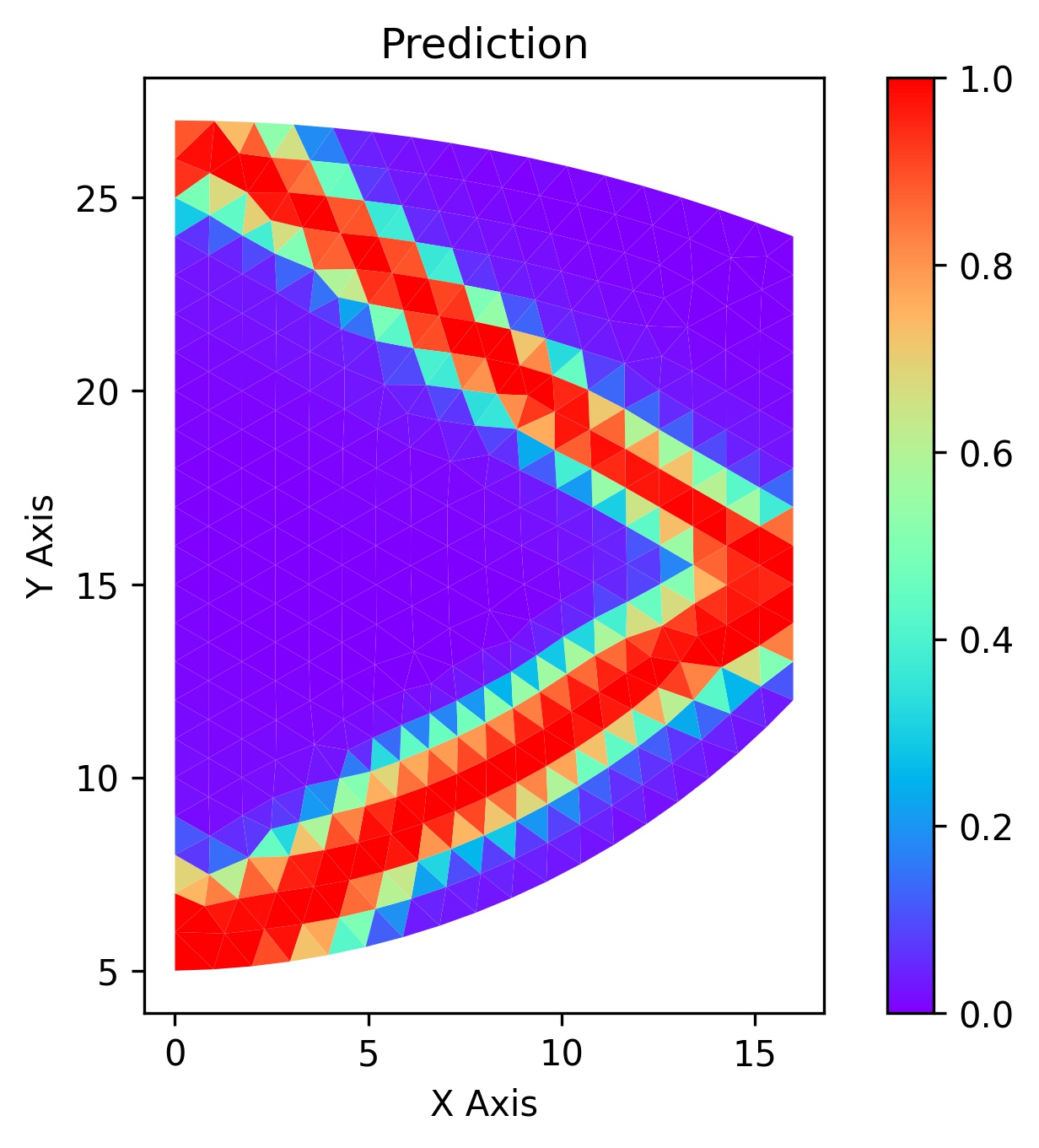}
  \caption{Shape \#\idx prediction}
\end{subfigure}
\begin{subfigure}{\figWidth\textwidth}
  \centering
  \includegraphics[width=\linWidthRatio\linewidth]{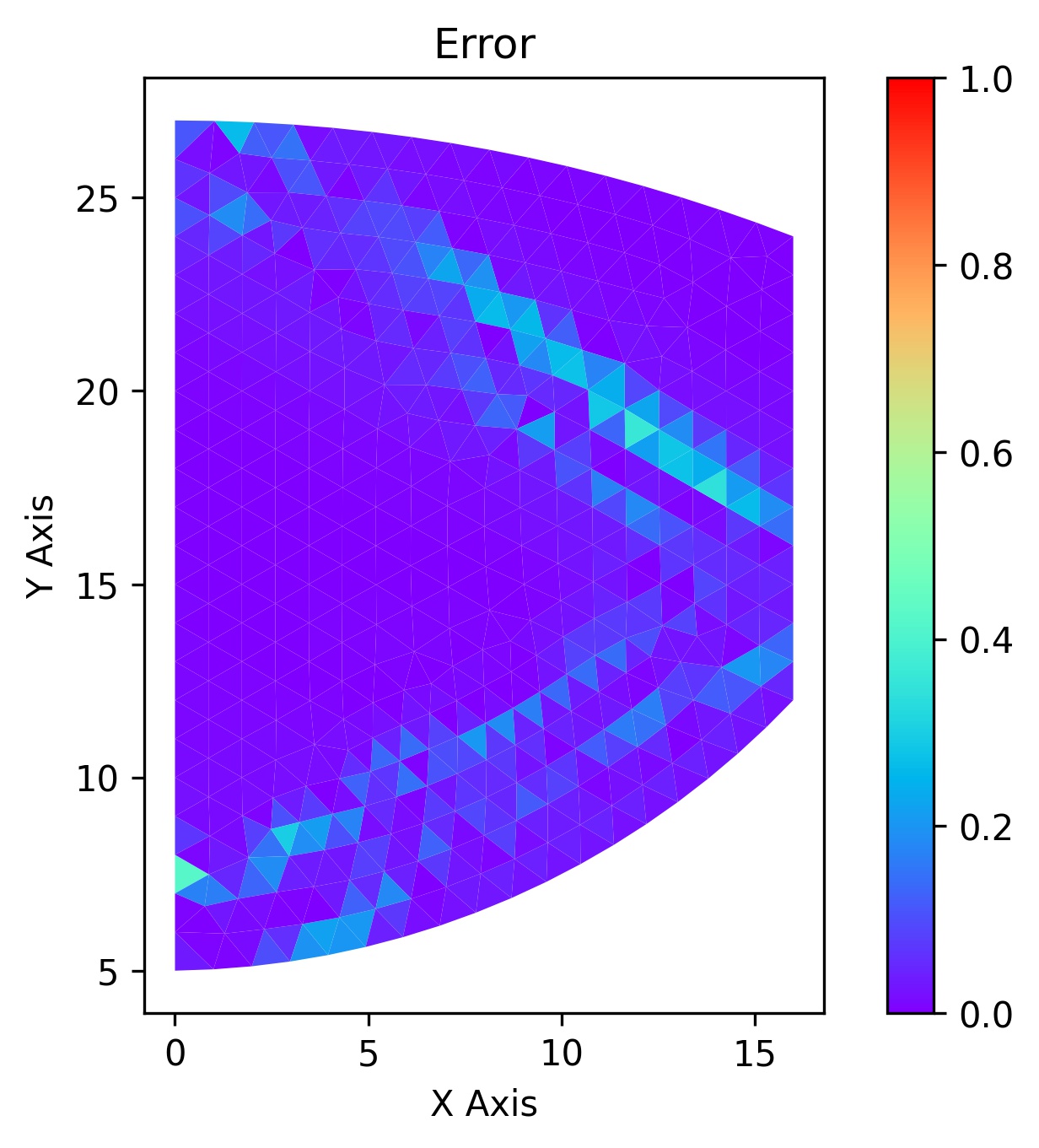}
  \caption{Shape \#\idx error}
\end{subfigure}
\end{figure}

\def \idx{8~} 
\begin{figure}[!htbp]\ContinuedFloat
\begin{subfigure}{\figWidth\textwidth}
  \centering
  \includegraphics[width=\linWidthRatio\linewidth]{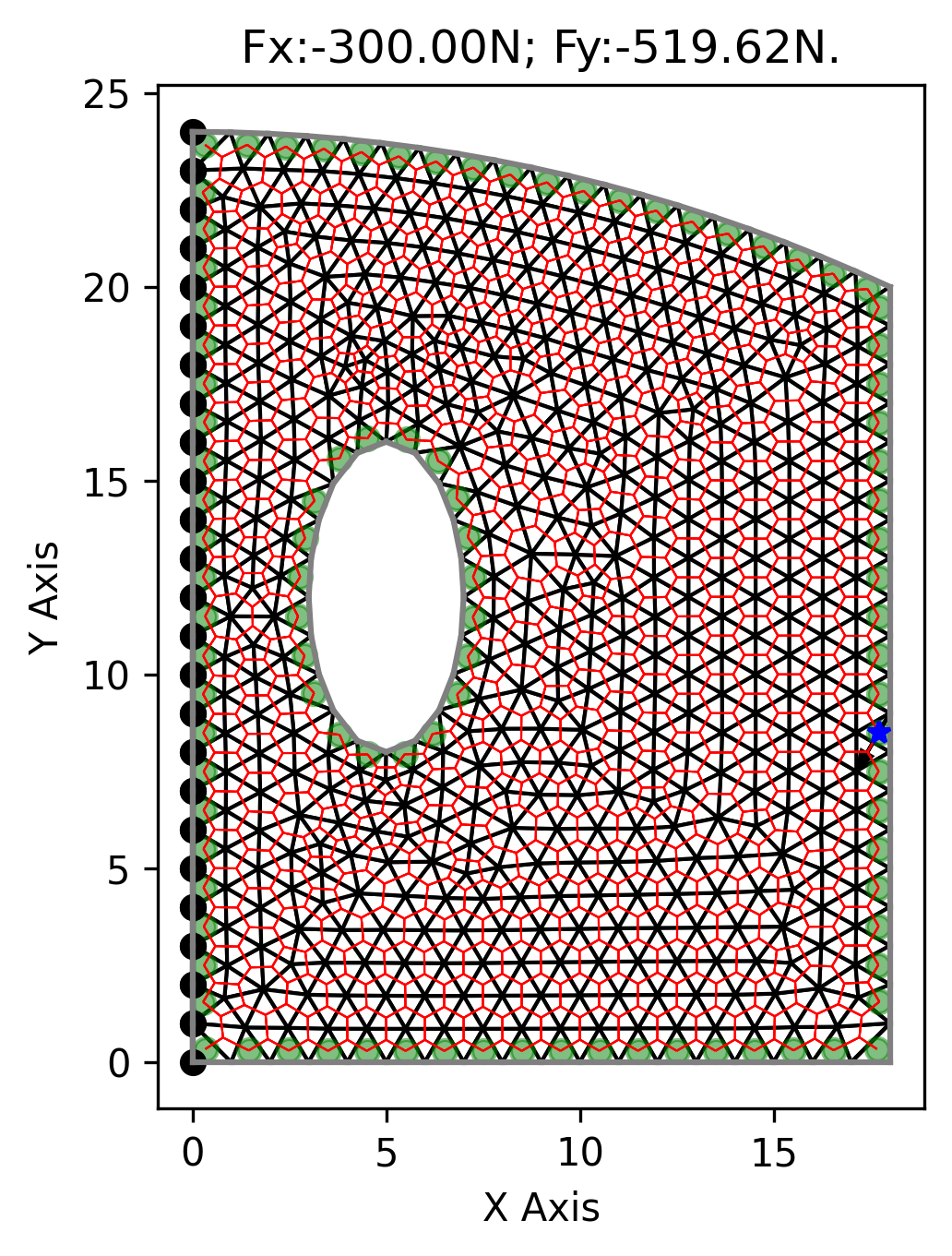}
  \caption{Shape \#\idx simulation settings}
\end{subfigure}%
\begin{subfigure}{\figWidth\textwidth}
  \centering
  \includegraphics[width=\linWidthRatio\linewidth]{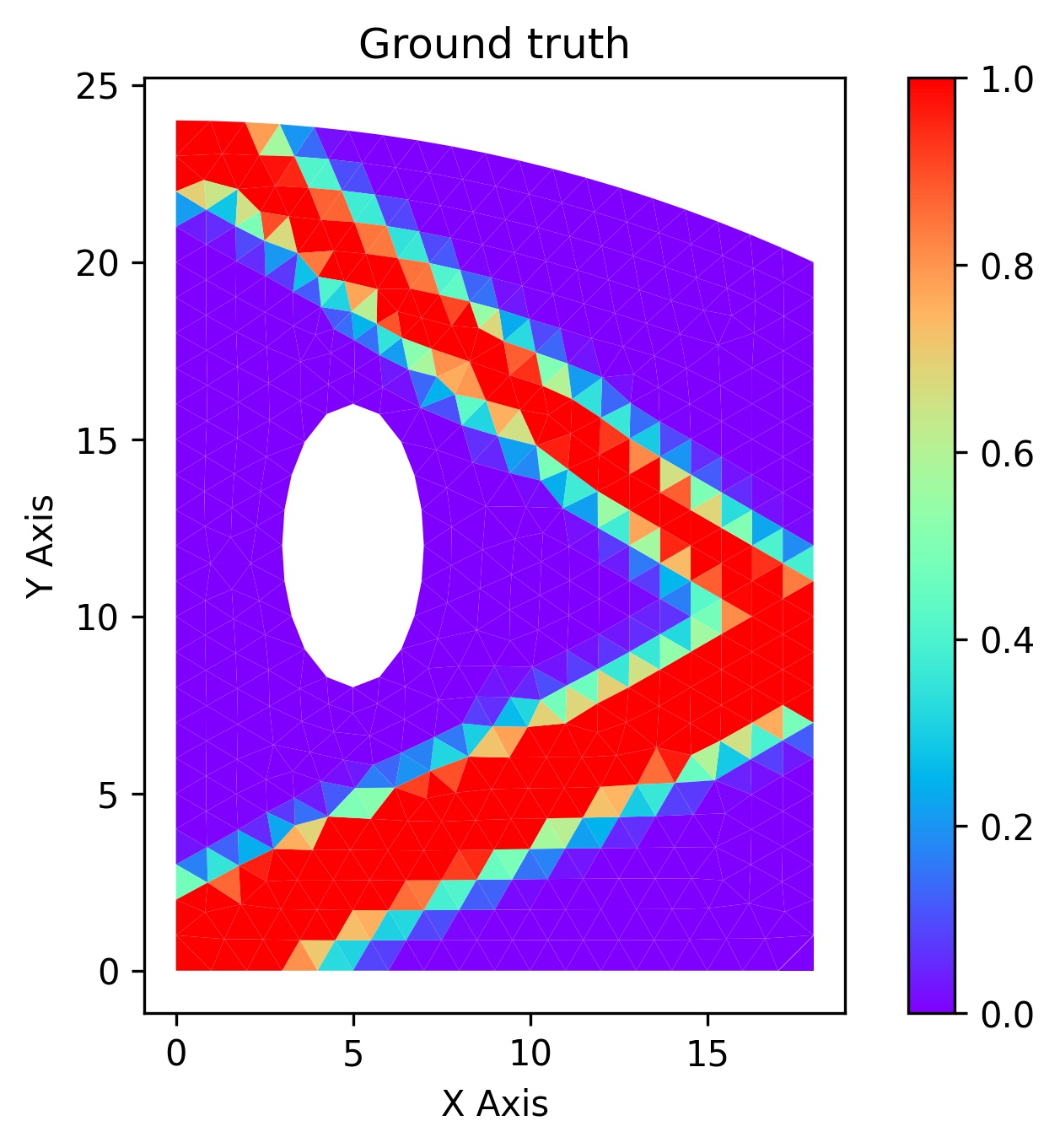}
  \caption{Shape \#\idx ground truth}
\end{subfigure}
\begin{subfigure}{\figWidth\textwidth}
  \centering
  \includegraphics[width=\linWidthRatio\linewidth]{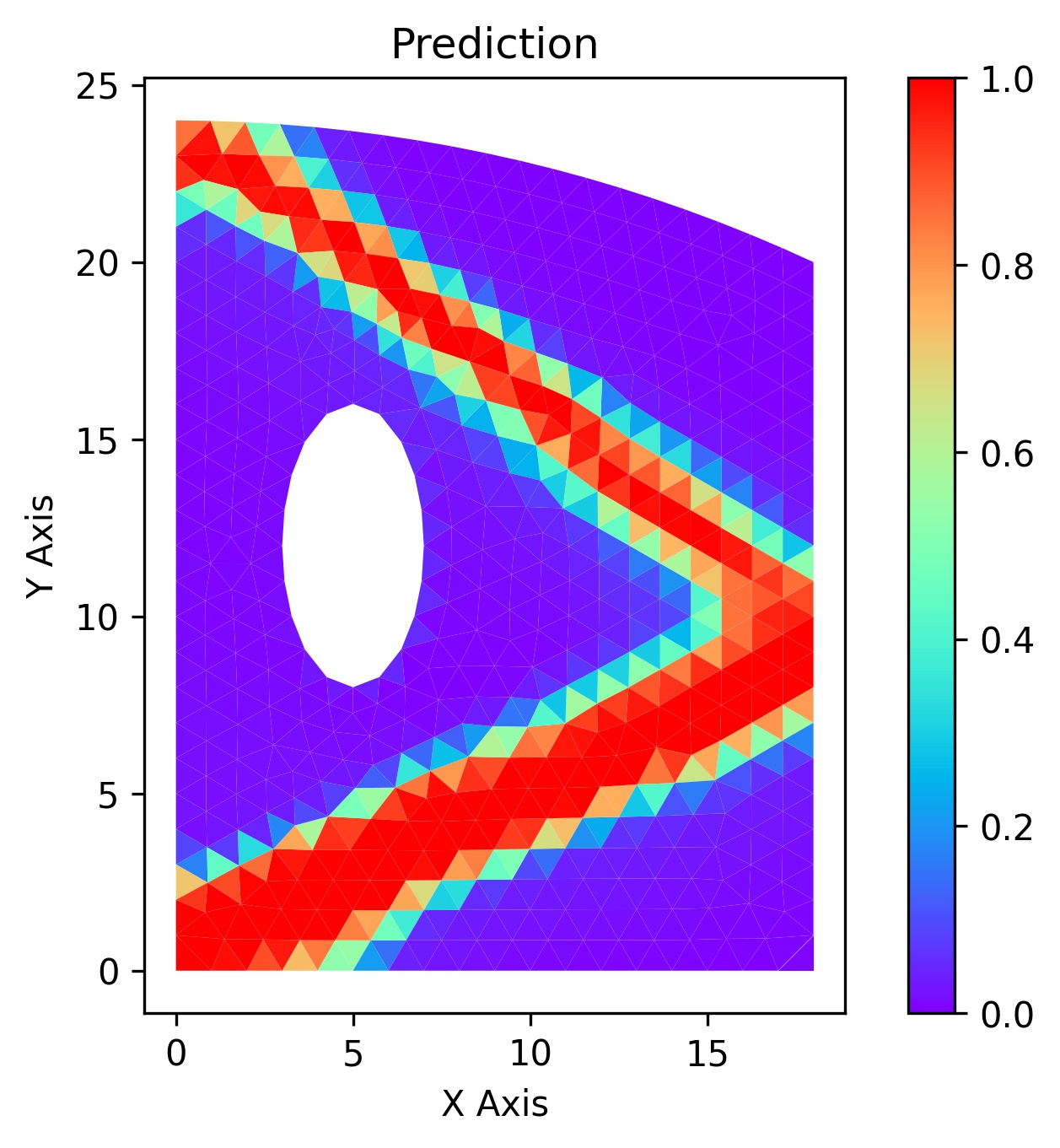}
  \caption{Shape \#\idx prediction}
\end{subfigure}
\begin{subfigure}{\figWidth\textwidth}
  \centering
  \includegraphics[width=\linWidthRatio\linewidth]{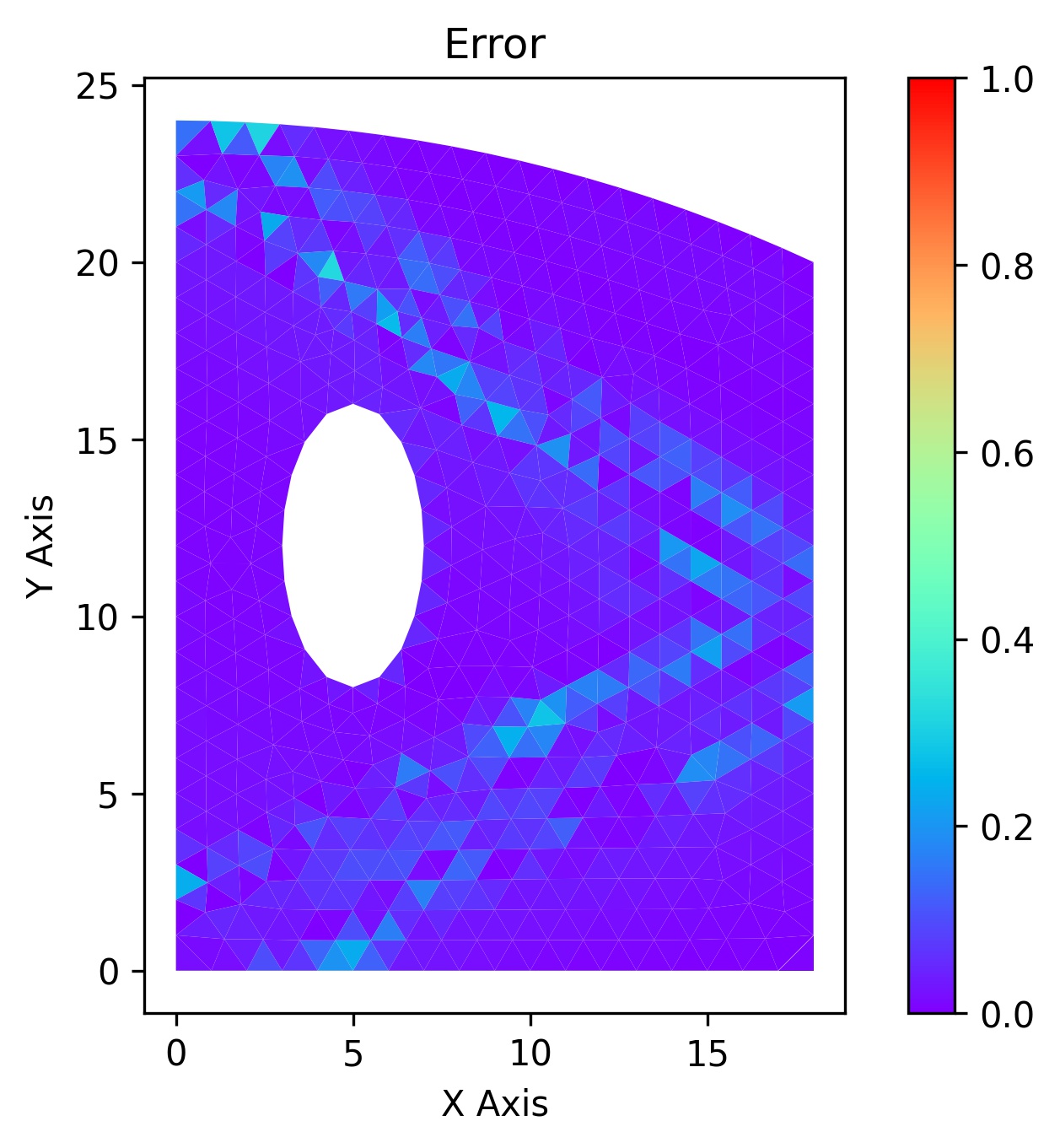}
  \caption{Shape \#\idx error}
\end{subfigure}
\end{figure}

\def \idx{9~} 
\begin{figure}[!htbp]\ContinuedFloat
\begin{subfigure}{\figWidth\textwidth}
  \centering
  \includegraphics[width=\linWidthRatio\linewidth]{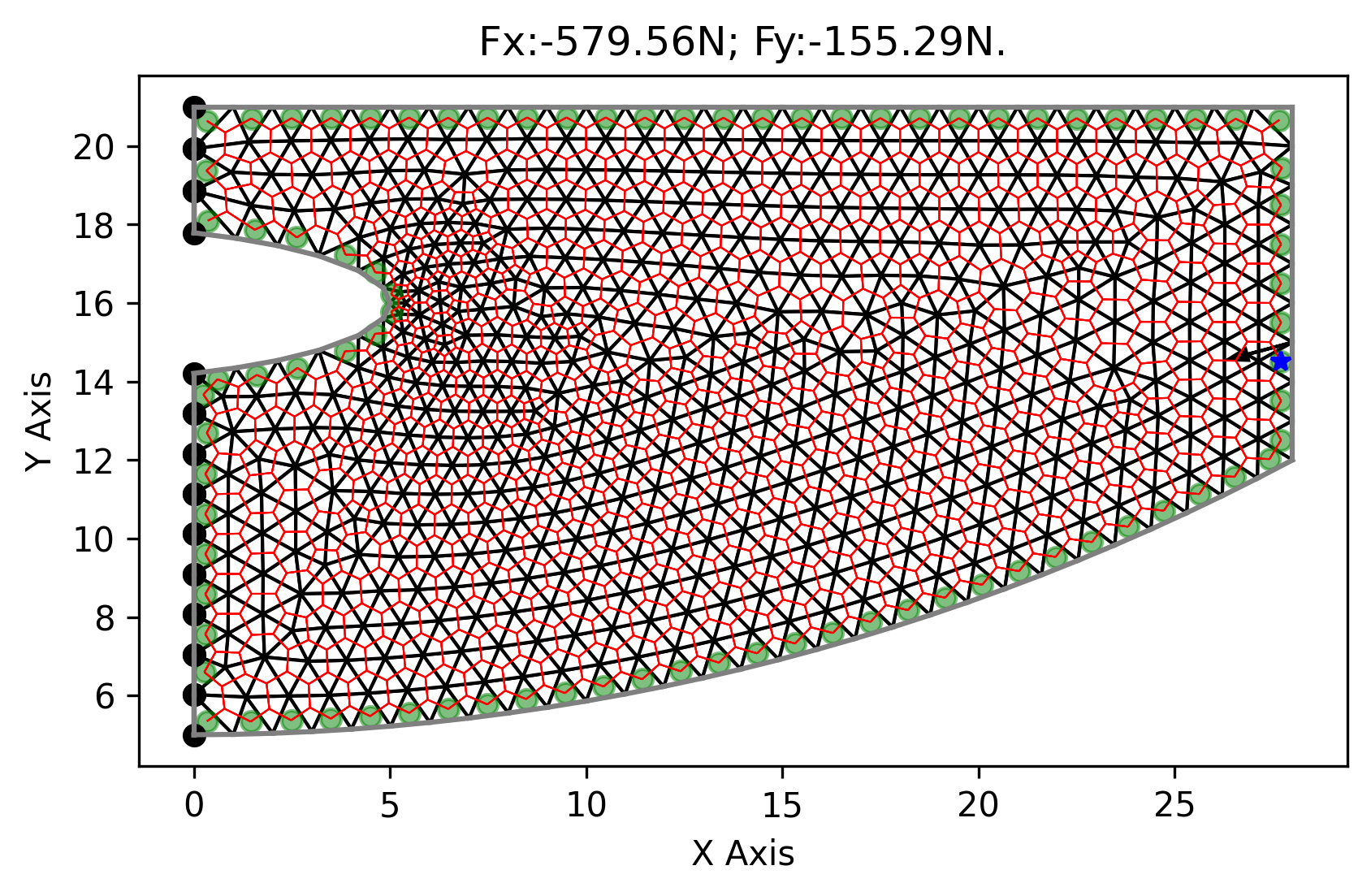}
  \caption{Shape \#\idx simulation settings}
\end{subfigure}%
\begin{subfigure}{\figWidth\textwidth}
  \centering
  \includegraphics[width=\linWidthRatio\linewidth]{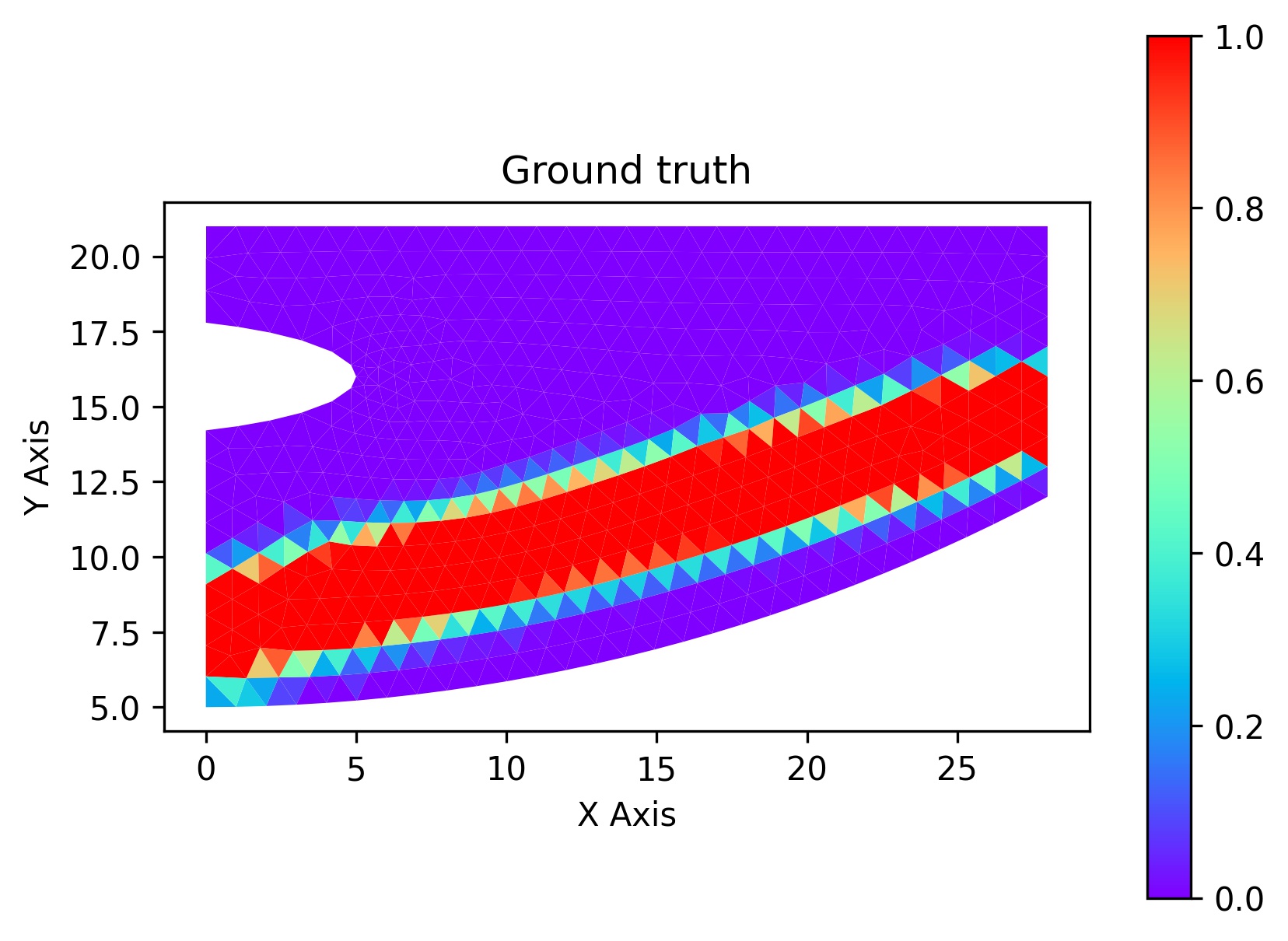}
  \caption{Shape \#\idx ground truth}
\end{subfigure}
\begin{subfigure}{\figWidth\textwidth}
  \centering
  \includegraphics[width=\linWidthRatio\linewidth]{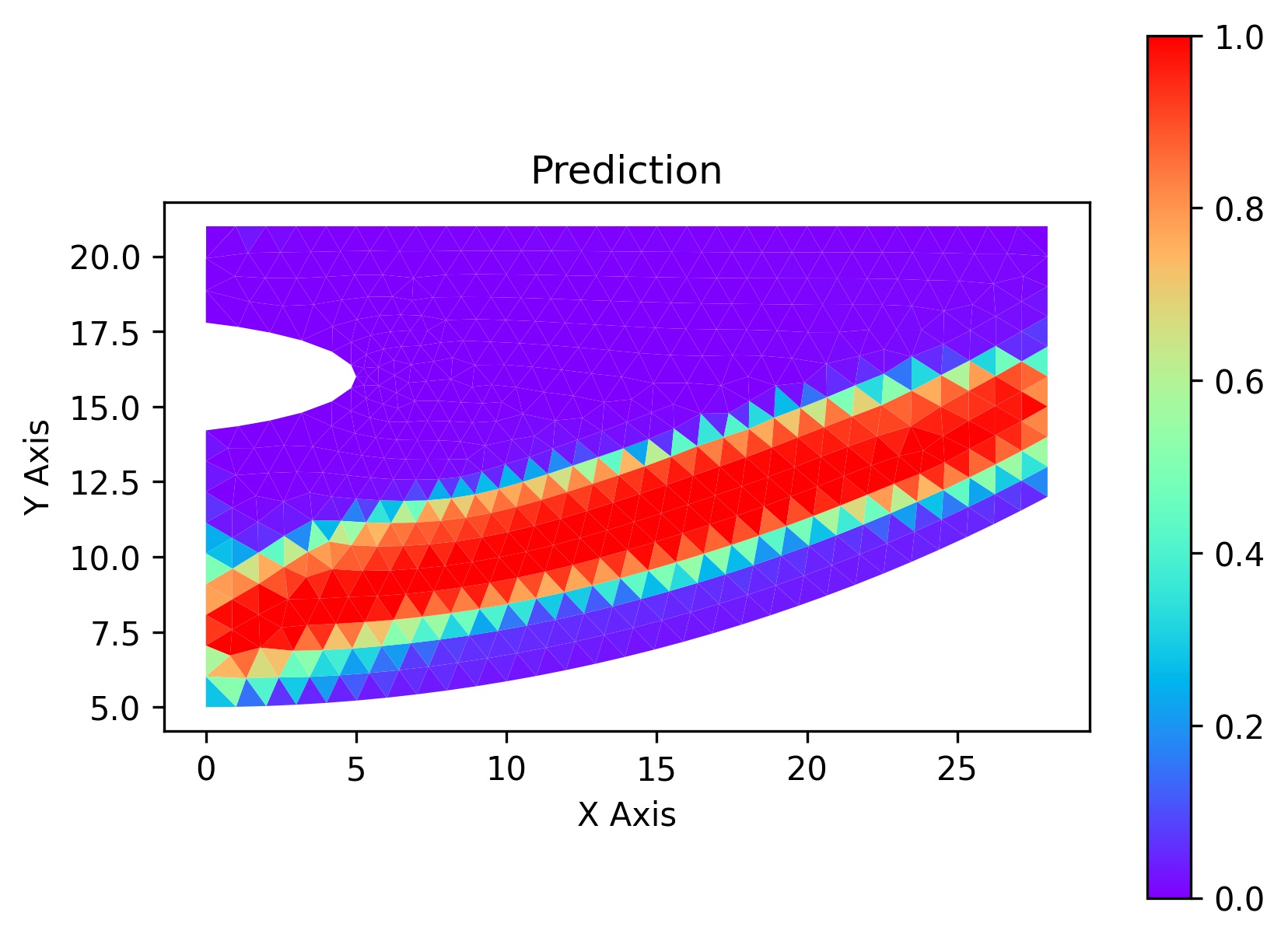}
  \caption{Shape \#\idx prediction}
\end{subfigure}
\begin{subfigure}{\figWidth\textwidth}
  \centering
  \includegraphics[width=\linWidthRatio\linewidth]{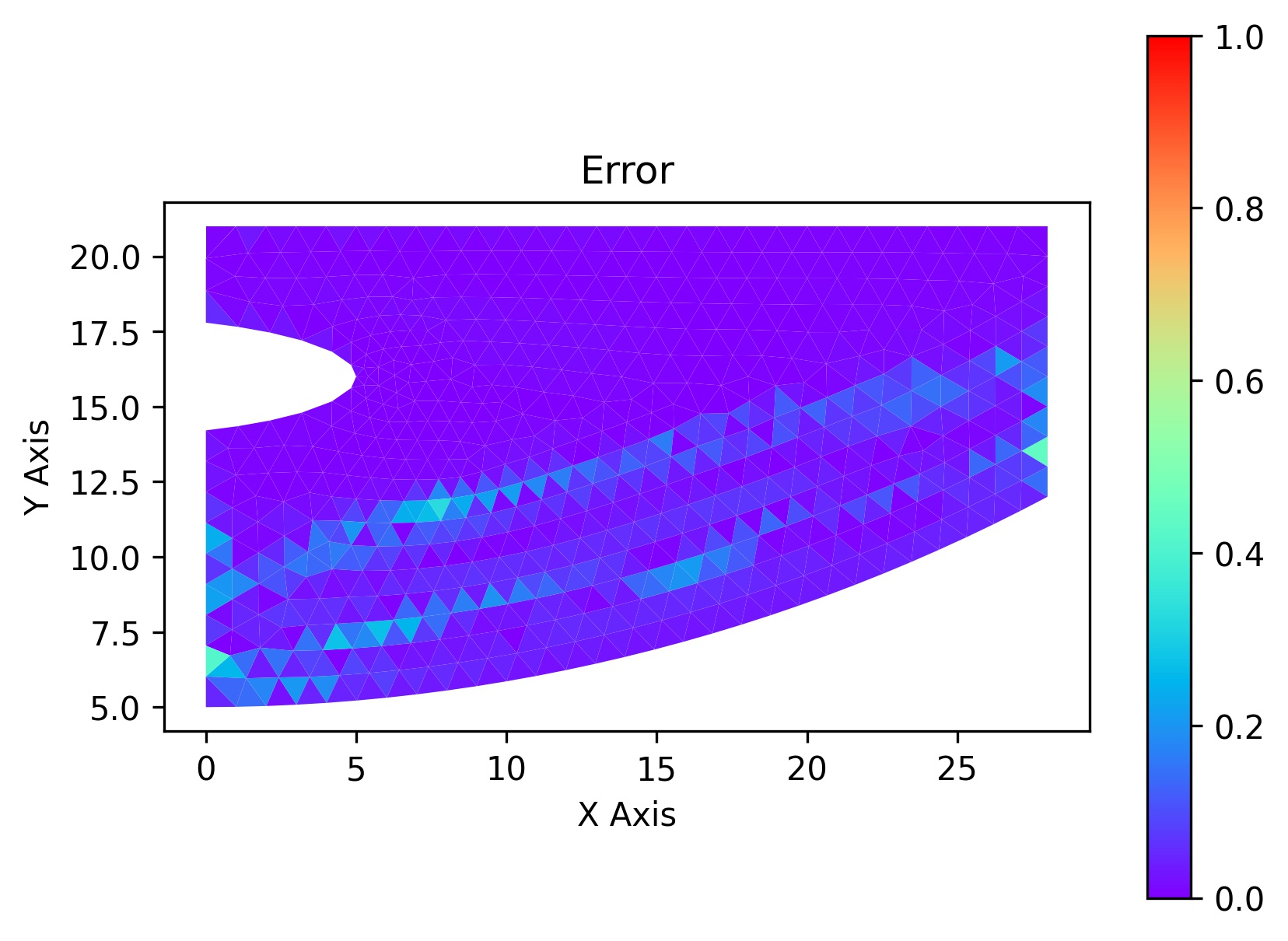}
  \caption{Shape \#\idx error}
\end{subfigure}
\caption{Prediction results for topology optimizations}
\label{fig:TopoResults}
\end{figure}

\section{Conclusion}\label{conclusion}

We develop a Boundary Oriented Graph Embedding (BOGE) approach for the GNN-based FEA surrogate model, especially for regressing triangular-mesh-based FEA simulation results. The BOGE approach bypasses the limitations of the MPNN framework and provides shortcuts for both global boundary information and local multi-hop elements, which allows shallow-layer GNN to regress boundary value problems with long-range graph-vertex interactions. Working on solving the cantilever beam problem, the BOGE approach with 3-layer DeepGCN model \textcolor{blue}{achieves the regression with MSE of 0.011706 (2.41\% MAPE) for stress field prediction and 0.002735 MSE (with 1.58\% elements having error larger than 0.01) for topological optimization.}, which validates the BOGE's efficiency in complex decision-making problems. The model has the potential to be applied to most physical-field fitting problems with any type of polygon meshes or mixed types of meshes. This study has potential limitations. The GNN performance is sensitive to the input-tensor shape and requires more theoretical analysis to investigate appropriate methods to improve the structure of GNN. Future work can focus on modifying GNN structures to obtain more accurate regression results and applying the BOGE approach to other boundary value problems.

\bibliographystyle{unsrt}  

\bibliography{refs}

\begin{thebibliography}{10}

\bibitem{axelsson2001finite}
Owe Axelsson and Vincent~Allan Barker.
\newblock {\em Finite element solution of boundary value problems: theory and
  computation}.
\newblock SIAM, 2001.

\bibitem{capuano2019smart}
German Capuano and Julian~J Rimoli.
\newblock Smart finite elements: A novel machine learning application.
\newblock {\em Computer Methods in Applied Mechanics and Engineering},
  345:363--381, 2019.

\bibitem{tamaddon2020data}
Hamid~Reza Tamaddon-Jahromi, Neeraj~Kavan Chakshu, Igor Sazonov, Llion~M Evans,
  Hywel Thomas, and Perumal Nithiarasu.
\newblock Data-driven inverse modelling through neural network (deep learning)
  and computational heat transfer.
\newblock {\em Computer Methods in Applied Mechanics and Engineering},
  369:113217, 2020.

\bibitem{nie2020stress}
Zhenguo Nie, Haoliang Jiang, and Levent~Burak Kara.
\newblock Stress field prediction in cantilevered structures using
  convolutional neural networks.
\newblock {\em Journal of Computing and Information Science in Engineering},
  20(1), 2020.

\bibitem{kantzos2019design}
Christopher Kantzos, Jacky Lao, and Anthony Rollett.
\newblock Design of an interpretable convolutional neural network for stress
  concentration prediction in rough surfaces.
\newblock {\em Materials Characterization}, 158:109961, 2019.

\bibitem{khadilkar2019deep}
Aditya Khadilkar, Jun Wang, and Rahul Rai.
\newblock Deep learning--based stress prediction for bottom-up sla 3d printing
  process.
\newblock {\em The International Journal of Advanced Manufacturing Technology},
  102(5):2555--2569, 2019.

\bibitem{kutz2017deep}
J~Nathan Kutz.
\newblock Deep learning in fluid dynamics.
\newblock {\em Journal of Fluid Mechanics}, 814:1--4, 2017.

\bibitem{guo2016convolutional}
Xiaoxiao Guo, Wei Li, and Francesco Iorio.
\newblock Convolutional neural networks for steady flow approximation.
\newblock In {\em Proceedings of the 22nd ACM SIGKDD international conference
  on knowledge discovery and data mining}, pages 481--490, 2016.

\bibitem{zhang2018application}
Yao Zhang, Woong~Je Sung, and Dimitri~N Mavris.
\newblock Application of convolutional neural network to predict airfoil lift
  coefficient.
\newblock In {\em 2018 AIAA/ASCE/AHS/ASC Structures, Structural Dynamics, and
  Materials Conference}, page 1903, 2018.

\bibitem{lee2020cnn}
Seunghye Lee, Hyunjoo Kim, Qui~X Lieu, and Jaehong Lee.
\newblock Cnn-based image recognition for topology optimization.
\newblock {\em Knowledge-Based Systems}, 198:105887, 2020.

\bibitem{nie2021topologygan}
Zhenguo Nie, Tong Lin, Haoliang Jiang, and Levent~Burak Kara.
\newblock Topologygan: Topology optimization using generative adversarial
  networks based on physical fields over the initial domain.
\newblock {\em Journal of Mechanical Design}, 143(3):031715, 2021.

\bibitem{zhang2019deep}
Yiquan Zhang, Bo~Peng, Xiaoyi Zhou, Cheng Xiang, and Dalei Wang.
\newblock A deep convolutional neural network for topology optimization with
  strong generalization ability.
\newblock {\em arXiv preprint arXiv:1901.07761}, 2019.

\bibitem{banga20183d}
Saurabh Banga, Harsh Gehani, Sanket Bhilare, Sagar Patel, and Levent Kara.
\newblock 3d topology optimization using convolutional neural networks.
\newblock {\em arXiv preprint arXiv:1808.07440}, 2018.

\bibitem{abaqus2011abaqus}
G~Abaqus.
\newblock Abaqus 6.11.
\newblock {\em Dassault Systemes Simulia Corporation, Providence, RI, USA},
  2011.

\bibitem{alet2019graph}
Ferran Alet, Adarsh~Keshav Jeewajee, Maria~Bauza Villalonga, Alberto Rodriguez,
  Tomas Lozano-Perez, and Leslie Kaelbling.
\newblock Graph element networks: adaptive, structured computation and memory.
\newblock In {\em International Conference on Machine Learning}, pages
  212--222. PMLR, 2019.

\bibitem{wang2020kalibre}
Ruihang Wang, Xin Zhou, Linsen Dong, Yonggang Wen, Rui Tan, Li~Chen, Guan Wang,
  and Feng Zeng.
\newblock Kalibre: Knowledge-based neural surrogate model calibration for data
  center digital twins.
\newblock In {\em Proceedings of the 7th ACM International Conference on
  Systems for Energy-Efficient Buildings, Cities, and Transportation}, pages
  200--209, 2020.

\bibitem{belbute2020combining}
Filipe de~Avila Belbute-Peres, Thomas Economon, and Zico Kolter.
\newblock Combining differentiable pde solvers and graph neural networks for
  fluid flow prediction.
\newblock In {\em International Conference on Machine Learning}, pages
  2402--2411. PMLR, 2020.

\bibitem{guo2020semi}
Kai Guo and Markus~J Buehler.
\newblock A semi-supervised approach to architected materials design using
  graph neural networks.
\newblock {\em Extreme Mechanics Letters}, 41:101029, 2020.

\bibitem{ogoke2020graph}
Francis Ogoke, Kazem Meidani, Amirreza Hashemi, and Amir~Barati Farimani.
\newblock Graph convolutional neural networks for body force prediction.
\newblock {\em arXiv preprint arXiv:2012.02232}, 2020.

\bibitem{pfaff2020learning}
Tobias Pfaff, Meire Fortunato, Alvaro Sanchez-Gonzalez, and Peter~W Battaglia.
\newblock Learning mesh-based simulation with graph networks.
\newblock {\em arXiv preprint arXiv:2010.03409}, 2020.

\bibitem{sanchez2020learning}
Alvaro Sanchez-Gonzalez, Jonathan Godwin, Tobias Pfaff, Rex Ying, Jure
  Leskovec, and Peter Battaglia.
\newblock Learning to simulate complex physics with graph networks.
\newblock In {\em International Conference on Machine Learning}, pages
  8459--8468. PMLR, 2020.

\bibitem{gilmer2017neural}
Justin Gilmer, Samuel~S Schoenholz, Patrick~F Riley, Oriol Vinyals, and
  George~E Dahl.
\newblock Neural message passing for quantum chemistry.
\newblock In {\em International Conference on Machine Learning}, pages
  1263--1272. PMLR, 2017.

\bibitem{li2020deepergcn}
Guohao Li, Chenxin Xiong, Ali Thabet, and Bernard Ghanem.
\newblock Deepergcn: All you need to train deeper gcns.
\newblock {\em arXiv preprint arXiv:2006.07739}, 2020.

\bibitem{yang2020revisiting}
Chaoqi Yang, Ruijie Wang, Shuochao Yao, Shengzhong Liu, and Tarek Abdelzaher.
\newblock Revisiting" over-smoothing" in deep gcns.
\newblock {\em arXiv preprint arXiv:2003.13663}, 2020.

\bibitem{koeppe2020intelligent}
Arnd Koeppe, Franz Bamer, and Bernd Markert.
\newblock An intelligent nonlinear meta element for elastoplastic continua:
  deep learning using a new time-distributed residual u-net architecture.
\newblock {\em Computer Methods in Applied Mechanics and Engineering},
  366:113088, 2020.

\bibitem{zhang2018featurenet}
Zhibo Zhang, Prakhar Jaiswal, and Rahul Rai.
\newblock Featurenet: Machining feature recognition based on 3d convolution
  neural network.
\newblock {\em Computer-Aided Design}, 101:12--22, 2018.

\bibitem{peddireddy2021identifying}
Dheeraj Peddireddy, Xingyu Fu, Anirudh Shankar, Haobo Wang, Byung~Gun Joung,
  Vaneet Aggarwal, John~W Sutherland, and Martin Byung-Guk Jun.
\newblock Identifying manufacturability and machining processes using deep 3d
  convolutional networks.
\newblock {\em Journal of Manufacturing Processes}, 64:1336--1348, 2021.

\bibitem{kipf2016semi}
Thomas~N Kipf and Max Welling.
\newblock Semi-supervised classification with graph convolutional networks.
\newblock {\em arXiv preprint arXiv:1609.02907}, 2016.

\bibitem{velivckovic2017graph}
Petar Veli{\v{c}}kovi{\'c}, Guillem Cucurull, Arantxa Casanova, Adriana Romero,
  Pietro Lio, and Yoshua Bengio.
\newblock Graph attention networks.
\newblock {\em arXiv preprint arXiv:1710.10903}, 2017.

\bibitem{gao2019graph}
Hongyang Gao and Shuiwang Ji.
\newblock Graph u-nets.
\newblock In {\em international conference on machine learning}, pages
  2083--2092. PMLR, 2019.

\bibitem{li2019deepgcns}
Guohao Li, Matthias Muller, Ali Thabet, and Bernard Ghanem.
\newblock Deepgcns: Can gcns go as deep as cnns?
\newblock In {\em Proceedings of the IEEE/CVF International Conference on
  Computer Vision}, pages 9267--9276, 2019.

\bibitem{fey2019fast}
Matthias Fey and Jan~Eric Lenssen.
\newblock Fast graph representation learning with pytorch geometric.
\newblock {\em arXiv preprint arXiv:1903.02428}, 2019.

\bibitem{rong2019dropedge}
Yu~Rong, Wenbing Huang, Tingyang Xu, and Junzhou Huang.
\newblock Dropedge: Towards deep graph convolutional networks on node
  classification.
\newblock {\em arXiv preprint arXiv:1907.10903}, 2019.

\bibitem{he2016deep}
Kaiming He, Xiangyu Zhang, Shaoqing Ren, and Jian Sun.
\newblock Deep residual learning for image recognition.
\newblock In {\em Proceedings of the IEEE conference on computer vision and
  pattern recognition}, pages 770--778, 2016.

\bibitem{loukas2019graph}
Andreas Loukas.
\newblock What graph neural networks cannot learn: depth vs width.
\newblock {\em arXiv preprint arXiv:1907.03199}, 2019.

\bibitem{alon2020bottleneck}
Uri Alon and Eran Yahav.
\newblock On the bottleneck of graph neural networks and its practical
  implications.
\newblock {\em arXiv preprint arXiv:2006.05205}, 2020.

\end{thebibliography}

\end{document}